\documentclass[11pt]{article}

\usepackage{amscd,amsmath,amsxtra,amssymb,latexsym,amsfonts,graphics}
\usepackage[usenames,dvipsnames]{color}

\usepackage[small]{caption}

\usepackage{ifthen}
\newboolean{PrivateVersion}
\setboolean{PrivateVersion}{false} 

\newcommand{\private}[1]{ \ifthenelse{\boolean{PrivateVersion}}
{#1}{} }

\newtheorem{Thm}{Theorem}[section]
\newtheorem{Lemma}[Thm]{Lemma}

\newtheorem{Definition}[Thm]{Definition}

\newcommand{\R}{{\mathbb R}}

\newcommand{\C}{{\mathbb C}}
\newcommand{\Z}{{\mathbb Z}}

\renewcommand{\Im}{{\rm Im \ }}
\renewcommand{\Re}{{\rm Re \ }}

\newcommand{\xu}{{\underline{x}}}

\makeatletter
\newif\if@caption@empty
\newcommand{\captionempty}
{%
  \@caption@emptytrue
  \caption{}%
}
\renewcommand{\@makecaption}[2]
{%
  \centering
  \itshape
  \rm #1%
  \if@caption@empty
    \global\@caption@emptyfalse
  \else
    \textbf{:} #2%
  \fi
}
\makeatother

\newcommand{\Ln}{{\rm Ln \ }}

\renewcommand{\Im}{{\rm Im \ }}

\begin{document}

\title{ Shatalov-Sternin's construction of complex WKB solutions and the associated Riemann surface. 
\ifthenelse{\boolean{PrivateVersion}}
{\textcolor{blue}{ \\ !!!PRIVATE VERSION!!!}}{}  }

\author{Alexander GETMANENKO \\
{\footnotesize Institute for the Physics and Mathematics of the Universe,} \\
{\footnotesize The University of Tokyo, 5-1-5 Kashiwanoha, Kashiwa, 277-8568, Japan} \\ 
{\tt \small Alexander.Getmanenko@ipmu.jp}}

\maketitle

\begin{abstract}
We re-examine Shatalov-Sternin's proof of existence of resurgent solutions of a linear ODE. In particular, we take a closer look at the ``Riemann surface" (actually, a two-dimensional complex manifold) whose existence, endless continuability and other properties are claimed by those authors. We present a detailed argument for a part of the ``Riemann surface" most relevant for the exact WKB method. 
\end{abstract}

\section{Introduction.}

\subsubsection*{Resurgent analysis.}

Resurgent analysis is a method of studying hyperasymptotic expansions 
\begin{equation}  \sum_{k,j} e^{-c_k/h} a_{k,j}h^j , \ \ \ h\to 0+  \label{hyperasexpn}  \end{equation}
and those of similar kind by treating such expansions as asymptotics obtained from a Laplace integral 
\begin{equation} \int_\gamma \Phi (s) e^{-s/h} ds, \label{LaplaceTransf} \end{equation} 
where $\Phi$ is a ramified analytic function in the complex domain with a discrete set of singularities and $\gamma$ is an infinite path on the Riemann surface of $\Phi$. The crucial observation is that the terms of (\ref{hyperasexpn}) can be recovered from studying the singularities of $\Phi$, see ~\cite{V83}, ~\cite{E81}, ~\cite{CNP}, ~\cite{DP99}, as well as ~\cite{G} for this author's preferred terminology. 

The methods of resurgent analysis have been used, in particular, to study asymptotics of solutions of linear ODE with a small parameter, especially the Schr\"odinger equation in the semiclassical approximation, see, e.g. ~\cite{DDP97}; this technique is a refinement of what is known as the {\it complex WKB method}. More specifically, consider an equation of the type
\begin{equation} -h^2 \partial_x^2 \varphi(h,x) + V(x) \varphi(h,x) = 0 \label{SchroeEq} \end{equation}
where $x$ ranges over $\C$, $h$ is a small complex asymptotic parameter, and $V(x)$ is an entire function often assumed to be a polynomial.  Under the transformation (\ref{LaplaceTransf}), this equation becomes an equation on an unknown ramified analytic function of two variable $\Phi(s,x)$ of the form
\begin{equation} - \partial_x^2 \Phi(s,x) + \partial_s^2 V(x)\Phi(s,x) \ = \ 0. \label{MainEqu} \end{equation}
The equation \eqref{MainEqu} only needs  to be satisfied modulo functions that are entire with respect to $s$ for every value of $x$ since such functions correspond to zero under a properly (~\cite[Pr\'e I.2]{CNP}) understood Laplace transform \eqref{LaplaceTransf}.  Since the beginnings of resurgent analysis in the early 1980s there has been no real doubt that (\ref{MainEqu}) possesses two linearly independent (in an appropriate sense) solutions that are endlessly analytically continuable with respect to $s$ and satisfy certain growth conditions at infinity.

The manifold on which $\Phi(s,x)$ is defined is usually quite complicated. In the special cases when $V(x)=x$ and $V(x)=x^2$, the function $\Phi(s,x)$ can be written down by an explicit formula and $\varphi(h,x)$ is expressible in terms of Airy or Weber function, see \cite{J94}. For more complicated potentials, say, when $V(x)$ is a generic polynomial of degree $\ge 4$, the function $\Phi(s,x)$ is expected to be defined on a highly transcendental manifold, see \cite{DDP93} and \cite{D92}: if for a fixed $x$ one projects {\it all} singularities on {\it all} sheets of the Riemann surface of $\Phi(s,x)$ to the complex plane of $s$, one expects to obtain an everywhere dense set. Thus, there is no hope that the manifold in question is a universal cover of $\C^2$ minus a discrete family of complex curves.

\begin{figure} \includegraphics{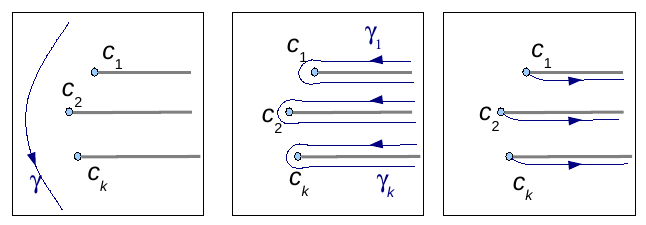} \caption{Deformation of the integration contour and the calculation of the hyperasymptotic expansion of $\varphi(h,x)$} \label{Paper3p3}
\end{figure}

Singularities of $\Phi$ and the precise structure of the manifold on which $\Phi$ is defined are important because they allow us to obtain the hyperasymptotic expansion of $\varphi(h,x)$ for $h\to 0+$ as follows (cf. ~\cite[p.218]{V83}, \cite{CNP}).  Fix $x$ and identify one of the sheets of the Riemann surface of $\Phi(x,s)$ with a complex plane of $s$ minus countably many cuts $c_1+\R_{\ge 0}, c_2+\R_{\ge 0}, ..., c_k+\R_{\ge 0}$ in the positive real direction.    Draw an infinite integration path $\gamma$ in $\C$ to the left of $c_1,c_2,...,c_k,...$, fig.\ref{Paper3p3},left,  so that, at least morally, $\varphi(h,x)=\int_\gamma \Phi(s,x) e^{-s/h} ds$. Using analyticity of $\Phi(s,x)$ and under appropriate conditions on its growth at infinity one can push the integration contour $\gamma$ to the right and rewrite 
$$ \varphi(h,x) \ = \ \sum_k \int_{\gamma_k} \Phi(s,x) e^{-s/h} ds, $$
where infinite integration paths $\gamma_k$ ``hang" on the singularities $c_k$, fig.\ref{Paper3p3},middle. Finally, one deforms each $\gamma_k$ so that both infinite branches lie on different sheets of the Riemann surface right on top of each other, and rewrites 
\begin{equation} \int_{\gamma_k} \Phi(s,x) e^{-s/h} ds \  = \ \int_{[c_k, c_k+\infty)}  (\Delta_{c_k}\Phi(s,x) )e^{-s/h} ds, \label{eq442} \end{equation}
where $\Delta_{c_k}\Phi$ denotes the jump of $\Phi$ across the cut starting at $c_k$. The integrals on the R.H.S. of \eqref{eq442} are taken along semi-infinite real analytic paths similar to those on fig.\ref{Paper3p3},right. The  asymptotic expansions of these integrals can now be calculated using Watson's lemma and combined to a hyperasymptotic expansion \eqref{hyperasexpn}. 

In ~\cite{CNP}, ~\cite{ShSt}, the following point of view is developed. For each fixed $x$, $\Phi(s,x)$ as a function of $s$ is assumed in the beginning to be a holomorphic function on a sectorial neighborhood of infinity $\Omega_0=\{ s\in \C \ : \ \arg s \in (\frac{\pi}{2}-\beta, \frac{3\pi}{2}+\beta;\  |s|>N\}$ for some $\beta>0$ and $N>0$; the contour $\gamma$ appearing in \eqref{LaplaceTransf} is a contour along the boundary of $\Omega_0$. It is then assumed that for a discrete subset $\{c_1,..,c_k,.. \}\subset \C\backslash \Omega_0$, the function $\Phi(s,x)$ has an analytic continuation to the set $\Omega=\C\backslash \bigcup_k (c_k+\R_{\ge 0})$; this $\Omega$ is called {\it the first sheet} of the Riemann surface of $\Phi(s,x)$, and the points $ c_k$, $k=1,2,...$, are called the {\it the first sheet singularities} of $\Phi$. The Riemann surface of $\Phi(s,x)$ for every fixed $x$ is the Riemann surface of the analytic continuation of $\Phi(s,x)$ as a holomorphic function on $\Omega_0$. It is important that in order to  obtain $\Delta_{c_k} \Phi(s,x)$ in \eqref{eq442} as an analytic function of $s$, we define it as  $\Phi(s',x)-\Phi(s'',x)$ where $s',s''$ belong to the different sheets of the Riemann surface of $\Phi(x,s)$ and project to the same point of $s\in\C$; we need therefore an analytic continuation of $\Phi$ beyond the first sheet at least near the cuts $c_k+\R_{\ge 0}$. 
 
While the position of the singularities of $\Phi(s,x)$ is important for the calculation of the asymptotics, there is a good intuition where these singularities are located. Given an initial point $x_0$ and  two ramified analytic functions $f_0(s), f_1(s)$, let $\Phi(s,x)$ solve the Cauchy problem $\Phi(s,x_0)=f_0(s)$, $\frac{\partial}{\partial x}\Phi(s,x_0)=f_1(s)$ for the equation \eqref{MainEqu}. The general philosophy of PDE suggests that the singularities of the initial conditions should propagate along the integral curves of the vector fields $\frac{\partial}{\partial x} \pm \sqrt{V(x)}\frac{\partial}{\partial s}$. Using this intuition, Voros ~\cite{V83} studied the Stokes phenomenon -- appearance and disappearance of singularities from the first sheet of $\Phi(s,x)$ as $x$ varies, and described its consequences (``connection formulas") for the hyperasymptotic expansions of $\varphi(h,x)$.

Since so much relies on the properties of singularities and analytic continuation of $\Phi(s,x)$, proving that \eqref{MainEqu} has an endlessly analytically continuable solution is an important foundational question.  
The present work is a step in this direction.

\subsubsection*{Literature review}

The literature on this subject is extremely vast, so we can hope to at most indicate {\it some} sources which reflect the state of the field and main developments.

The problem of existence and singularities of complex-analytic solutions $\Phi$ of \eqref{MainEqu} appear in numerous classical works, notably ~\cite{Le}, ~\cite{Ha} and their sequels, but the solutions are shown to exist only locally, and the results do not guarantee existence of the analytic continuation of $\Phi(s,x)$ to the values of $x$ far away from an initial point $x_0$ where the Cauchy data are given. 

From ~\cite{DP99} we learned about the existence of a preprint ~\cite{E84} containing a sketch of a construction of endlessly continuable solutions $\Phi$ satisfying \eqref{MainEqu}, but at least according to ~\cite{DP99}, not all details are clear in that sketch. 

Lacking a general statement, one could still work out examples of potentials $V$ for which the function $\Phi$ can be given by a more or less explicit formula and singularities of $\Phi$ are possible to analyze from that explicit representation, see e.g. the easiest examples in  ~\cite{J94} and much more complicated one in the recent article ~\cite{FS}.  

In the terminology of resurgent analysis, the function $\Phi$ appearing in \eqref{MainEqu} is the ``major" of $\varphi$ appearing in \eqref{SchroeEq}. Many authors prefer to take a somewhat different Laplace integral and work with ``minors"; there is a technology of translating statements between the two setups, ~\cite{CNP}.  Working with minors, the authors of ~\cite{DLS93} present a proof that we expect to imply the existence of $\Phi(s,x)$ for values of $s$ on the first sheet minus the cuts and for $x$ confined to a region where no Stokes phenomenon occurs. 

The monograph ~\cite[Ch.3.1]{ShSt} and numerous works by the same authors, e.g. ~\cite{SS93}, ~\cite{SSS97},  contain another approach to the proof of existence of endlessly continuable solutions of \eqref{MainEqu} and of similar equations of higher order. From the parts of the argument that we were able to understand, 
the approach seems very natural and attractive. Discussion  of ~\cite{ShSt}'s proof is the content of this article. 

The topic has remained in the focus of many researchers. It may have been one of the motivations for development of the mould calculus, cf. ~\cite{Sa} and references therein.

Meanwhile the Kyoto school has been working on the idea of transforming the
Schr\"odinger equation with an arbitrary potential $V(x)$ to appropriately chosen canonical models, e.g. Airy,  Weber, Whittaker equations, e.g. ~\cite{AKT91},  ~\cite{KKKT10}; the language of ``minors" is used by these authors.
A breakthrough was announced in the autumn of 2010 by Kamimoto and Koike. Their result is expected to describe the first sheet singularities of $\Phi(s,x)$ as a function of $s$, as long as $x$ is close to a simple zero of a very general potential $V(x)$.  

Not only \eqref{SchroeEq}, but also other similar equations have been studied by means of complex WKB method; respectively, different equations in the Laplace-transformed picture take the place of \eqref{MainEqu}. E.g., higher order ODEs were studied semi-heuristically in ~\cite{AKSST05}, ~\cite{H08}, or rigorously in ~\cite{NNN91}; the  first order difference equations with a small parameter were studied in ~\cite{CG08}. 

In the present article we are re-examining certain details of the Shatalov-Sternin's proof. The idea of the argument presented in ~\cite{ShSt} differs significantly from what the approach of the Kyoto school and from that of other authors. Even in  view of the results announced by Kamimoto and Koike it remains important, for our understanding of the subject as well as for possible extensions and generalizations, to clarify the status of ~\cite{ShSt}'s very natural-looking argument. 

At the time when this version of the article is written, its ideas have been already used in ~\cite{GT}.

\subsubsection*{Contribution of this article.}

In ~\cite{ShSt}, Sternin and Shatalov solve \eqref{MainEqu} by reducing it to an integral equation and obtaining a resolvent. In other words, they try to represent a solution $\Phi(s,x)$ in terms of an infinite series 
\begin{equation} \Phi(s,x) = \sum_{n=0}^\infty \Psi_n(s,x) \label{Jan6} \end{equation}
where $\Psi_n(s,x)$ is, morally, the result of an $n$-fold application of some integro-differential operator to a ``0-th order approximation" $f(s)$. The actual formulas will be recalled in section \ref{SSconSec}.

Having formally obtained an expression \eqref{Jan6}, Sternin and Shatalov set out to prove that a) all functions $\Psi_n(s,x)$ are defined on the same endlessly continuable manifold of complex dimension two (which is still called a ``Riemann surface"), and that b) the series converges on compact sets of this ``Riemann surface".

In ~\cite[Prop.3.1, pp.204-207]{ShSt}, the construction of the ``Riemann surface" takes only three pages and is presented very intuitively; however, once we wanted to make a precise sense of how exactly the ``Riemann surface" is described and how exactly all functions $\Psi_n$ can be analytically continued to it by which specific deformations of integration contours, we found ourselves dealing with a rather complex situation. For now we restrict ourselves to constructing an open piece ${\cal S}$ of the ``Riemann surface". As a bit of an oversimplification, let us say that over each point $x$ in an appropriate region of the complex plane, the fiber of ${\cal S}$ consists of a the first sheet (i.e. the complex plane with finitely many cuts)  and small ``flaps"  attached on the sides along each cut, see section \ref{FSS}. Spelling out all the details is the content and the contribution of this work.

Thus, the statement and the proof of the following theorem are intended to make precise some things which we could not find in ~\cite{ShSt}. 

\begin{Thm} \label{ContribTh} For $V(x)$ satisfying assumptions of section \ref{NotationSec}, the countably many functions  \eqref{vNsimpl} possess an analytic continuation to the 2-dimensional complex manifold ${\cal S}$ defined in section \ref{StructureRS}. \end{Thm} 

A word of caution: The functions $\Psi_n(s,x)$ appearing in \eqref{Jan6} are more complicated than functions \eqref{vNsimpl}, but it will be obvious that the theorem implies that $\Psi_n$ also analytically continue to ${\cal S}$.

Here is what remains outside the scope of this article. The series \eqref{Jan6} is very likely to converge uniformly on compact subsets of ${\cal S}$. Unfortunately, in ~\cite[(3.14)]{ShSt} the derivative in the integrand of ~\eqref{Rj} is missing, and those authors end up proving convergence of a wrong and much better behaving series. A more delicate study of convergence will need to be performed in the future. The current paper makes the question more well-defined: before we study convergence of the series \eqref{Jan6} at a point $(s,x)$ of ${\cal S}$, we need to know first {\it how exactly} the functions $\Psi_n$ are analytically continued to the point $(s,x)$.
If the convergence is shown, that will provide an alternative both to the approach announced by Kamimoto and Koike and to the method of ~\cite[GT]. 

We will finish this introduction by indicating what is involved in the proof of theorem \ref{ContribTh}. As the ramified analytic functions  \eqref{vNsimpl} of variables $(s,x)$ are iterations of two integro-differential operators $R_1$ and $R_2$, in order to analytically continue these functions to a point $(s,x)$ we need to appropriately define two integration paths (one for $R_1$ and one for $R_2$) leading from $(s_0,x_0)$ to $(s,x)$; here $x_0$ is some fixed initial point and $s_0$ depends on $s$ and $x$. First we treat the case when $x$ is in the same Stokes region as $x_0$, and then describe in the lemmas of section \ref{modelcases} a method that allows us to draw the integration paths for $x$ belonging to further and further Stokes regions. As we take $x$ in Stokes regions further and further away from $x_0$, there appear more and more obstacles to drawing an integration path from $(s_0,x_0)$ to $(s,x)$; points $(s,x)$ that cannot be reached by an integration path give rise exactly to
 the  singularities of ${\cal S}$
 predicted by Voros. 

Unfortunately lemmas of the section \ref{modelcases} do not define an clear-cut inductive procedure, as we have not yet systematized many little irregular combinatorial details occupying section \ref{ApplyManyLemmas}. Still, by referring to section \ref{modelcases} we are able to construct all the paths of analytic continuation relevant for the proof of theorem \ref{ContribTh}. We believe that more complicated potentials $V(x)$ can be treated similarly by using section \ref{modelcases}; see also a remark on the combinatorial complexity of this problem on p.\pageref{Combina}.

\section{Shatalov-Sternin's construction.} \label{SSconSec}

The purpose of this section is to review the content of  ~\cite[pp.198-204]{ShSt} in the special case of the one-dimensional Schr\"odinger equation
\begin{equation}  [-h^2 \partial^2_x + V(x)]\varphi(h,x) = 0, \label{SchroeEq1} \end{equation}
where the variable $x$ takes values in $\C$ and $V(x)$ is an entire function. 

To describe the Laplace-transformed version of \eqref{SchroeEq1}, consider the following operation on the equivalence classes of germs of analytic functions at a point $(s_0,x_0)\in \C^2$ modulo functions entire with respect to $s$ for every $x$: 
$${\hat h}\Phi(s,x) \ := \  \partial_s^{-1} \Phi(s,x) = \int^s_{s_*(x)} \Phi(s',x)ds',$$
where the starting point of the integration $s_*(x)$ may depend on $x$ and  changing $s_*(x)$ will change the result by a function depending only on $x$. 

In this notation, the Laplace transform  \eqref{LaplaceTransf} turns \eqref{SchroeEq1} into 
\begin{equation} -{\hat h}^2 \partial_x^2 \Phi(s,x) + V(x)\Phi(s,x) \ = \ 0 \label{MainEqu1} \end{equation}
which has to be satisfied modulo functions that are entire with respect to $s$ for every $x$. We would like to find solutions $\Phi$ of \eqref{MainEqu1} that are holomorphic functions on a complex two-dimensional manifold ${\cal S}$ endowed with a locally biholomorphic projection $\Pi$ to $\C^2$ with coordinates $(s,x)$. We would also like, for every $\xu\in\C$, the connected components of $\Pi^{-1}(\{ (s,\xu):s\in \C\})$ to be endlessly continuable Riemann surfaces in the sense of resurgent analysis, e.g., ~\cite[R\'es I]{CNP}. In fact, ~\cite{ShSt} use the concept of a ``ramified analytic function" of several complex variables; we will replace it by a clearer notion of ``a germ of an analytic function" except in philosophical statements.

The Cauchy-Kowalewskaya theorem, e.g.~\cite[Th.3.1.1]{Sch}, or the related results of ~\cite{Le} and ~\cite{Ha}, for this equation fall far short of the statement that we need. Indeed, for the equation \eqref{MainEqu1} with an initial condition, say, $\Phi(s,x_0)=\frac{1}{2\pi i s}$, $\frac{\partial}{\partial x}\Phi(s, x_0)=0$ (corresponding to $\varphi(h,x_0)=1$, $\frac{\partial}{\partial x}\varphi(h,x_0)=0$) one would only get existence of solution $\Phi(s,x)$ in a small polydisc centered at $(s_0,x_0)$ for $s_0\ne 0$, and the size of that polydisc is hard to increase. Therefore a more explicit construction of $\Phi$ is proposed. 

Fix a point $x_0$ such that $V(x_0)\ne 0$ and a determination $p(x)$ of $\sqrt{V(x)}$ in a neighborhood of $x_0$. Let $p_1(x)=-p_2(x)=p(x)$; let further  $S_j(x) = \int_{x_0}^x p_j(y) dy$, $j=1,2$, and $S(x)=S_1(x)$. 
In this notation, the operator $-\hat h^2 \partial_x^2 +V(x)$ on the L.H.S. of \eqref{MainEqu1} can be rewritten as 
\begin{equation} \left(p^2(x) [-\frac{1}{p(x)} \hat h\partial_x - 1] - {\hat h}p'(x) \right)\left(\frac{1}{p(x)}{\hat h}\partial_x - 1\right) - {\hat h}p'(x). \label{splitOrd1} \end{equation} 

We will be able to make use of this representation once we are able to invert the operators $\pm \frac{1}{p(x)} \hat h \partial_x -1$. Namely, consider an equation
\begin{equation} [\frac{1}{p_j(x)} \hat h \partial_x -1] u(s,x) = b(s,x), \label{e308} \end{equation}
as an equation of germs at $(s_0,x_0)\in \C^2$ of analytic functions of $(s,x)$  modulo functions depending only on $x$. Then \eqref{e308} is satisfied by 
$$ u(s,x) \ = \ R_j b(s,x) + f(s+S_j(x)), $$
where $f(s)$ is any germ of an analytic function near $s_0$ and the operator $R_j$  is defined by the formula
\begin{equation} (R_j G)(s,x) = \int_{x_0}^x (D_1G)(s+S_j(x)-S_j(y),y) p_j(y) dy.  \label{Rj} \end{equation}
where $D_1$ stands for the derivative of the function with respect to the first argument. We consider $R_j$ as acting on germs of analytic functions $G(s,x)$ at a point $(s_0,x_0)$. In ~\cite{ShSt} this derivative is missing.

Let us start looking for a solution \eqref{MainEqu1} in the form
\begin{equation} \Phi(s,x) = R_1 \Phi_1(s,x) + f_1(s+S_1(x)).  \label{e327} \end{equation}
Substituting \eqref{e327} into \eqref{MainEqu1} and using the expression \eqref{splitOrd1}, we have
\begin{equation} \left\{\left(p^2(x) [-\frac{1}{p(x)} {\hat h}\partial_x - 1] - {\hat h}p'(x) \right) - {\hat h}p'(x) R_1\right\} \Phi_1 \ = \
 - {\hat h} p'(x) f_1(s+S_1(x)).  \label{e330} \end{equation}
Looking for a solution of \eqref{e330} in the form
$$ \Phi_1(s,x) = R_2\Phi_2(s,x) + f_2(s+S_2(x)), $$
we obtain
{ \small
$$ \left[ 1 - {\hat h}\frac{p'(x)}{p^2(x)}\{ R_2  + R_1 R_2 \} \right] \Phi_2  \ = \  - {\hat h} \frac{p'(x)}{p^2(x)}  \{ (1+R_1)  f_2(s+S_2(x)) + f_1(s+S_1(x)) \}.  $$
}

Formally, the last equation has a solution
\begin{equation} \Phi_2(s,x) \ = \  \sum_{j=0}^{\infty} (-1)^j {\hat h}^{j+1} [ (-\frac{p'(x)}{p^2(x)})(1+R_1) R_2 ]^j g_0(s,x) , 
\label{vNseries} \end{equation} 
where 
$$ g_0(s,x) \ = \ - {\hat h} \frac{p'(x)}{p^2(x)}  \{ (1+R_1)  f_2(s+S_2(x)) + f_1(s+S_1(x)) \}. $$
 On the R.H.S. of \eqref{vNseries} we see an infinite series of germs of analytic functions; only its partial sums are mathematically well-defined at this stage.

Assume that we are able to prove that the series on the right hand side of \eqref{vNseries} converges both for the choice a) $f_1(s)=\Ln s$, $f_2=0$, and for the choice b) $f_1=0$, $f_2(s)=\Ln s$, and in both cases defines analytic functions $\Phi_2(s,x)$ and $\Phi(s,x)$ on a sufficiently large complex two-dimensional manifold. Then we can perform a Laplace integral as in \eqref{eq442}; as a result, we expect to obtain two formal WKB solutions of \eqref{SchroeEq1} for $x$ in a neighborhood of $x_0$, namely $A_+(h,x) e^{S(x)/h} + A_-(h,x)e^{-S(x)/h}$ for the choice a), and $B_+(h,x) e^{S(x)/h} + B_-(h,x)e^{-S(x)/h}$ for the choice b). Here $A_\pm(h,x), B_\pm(h,x)$ are expected to be formal (actually, Gevrey) power series in $h$ with $x$-dependent coefficients. We expect further that the two vectors $[A_+(h,x_0), A_-(h,x_0)]$ and $[B_+(h,x_0), B_-(h,x_0)]$ in $\C[[h]]^2$ will be linearly independent over $\C[[h]]$, thus yielding two linearly independent resurgent solutions of \eqref{SchroeEq1} in every reasonable definition of this notion.

The first task is therefore to construct a ``Riemann surface'' -- a two dimensional complex manifold on which all summands in the R.H.S. of \eqref{vNseries} are defined for the choices a) and b) from the previous paragraph. It is easy to see that an equivalent question is to construct a ``Riemann surface'' on which all functions 
\begin{equation} R_{j_k}...R_{j_2}R_{j_1} f(s,x) , \ \ \ j_i=1,2, \ \ \ k\ge 0  \label{vNsimpl} \end{equation}
are defined for $f(s,x)=\Ln(s\pm S(x))$. 

This is the question we are dealing with in this article. The second task would be to show that the infinite series converges on this ``Riemann surface". Unfortunately, a derivative in the integrand is missing in ~\cite{ShSt}'s definition of operators $R_j$ and we cannot suggest an easy way to repair their convergence argument, but hope to give (or read!) an alternative proof elsewhere.


\section{Analytic continuation and integration paths} \label{ACIP}

In sections \ref{Notation} , \ref{StructureRS} we are going to precisely define the ``Riemann surface" ${\cal S}$ to which we will then be able to analytically continue the functions \eqref{vNsimpl}. The section \ref{RedToPaths} exposes the main idea of this article; its content will make precise sense after reading sections  \ref{Notation} and \ref{StructureRS}.  For now we will think of ${\cal S}$ as some complex two-dimensional manifold with a locally biholomorphic projection ${\cal S}\to \C_{s}\times\tilde {\cal O}$, where $\tilde {\cal O}$ is a complex one-dimensional manifold with a locally biholomorphic projection to $\C_x$, and $\C_s,\C_x$ denote the complex planes of the variables $s,x$, respectively.  We will freely use $(s,x)$ as local coordinates on ${\cal S}$.

\subsection{Reduction of the problem to construction of the integration paths.} \label{RedToPaths}

Recall that we denote  $p_1(x)=-p_2(x)=p(x)$,  $S_j(x) = \int_{x_0}^x p_j(y) dy$, $j=1,2$, and the operators $R_j$ were defined by \eqref{Rj} as operators acting on germs of analytic functions. 

It will be obvious from the construction of ${\cal S}$ that the functions $\Ln (s\pm S(x))$ have analytic continuations to ${\mathcal S}$. Existence of analytic continuation of all terms of \eqref{vNsimpl} to ${\mathcal S}$ will follow by induction from the following

\begin{Thm} If $G(s,x)$ is defined and analytic on ${\mathcal S}$, then $R_jG$, $j=1,2$ have analytic continuations to ${\mathcal S}$. \label{RjGmainTh}
\end{Thm}


A detailed proof of this theorem will be given in section \ref{PfMresSec}. In this section \ref{RedToPaths} we will introduce some of the terminology used in the proof; in section \ref{AppearsStokes} we will informally  explain the idea on which the proof is based.

If $G(s,x)$ were an analytic function on the whole $\C\times \tilde{\cal O}$, we could define $(R_jG)(s,\xu)$ by the formula 
\begin{equation} (R_j G)(s,\xu) = \int_{x_0}^\xu (D_1G)(s+S_j(\xu)-S_j(y),y) p_j(y) dy\label{Rj1} \end{equation}
 where the integral is taken along {\it any} path from $x_0$ to $\xu$ in $\tilde{\mathcal O}$. Since, however, $G(s,\xu)$ is defined on a complicated manifold ${\cal S}$, we need to find for each $(s,\xu)\in{\cal S}$ a path $y(t)$ in $\tilde{\mathcal O}$ from $x_0$ to $\xu$ satisfying the following 

{\bf Definition.} We say that a path $y(t)$ in $\tilde {\cal O}$  {\it can be  lifted to ${\cal S}$ parallel to $-S_j$ with endpoint $(s,\xu)$} if  $(s+S_j(\xu)-S_j(y(t)),y(t))$ defines a path in ${\mathcal S}$.

Intuitively, this condition means that the point $(s+S_j(\xu)-S_j(y(t)),y(t))$  does not ``leave" ${\cal S}$ and does not hit any of its singularities.

We will call such a $y(t)$  {\it an integration path for $(s,\xu)$ and $R_j$} and draw it in green on our figures; let us stress that the choice of a path $y(t)$ depends on $s$ in the fiber ${\cal S}_\xu$ of ${\cal S}$ over $\xu$. If the integration paths $y(t)$ continuously depend on $(s,\xu)$, using them in \eqref{Rj1} yields an analytic function $R_jG(s,x)$; construction of  $R_jG$ from $G$ is thus reduced to finding a family of  integration paths continuously depending on $(s,\xu)$. \label{LiftingDefined}

\private{\textcolor{blue}{It seems that this paragraph is no longer needed.} Note that if $\xu$ is contained in a compact subset $K\subset\tilde{\mathcal O}$ and $x_0\in K$, then there is a constant $N_K\in \R$ such that any integration path $y(t)$ from $x_0$ to $\xu$ contained in $K$ will satisfy the desired property \textcolor{blue}{which property?} for any $s$ with $\Re s<N_K$. \textcolor{blue}{HEAVY!} }


In Section \ref{StructureRS} we will describe the fibers  ${\cal S}_x$ of ${\cal S}$ over every $x\in \tilde {\cal O}$; for each $x\in \tilde{\cal O}$, we will define in \eqref{ListOfSing} a list of {\it singularities} in ${\cal S}$ each of which will be of the form $s=S_j(x)+c$, $j=1,2$, $c\in \C$, for appropriate constants $c$.

For $U\subset{\mathcal S}_\xu$ let us try to construct an integration path $y(t)$ which for any endpoint $(s,\xu)$, $s\in U$, can be lifted to ${\cal S}$ parallel to $-S_j$. We want to make sure that $s+S_j(\xu)-S_j(y(t))$ avoids the singularity $S_j(y(t))+c$, i.e. we want the equality 
$$ 2S_j(y(t)) \ = \ s+S_j(\xu)-c $$
to hold for no point $y(t)$ along the integration path and for no point $s\in U$. That is to say, we want the integration path $y(t)$ to avoid the set 
\begin{equation} S_j^{-1}\left( \frac{U+S_j(\xu)-c}{2} \right)\subset \tilde{\mathcal O}. \label{SjInverse} \end{equation}
 We will need to carefully keep track of the appropriate branches of the functions involved in this expression.

On our figures we will draw the boundary of $U$ in red and the boundaries of the sets of type  $V=S_j^{-1}\left( \frac{U+S_j(\xu)-c}{2} \right)$ in purple.

It will turn out a posteriori that the condition that $(s+S_j(\xu)-S_j(y(t)),y(t))$ does not coincide with any of the singularities of ${\cal S}_{y(t)}$  is enough to guide us through the choice of the integration paths $y(t)$ for the point $(s,\xu)$. Once a choice of an integration path $y(t)$ is proposed, it is an extra logical step to check that its lifting parallel to $-S_j$ stays within ${\cal S}$; this however will always be obvious by inspection and not mentioned explicitly.

When constructing integration paths $y(t)$ for $R_j$, we found it convenient to construct the parallel transport of the set $U\in {\mathcal S}_{\xu}$ by defining $U(y(t))=U+S_j(\xu)-S_j(y(t))$ (in terms of the projection of ${\cal S}_{y(t)}$ to the complex $s$-plane). Then, as $t$ varies, the set $U(y(t))$ and the singularities  of type $-S_j(y(t))+const$ move with respect to the $s$-coordinate parallel to each other, and differently from the singularities of type $S_j(y)+const$.  For this reason, we we will introduce the following terminology. 

\begin{Definition} \label{StatMovingDef} {\rm For a given choice of $j\in \{1,2\}$, we call singularities of type $-S_j(y)+const$ {\it stationary singularities} and the singularities of type $S_j(y)+const$ {\it moving singularities}. When the index $j$ changes, the roles of moving and stationary singularities reverse. }
\end{Definition}



\subsection{Appearance of the Stokes curves in the construction of the integration paths.} \label{AppearsStokes}

This subsection \ref{AppearsStokes} is written informally and included for illustrative purposes only; the precise argument in the rest of the paper does not logically depend on it.

Recall that $\C_s,\C_x$ denote the complex planes of the variables $s,x$, respectively.

As the functions $S_j$ enter into the definition of the operators $R_j$, it is natural to choose $\tilde {\cal O}$ from the introductory paragraph of section \ref{ACIP} in such a way that $\tilde {\cal O}\to \C_x$ factors through the universal cover of $\C_x \backslash V^{-1}(0)$ with the base point $x_0$.

Let us discuss the construction of $R_1G_0$ and $R_2G_0$ for the function $G_0(s,x)=\Ln(s+S(x))$ (compare to \eqref{vNsimpl}). This function $G_0(s,x)$ is naturally defined on a ``Riemann surface" ${\cal S}_0$ whose fiber over any $x\in \tilde{\cal O}$ is a universal cover of $\C_s\backslash \{ -S(x)\}$. All iterations $R^n_1G_0$, $n\ge 1$, are also defined on ${\cal S}_0$: arbitrary paths in $\tilde{\cal O}$ can be chosen as integration paths for defining $R_1G_0$ {\it for this specific function $G_0$}. On the contrary, the ``Riemann surface" ${\cal S}_2$ of $R_2G_0$ necessarily has (at least) an additional singularity at $s=S(x)$: for any integration path $y(t)$ from $x_0$ to $\xu$ for $R_2$ and $(s=S(\xu),\xu)$, the integrand of \eqref{Rj1} is singular for $y=x_0$ because $S(x_0)=0$ and  $G_0$ has a singularity at $(x_0,0)$. 

Thus, the common ``Riemann surface" ${\cal S}$ of all the functions \eqref{vNsimpl} necessarily has singularities at $s=S(x)$ and $s=-S(x)$ on its first sheet.

\begin{figure} \includegraphics{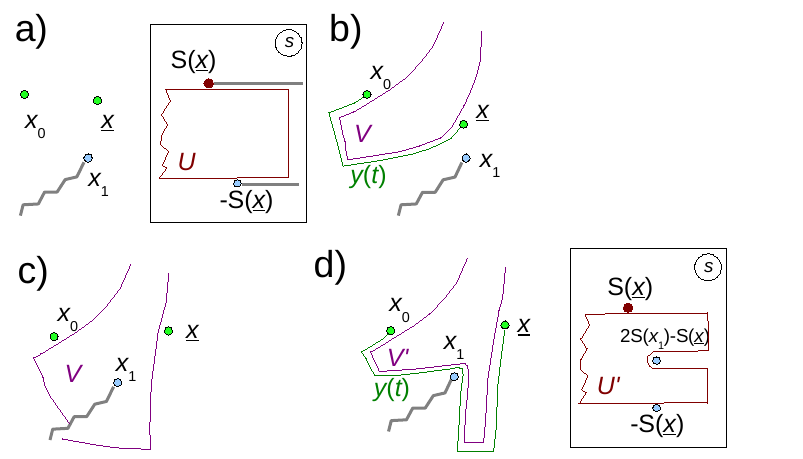}  \caption{Section \ref{AppearsStokes}. a) Projections of $x_0$ and $\xu$ to the complex plane of $x$; the set $U$ in the fiber of ${\cal S}_2$ over $\xu$; b) The set $V_\C$ and the integration path $y(t)$ in the case $\Im S(\xu) < \Im S(x_1)$; c) The set $V_\C$ in the case $\Im S(\xu) > \Im S(x_1)$; d) In the situation of c), the set $V'$ and the integration path $y(t)$ in the complex plane of $x$ and the set $U'$ in the fiber of ${\cal S}_2$ over $xu$.
The branch cut starting from $x_1$ reminds us that the function $S(x)$ has a branch point at $x_1$.} \label{Paper3p44} \end{figure}

 Suppose $x_1\in \C$ is a zero of $V(x)$, and all other zeros of $V(x)$ are far enough from $x_0$ so as not to affect our reasoning here; let $\Im S(x_1)>0$. Take a point $\xu\in {\cal O}$ such that $\Im S(\xu)>0$. Assume that the function $G_2(s,x)=(R_2G_0)(s,x)$ is defined on the Riemann surface ${\cal S}_2$ with singularities at $s=\pm S(x)$  on the first sheet. Let us study whether $R_1G_2(s,x)$ can be analytically continued to the set $U=\{ s\in \C : \Im [-S(\xu)]<\Im s <\Im S(\xu), \ \Re s < N \}$ identified with a subset in the fiber ${\cal S}_{2}$ over $\xu$, where $N\in \R$ is a large positive number, fig. \ref{Paper3p44},a). 

The reasoning of \eqref{SjInverse} with $S(x)$ playing the role of $S_j(x)+c$ leads us to considering the set $V= S^{-1}(\frac{U+S(\xu)}{2})\subset \tilde {\cal O}$; let $V_{\C}$ denote the subset of $\C$ given by the same formula. If $\Im S(\xu) < \Im S(x_1)$, fig. \ref{Paper3p44},b),  then it is possible to draw an integration path $y(t)$  in $\tilde{\cal O}$ from $x_0$ and $\xu$. In a careful treatment, one sees that the set $U$ can indeed be transported along $y(t)$ parallel to $-S_1(x)$. 

If, on the contrary, $\Im S(\xu) > \Im S(x_1)$ and $N$ is large enough, $x_0$ and $\xu$ belong to different connected component of $\tilde {\cal O}\backslash V_\C$, fig. \ref{Paper3p44},c), and an integration path $y(t)$ cannot be drawn. The situation is however remedied if instead of $U$ one considers a smaller subset $U'=U\backslash B_\varepsilon(2S(x_1)-S(\xu)+\R_{\ge 0})$ where $B_\varepsilon$ denotes an $\varepsilon$-neighborhood of a subset of the complex plane of $s$, for $\varepsilon>0$ small enough; the set $U'$ and the corresponding set $V'=S^{-1}(\frac{U'+S(\xu)}{2})\subset \tilde{\cal O}$ and a possible path $y(t)$ are shown on fig.\ref{Paper3p44},d). This strongly suggests that for $\Im S(x)> \Im S(x_1)$ the first sheet of the Riemann surface of $R_1G_2$ contains a singularity at $s=2S(x_1)-S(\xu)$. We immediately recognize the curve $\Im S(x)=\Im S(x_1)$ as the Stokes curve and appearance of the new singularity as the Stokes phenomenon known, e.g., from ~\cite{V83}.


\section{Notation and terminology} \label{Notation}

\subsection{The potential under consideration. Stokes curves and Stokes regions.} \label{NotationSec}

 In ~\cite{DDP97}, Schr\"odinger operators with many different potentials $V(x)$ have been studied using the exact WKB method. Every example of $V(x)$ gives rise to its own pattern of turning points and Stokes curves. In this article we will confine our attention to a piece of the complex plane of $x$ surrounding a commonly occurring piece of the Stokes pattern: two simple turning points and the total of six unbounded Stokes curves starting from them. In this section, we will formally describe such a situation. Our analysis in further chapters suggests that one should look for solutions of \eqref{MainEqu} on the universal cover of the complex plane with the turning points removed; we will define an appropriate piece of this universal cover. Finally, we will formulate the requirement that various Stokes curves are not too close to one another compared to some number $\delta>0$.

From this section and in the next we present a careful description of the particular case of a potential and of the ``Riemann surface" ${\cal S}$ that we are going to study.

On all the figures below  thick gray lines indicate branch cuts of respective Riemann surfaces.

\paragraph*{Assumptions on $V(x)$.} 
In this section we will describe a typical potential well in a potential $V(x)$ and two simple turning points $x_1$ and $x_2$. We will draw the total of six Stokes curves emanating from $x_1$ and $x_2$ and consider  their neighborhood in the complex plane of $x$. It is in this neighborhood that the summands of \eqref{vNsimpl} will be constructed. Let us now say this more formally.

Let $V(x)$ be a function analytic on the closure of a domain ${\mathcal O}_0\subset \C$ which is simply connected and such that $\C\backslash{\mathcal O}_0$ has four connected components $B_1$,..., $B_4$ numbered in a clockwise order. Let $V(x)$ have exactly two distinct zeros in ${\mathcal O}_0$ at points $x_1$ and $x_2$, and both zeros are simple. 

\begin{figure}
\includegraphics{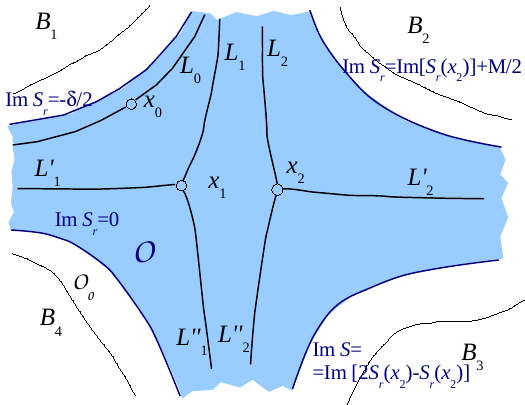} \caption{Domain ${\mathcal O}$. } \label{Paper3p2} \end{figure}

For $j=1,2$, let  $L_j$, $L'_j$, $L''_j$ be curves given by the equation $\Im \int_{x_j}^x \sqrt{V(y)}dy=0$ for $x$ on any of these curves (this definition does not depend on the choice of the square root). Suppose all these curves go off to infinity inside ${\mathcal O}_0$:  $L_1, L_2$  between $B_1$ and $B_2$, $L'_2$ between $B_2$ and $B_3$, $L''_2, L''_1$ between $B_3$ and $B_4$ and $L'_1$ between $B_4$ and $B_1$, fig. \ref{Paper3p2}.

The curves $L_j$, $L'_j$, $L''_j$, $j=1,2$, as well as their preimages on the universal cover of ${\mathcal O}_0\backslash\{x_1,x_2\}$ are called {\it Stokes curves}. 

Fix a determination $p(x)$ of $\sqrt{V(x)}$ and denote by $p_r(x)$ its restriction to ${\mathcal O}_0\backslash(L'_1\cup L'_2)$,; assume that $\Re \int^x p_r(y)dy$ decreases along $L_j$ and increases along $L'_j$ in the direction away from $x_j$, $j=1,2$.

Let $x_0$ be a point in the part of ${\mathcal O}_0$ bounded by $L_1$, $L'_1$, and $\partial B_1$. The function $S(x) = \int_{x_0}^x p(y)dy$ is well-defined on the universal cover of ${\mathcal O}_0\backslash\{x_1,x_2\}$ (with base point $x_0$; denote by $S_r$ its restriction to ${\mathcal O}_0\backslash (L'_1\cup L'_2)$ continuously extended to the points $x_1$ and $x_2$.

Fix a number $\delta>0$. Assume there is a constant $M>0$ such that: $\Im S_r(x)<-\delta/2$ on $\partial B_1$, $\Im S_r(x)-S_r(x_2) > M/2$ on $\partial B_2$, $\Im S_r(x)<0$  for $x\in\partial B_4$, and that $\Im [S_r(x_1)+S_r(x_2)-S_r(x)] <0$ on $\partial B_3$.  
Assume that
\begin{equation}
\delta \ < \ \frac{1}{3} \min\{ \Im S(x_1);  \  \frac{1}{3}\Im[S(x_2)-S(x_1)], \ M \}.
\label{DeltaIsSmall}
\end{equation}

The assumptions on $V(x)$ have now been listed completely.

\paragraph*{Further notation.} Now let us consider a subdomain ${\mathcal O}$ of ${\mathcal O}_0$ (figure \ref{Paper3p2}) bounded by the curve $\Im S_r(x)=-\delta/2$ in the ``quadrant'' defined by $L_1$ and $L'_1$, bounded by $\Im S_r(x)=\Im S_r(x_2) + M/2$ in the ``quadrant'' defined by $L_2$ and $L'_2$, bounded by $\Im S_r(x)=\Im [2S_r(x_2)-S_r(x_1)]$ in the ``quadrant'' defined by $L'_2$ and $L''_2$, bounded by $\Im S_r(x)=0$ in the ``quadrant'' defined by $L'_1$ and $L''_1$. 

Denote by $L_0$ the curve $\Im S(x)=0$ passing through $x_0$.

\begin{figure}
\includegraphics{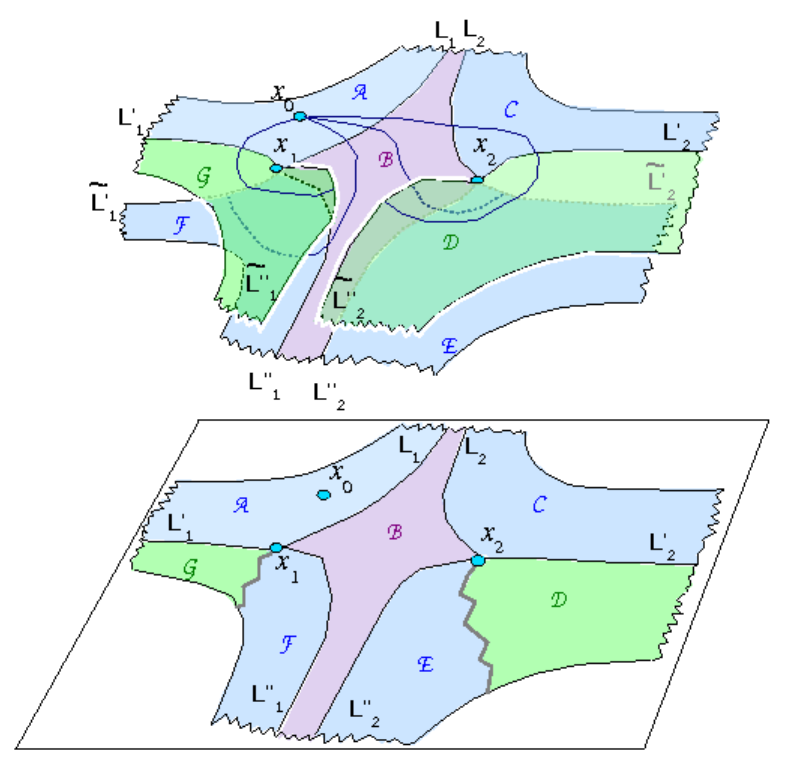} \caption{The Stokes regions ${\cal A}$,${\cal B}$,${\cal C}$,${\cal D}$,${\cal E}$,${\cal F}$,${\cal G}$ as subsets of $\tilde{\mathcal O}$ and their projections on the complex plane of $x$. Also shown are paths from $x_0$ to the curves ${\tilde L}'_1$,${\tilde L}''_1$, ${\tilde L}'_2$,${\tilde L}''_2$. } \label{Paper3p13} \end{figure}

In the universal cover of ${\mathcal O}\backslash\{x_1,x_2\}$ (with base point $x_0$) consider preimages ${\tilde L}'_1$, ${\tilde L}''_1$, ${\tilde L}'_2$, ${\tilde L}''_2$ of Stokes curves $L'_1, L''_1, L'_2, L''_2$ lying on further sheets, fig.\ref{Paper3p13}. Consider the open subset $\tilde{\mathcal O}$ of the universal cover of ${\mathcal O}\backslash\{x_1,x_2\}$ bounded by ${\tilde L}'_1$, ${\tilde L}''_1$, ${\tilde L}'_2$, ${\tilde L}''_2$. The curves ${\tilde L}'_1$, ${\tilde L}''_1$, ${\tilde L}'_2$, ${\tilde L}''_2$ will be called {\it external Stokes curves}.


Denote the sets of curves ${\cal L}_{iS}= \{ L_j, L'_j, L''_j \}_{j=1,2} $, ${\cal L}_{eS}={\tilde L}'_1, {\tilde L}''_1, {\tilde L}'_2, {\tilde L}''_2\}$, ${\cal L}_S ={\cal L}_{iS}\cup {\cal L}_{eS}$, ${\cal L} = \{ L_0 \}\cup {\cal L}_S$.

 The open subsets of $\tilde{\cal O}$ bounded by curves in ${\cal L}_{iS}$ are called {\it Stokes regions} and denoted ${\cal A}$, ${\cal B}$, ${\cal C}$,${\cal D}$, ${\cal E}$, ${\cal F}$, ${\cal G}$ as on fig.\ref{Paper3p13}.  As the Stokes regions are subsets of $\tilde {\cal O}$, a closure of a Stokes region contains the bounding Stokes curves but does not contain any turning points.  We will think of $\tilde {\mathcal O}$ as of a subset of the universal cover of ${\cal O}_0\backslash \{ x_1,x_2\}$. 
 It is for $x$ in this set $\tilde{\mathcal O}$ that we will be discussing the construction of the solution for the equation (\ref{MainEqu1}).  

As usual, the {\it canonical distance} between points $\xu_1$ and $\xu_2$ of $\tilde{\mathcal O}$  is $\inf \int_\pi |p(y)dy| $ where the infimum is taken over all paths $\pi$ in $\tilde{\mathcal O}$ connecting $\xu_1$ and $\xu_2$.

If $\delta>0$, $L\in {\cal L}_S$ let
$$ \begin{array}{lll} {\cal U}_{L,\delta} = \{ x\in \tilde {\cal O} \ : & 1) \ can.dist.(x,L)<\delta ; &  
\\ & 2) \ \Re S(x)\ge \Re S(x_t), \ \text{resp.} \Re S(x)\le \Re S(x_t) & \}, \end{array} $$
where $x_t\in \{ x_1,x_2\}$  is the starting point of $L$ and the sign in the second condition is $\ge$ if $\Re S(x)$ is increasing along $L$ away from $x_*$, and $\le$ otherwise. 
Also let
$$ {\cal U}_{L_0,\delta} = \{ x\in \tilde {\cal O} \ : \ can.dist.(x,L_0)<\delta \}. $$


For a Stokes region ${\cal X}$, consider the part of ${\cal X}$ that is separated away from the internal Stokes curves: 
$${\cal X}_{int}= {\cal X}\backslash \bigcup_{L\in {\cal L}_{iS}; \ L\subset \overline{\cal X}} {\cal U}_{L,\delta/2}.$$
In the collection of sets 
\begin{equation}{\mathbb S}=\{ {\cal A}_{int},...,{\cal G}_{int}\}\cup \{ {\cal U}_{L,\delta/2} \}_{L\in {\cal L}_{iS}} \label{tiling} \end{equation}
 there is a partial order: we say that$A\in {\mathbb S}$ {\it is closer to $x_0$}, or {\it comes earlier} than $B\in {\mathbb S}$ if any path in $\tilde{\cal O}$ from $x_0$ to a point in $B$ has to pass through $A$.


\subsection{Flaps, strips, and slots.} \label{FSS}

 In the further chapters we will start from a complex plane with a few cuts and enlarge it a little by attaching flaps along the cuts. Here we are formally introducing the notation to express this idea.

\paragraph*{Flaps.} Let $k\in \Z_{\ge 0}$, $s_1,...,s_k \in \C$ and $\Im s_{j+1} < \Im s_{j}$. 
 Consider the set $U_0$ obtained from $\C$ by removing horizontal cuts starting at $s_1,...,s_k$:
$$ U_0 = \C \backslash \bigcup_{j=1}^k (s_k + \R_{\ge 0}). $$
Let $\varepsilon>0$ be such that $|\Im(s_{j}-s_{j'})|>\varepsilon$ if $j\ne j'$. 

Fix $j\in \{1,...,k\}$. For any $\eta\ge 0$, let $F_{j,\eta}^a=\{ s\in \C \ : \ \Re s > \Re s_j, \Re s_j+\eta > \Im s > \Re s_j-\varepsilon\}$. Define $U$ by identifying $F_{j,\eta}^a$ and $U_0$ along $F_{j,0}^a$. We say that $U$ is obtained from $U_0$ by {\it attaching a flap of size $\eta$ above $s_j$ along the cut $(s_j,+\infty)$}. We refer to the set $U\backslash U_0$ as the {\it flap} itself.

We can analogously attach a flap of size $\eta$ {\it below $s_j$} along the cut $(s_j,+\infty)$, or simultaneously several flaps above and/or below some singularities among $s_1,...,s_k$. The obvious map $U\to \C$ is then locally biholomorphic. 

\begin{figure} \includegraphics{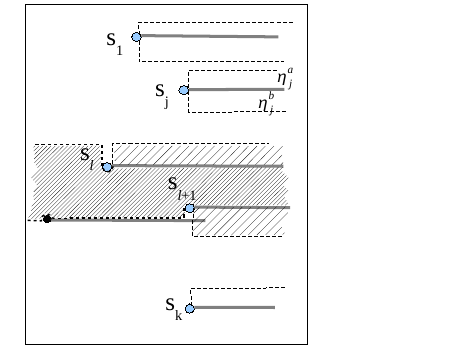} \caption{Strip between $s_\ell$ and $s_{\ell+1}$} \label{Paper3p14} \end{figure}

\paragraph*{Strips.} Let $U$ be obtained from $U_0$ by attaching flaps of sizes $\eta_j^a\ge 0$ above and $\eta_j^b$ below the cuts $(s_j,+\infty)$, fig.\ref{Paper3p14}

We are going to define the {\it strip between singularities $s_\ell$ and $s_{\ell+1}$}, $1\le \ell \le k-1$ as the intersection of $A_{\ell+1}\cap B_{\ell}$,
$$ \begin{array}{clc} A_\ell=\{ s\in \C \ : & \Im s > \Im s_\ell - \eta^a_\ell \\ 
                                            & (\Re s=\Re s_\ell) \Rightarrow (\Im s > \Im s_{\ell}) \\
                                            & ( \Re s< \Re s_\ell \ \text{and} \ \Re s_{\ell'}<\Re s_{\ell} \ \text{and} \ \ell'>\ell ) \Rightarrow (\Im s>\Im s_{\ell'}) & \} \end{array}; $$
$$ \begin{array}{clc} B_\ell=\{ s\in \C \ : & \Im s < \Im s_\ell + \eta^a_\ell \\ 
                                            & (\Re s=\Re s_\ell) \Rightarrow (\Im s < \Im s_{\ell}) \\
                                            & ( \Re s< \Re s_\ell \ \text{and} \ \Re s_{\ell'}<\Re s_{\ell} \ \text{and} \ \ell'<\ell ) \Rightarrow (\Im s<\Im s_{\ell'}) & \} \end{array}. $$
The {\it (semi-infinite) strip above $s_1$} is defined to be $A_1$, the {\it (semi-infinite) strip below $s_k$} is defined to be $B_k$. Strips can be viewed as subsets of $\C$ or as subsets of $U$, as will be clear from the context.

If $s_0\in \C$ and $\Im s_\ell<\Im s_0<\Im s_{\ell +1}$, the strip between $s_0$ and $s_{\ell+1}$ is defined as $A_{\ell+1} \cap \{ s\in \C: \Im s<\Im s_0\}$. Strips between $s_\ell$ and $s_0$, resp., between $s_0$ and $s_1$, resp, between $s_k$ and $s_0$ are defined analogously.

\paragraph*{Charts $(a\pm i0,b\pm i0)$, $(a\pm i0, b\pm \underline{i0})$, $(a\pm i0)$, $(a\pm \underline{i0})$.} \label{PMiCharts} 
 For a fixed number $\delta>0$, suppose $a,b\in \C$ such that $|\Im (a-b)|<\delta$, $\Re a\le \Re b$.
Define the following sets which we will call the {\it charts $(a+ i0,b+i0)$}, etc.
$$ U_{(a-i0, b-i0)} = \left\{ s\in\C \ : \begin{array}{l} \Im(s) < \Im (a) + \delta \\
                                                          (\Re s > \Re b) \Rightarrow (\Im (s)< \Im(b) + \delta) \\
                                                           \Im (s) > \max\{\Im(a),\Im(b)\}-\delta  \end{array}  \right\}\backslash ((a+i\R_{\ge 0})\cup (b+i\R_{\ge 0})), $$
$$ U_{(a-i0, b-\underline{i0})} = U_{(a-i0, b-i0)}\backslash (b+ (\R_{\ge 0}\times i\R_{\ge 0})), $$
$$U_{(a+i0, b+i0)}= U^*_{(a^*-i0, b^*-i0)}, \ \ U_{(a+i0, b+\underline{i0})}= U^*_{(a^*-i0, b^*-\underline{i0})}.  $$
where a star means the complex conjugate. Further,
$$ U_{(a+i0, b-i0)} = \left\{ s\in\C \ : \begin{array}{l} \Im(s) < \min \{\Im (a), \Im (b)\} + \delta \\
                                                          \Im(s) > \max \{\Im (a), \Im (b)\} - \delta \\
                                                           (a\in b+i\R_{\ge 0}) \Rightarrow (\Re s< \Re a)   \end{array}  \right\}\backslash ((a-i\R_{\ge 0})\cup (b+i\R_{\ge 0})), $$
$$ U_{(a+i0, b-\underline{i0})} = U_{(a+i0, b-i0)}\backslash (b+ (\R_{\ge 0}\times i\R_{\ge 0})), $$
$$U_{(a-i0, b+i0)}= U^*_{(a^*+i0, b^*-i0)}, \ \ U_{(a-i0, b+\underline{i0})}= U^*_{(a^*+i0, b^*-\underline{i0})}. $$
Finally,
$$ U_{(a\pm i0)} = U_{(a\pm i0,a\pm i0)}; \ \ \ U_{(a\pm \underline{i0})} = U_{(a\pm i0,a\pm \underline{i0})}. $$

\paragraph*{Slots.} For $s_0\in\C$, $\rho\ge 0$, define the {\it upward-facing slot of size $\rho$ around $s_0$}
$$Sl^\cup_\rho = \{ s\in \C : |s-s_0|\le \rho \} \cup \{ s\in \C : |\Re (s - s_0)|\le \rho \ \text{and} \ \Im s \ge \Im s_0 \} $$ 
and the {\it downward-facing slot of size $\rho$ around $s_0$} as
$$Sl^\cap_\rho = \{ s\in \C : |s-s_0|\le \rho \} \cup \{ s\in \C : |\Re (s - s_0)|\le \rho \ \text{and} \ \Im s \le \Im s_0 \}. $$



\section{Structure of the Riemann surface.} \label{StructureRS}

The structure of a manifold ${\cal S}$ on which solutions of the equation \eqref{MainEqu1} should be defined is predicted in the  heuristic argument of Voros ~\cite{V83} or, alternatively, by a reasoning generalizing section \ref{AppearsStokes} . We will now take that prediction as an {\it a priori} definition of ${\cal S}$. 

The manifold ${\cal S}$ will be endowed with a projection $\pi:{\cal S} \to \tilde{\cal O}$; we will usually denote by $x$ a point of $\tilde {\cal O}$, by $s$ a point of ${\cal S}_{x}=\pi^{-1}(x)$ and by $(s,x)$ a point of ${\cal S}$. 

In this section we will describe the fibers ${\cal S}_x$ for all $x\in \tilde{\cal O}$ by gluing them from subsets of $\C$; we will take the obvious structure of a complex two-dimensional manifold on ${\cal S}=\cup_{x\in \tilde{\cal O}}{\cal S}_x$. In the rest of the paper, subsets of ${\cal S}$, resp., of ${\cal S}_x$, will often be identified with subsets of $\C^2$ with coordinates $(s,x)$, resp., with subsets of $\C$ with the coordinate $s$; we hope the identification will be clear in each instance.

We will begin by describing a subset ${\mathcal S}'\subset {\mathcal S}$.

\begin{figure} \includegraphics{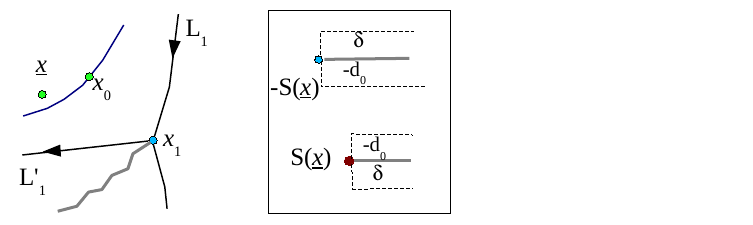} \caption{
Fiber of ${\mathcal S}'$ over $\xu$ when $\xu\in {\cal U}_{L_0,\delta/2}$ and $\Im S(\xu) \le 0$.} \label{ShStFixp19L0} 
\end{figure}

\begin{figure} \includegraphics{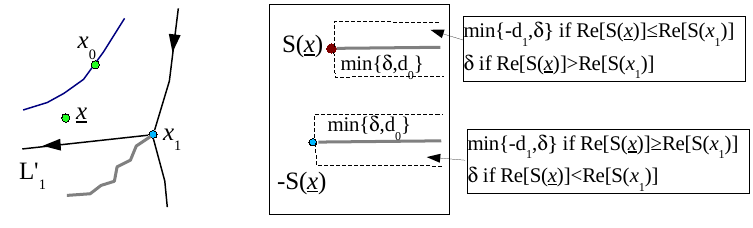} \caption{
Fiber of ${\mathcal S}'$ over $\xu$ when $\xu\in \overline{\cal A}$ and $\Im S(\xu) \ge 0$.} \label{ShStFixp19a} 
\end{figure}

\begin{figure} \includegraphics{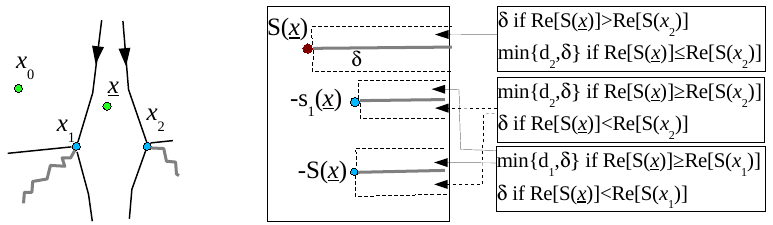} \caption{
Fiber of ${\mathcal S}'$ over $\xu$ when $\xu\in \overline{\cal B}$.} \label{ShStFixp19b} 
\end{figure}

\begin{figure} \includegraphics{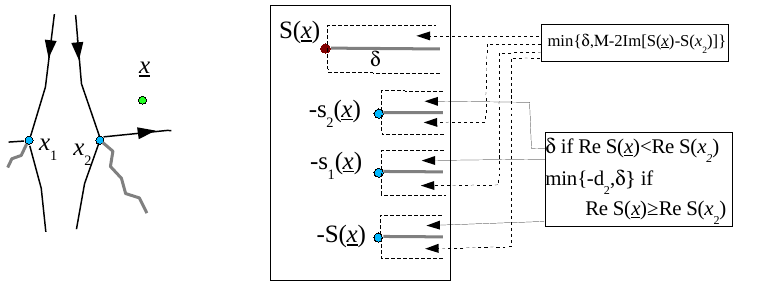} \caption{
Fiber of ${\mathcal S}'$ over $\xu$ when $\xu\in \overline{\cal C}$.} \label{ShStFixp19c} 
\end{figure}

\begin{figure} \includegraphics{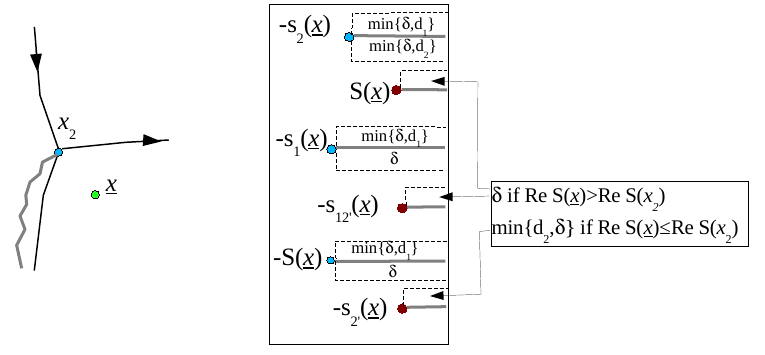} \caption{
Fiber of ${\mathcal S}'$ over $\xu$ when $\xu\in \overline{\cal D}$.} \label{ShStFixp20d} 
\end{figure}

\begin{figure} \includegraphics{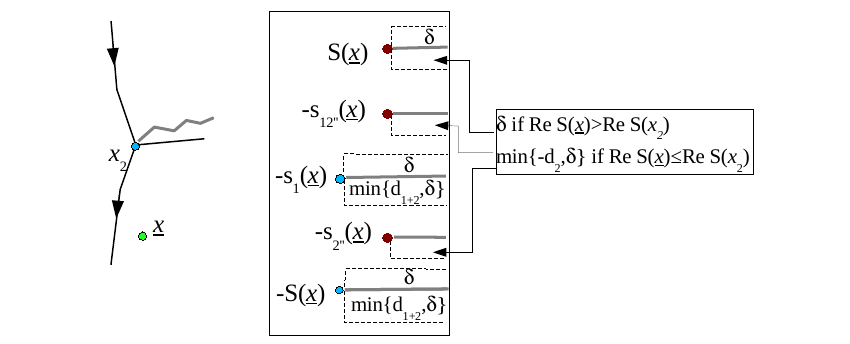} \caption{
Fiber of ${\mathcal S}'$ over $\xu$ when $\xu\in \overline{\cal E}$.} \label{ShStFixp20e} 
\end{figure}

\begin{figure} \includegraphics{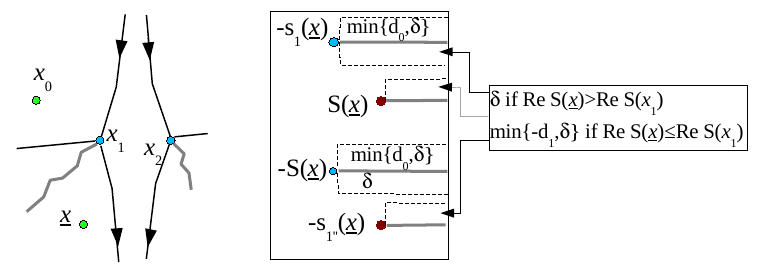} \caption{
Fiber of ${\mathcal S}'$ over $\xu$ when $\xu\in \overline{\cal F}$.} \label{ShStFixp20f} 
\end{figure}

\begin{figure} \includegraphics{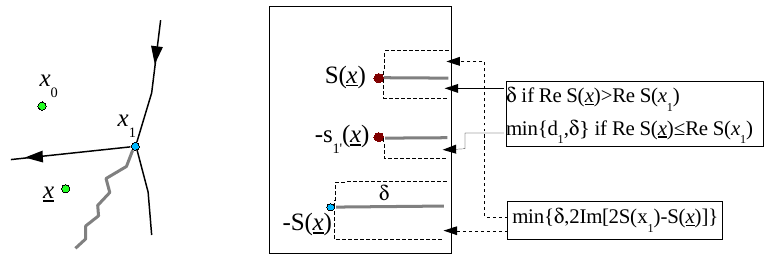} \caption{
Fiber of ${\mathcal S}'$ over $\xu$ when $\xu\in \overline{\cal G}$.} \label{ShStFixp20g} 
\end{figure}

For $\xu\in \tilde{\mathcal O}$ define the fiber of ${\mathcal S}'$ over $\xu$, as on figures  \ref{ShStFixp19L0}--\ref{ShStFixp20g}. 

On these pictures $d_1=2\Im [S(\xu)-S(x_1)]$ and $d_2=2\Im[S(x_2)-S(\xu)]$, 
  $d_{1+2}=2\Im[2S(x_2)-S(x_1)-S(\xu)]$, $d_0=2\Im S(x)$ for the determinations of $S$ in the corresponding Stokes regions. The arrows on the Stokes curves indicates the direction in which $\Re S$ grows.

On these pictures the corresponding fibers are given as complex planes with a few singularities, the locations of the singularities are marked. There are two groups of singularities: 
\begin{equation} \text{\parbox{13cm}{\begin{itemize}
 \item the singularities of the form $-S(\xu)+const$, or ``blue'' singularities : $-S(\xu)$, $-s_1(\xu)=2S(x_1)-S(\xu)$, $-s_2(\xu)=2S(x_2)-S(\xu)$
 \item the singularities of the form $S(\xu)+const$, or ``red'' singularities $S(\xu)$, $-s_{2'}(\xu)=-2S(x_2)+S(\xu)$,  $-s_{12'}(\xu)=2S(x_1)-2S(x_2)+S(\xu)$, $-s_{2''}(\xu)=-2S(x_2)+S(\xu)$, $-s_{12''}(\xu)=2S(x_1)-2S(x_2)+S(\xu)$, $-s_{1'}(\xu)=-2S(x_1)+S(\xu)$.
\end{itemize}}} \label{ListOfSing} \end{equation}
The fibers ${\cal S}'_\xu$ are complex planes with cuts made in the positive real direction and with flaps attached (section \ref{FSS}) on both sides of most of the cuts; the sizes of the flaps are specified on the pictures. 

Note that unless $\xu\in {\cal U}_{L,\delta/2}$ for some $L\in {\cal L}$, we have chosen the sizes of all flaps in ${\cal S}'$ to be $\delta$. If $\xu \in {\cal U}_{L,\delta/2}$, there are one or several pairs of singularities from the list \eqref{ListOfSing}, one red and one blue, near the same horizontal line in the $s$-plane; attaching a rectangular flap of size $\delta$ would lead to a wrong definition of ${\cal S}_\xu$. Instead, for $\xu\in {\cal U}_{L,\delta/2}$, $L\in \{L_0\}\cup {\cal L}_{iS}$, we will  glue to  ${\mathcal S}'_\xu$ additional subsets as described below, and the resulting manifold will be our ${\mathcal S}_\xu$.

\private{\textcolor{blue}{REPHRASE and include in the text if possible.} 
The sizes of the \textcolor{Fuchsia}{flaps} are chosen according to the following principle: they are $\delta$ except for $\xu\in {\cal U}_{L,\delta/2}$, $L\in {\cal L}$, in which case some another singularity approaches the cut on the second sheet from the given side, in which the \textcolor{Fuchsia}{flap} is drawn right up to the singularity \textcolor{blue}{\bf clumsy sentence, maybe an example will help; what means second sheet? what means right up to? they are of {\it size} $\delta$}. \\
With these definitions, when a $\xu\in {\cal U}_{L,\delta/2, \textcolor{red}{\Box}}$, $L\in {\cal L}$, the fiber of ${\mathcal S}'$ over $\xu$ has flaps that are too thin for our purposes. E.g., ${\mathcal S}'$ does not ``see'' the singularity $-s_1(\xu)$ on the second sheet just before $\xu$ crosses $L_1$ from the region A to the region B. In order to include this information, for $\xu\in {\cal U}_{L,\delta/2, \textcolor{red}{\Box}}$, $L\in \{L_0\}\cup {\cal L}_{iS}$, we will  glue to the fiber of ${\mathcal S}'$ over $\xu$ additional subsets as described below, and the resulting open manifold will be our ${\mathcal S}$. }

Let $L\in {\cal L}_{iS}$ be a Stokes curve between two Stokes regions ${\cal X}$ and ${\cal Y}$ and let ${\cal Y}$ be closer to $x_0$ than ${\cal X}$. Suppose $-s_a(\xu)=a-S(\xu)$ is a ``blue" and $-s_b(\xu)=b+S(\xu)$ is a ``red" singularity on the first sheet of ${\mathcal S}'_\xu$ for $\xu\in {\cal Y}$. Suppose for definiteness that  $\Re [-s_a(\xu)] < \Re [-s_b(\xu)]$ and $\Im [-s_a(\xu)] = \Im [-s_b(\xu)]$ when $\xu$ is on $L$. In this situation we say that  {\it the singularity $-s_b(\xu)$ appears from under (or disappears under) the cut starting at $-s_a(\xu)$ when $\xu$ crosses $L$.} For example, the red singularity $-s_{1'}$ appears from under the cut starting at the blue singularity $-S(\xu)$ when $\xu$ crosses $L'_1$ from ${\cal A}$ to ${\cal G}$, or the blue singularity $-s_1(x)$ appears from under the cut starting at the red singularity $S(\xu)$ when $\xu$ crosses $L_1$ from ${\cal A}$ to ${\cal B}$. If a singularity at $-s_b(\xu)$ is present on the first sheet for $\xu$ on both sides of $L$, we say that one singularity located at $-s_b(\xu)$ disappears under and another singularity located at $-s_b(\xu)$ appears from under the cut starting at $-s_a(\xu)$ when $\xu$ crosses $L$; take, for example, two singularities located at $S(\xu)$, one appearing and one disappearing under the cut starting at $-s_1(\xu)$ when $\xu$ crosses $L''_1$. If $L=L_0$, one can use similar terminology.

\begin{figure} \includegraphics{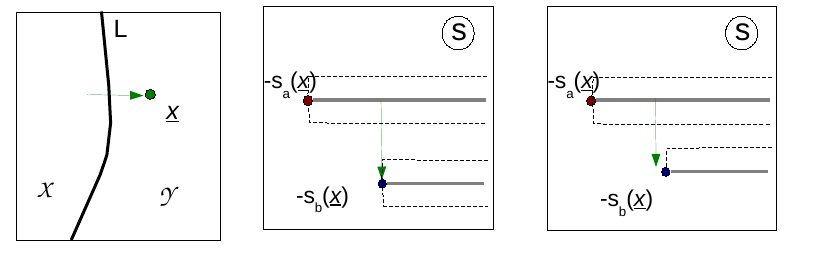} \caption{When $\xu$ crosses a Stokes curve $L$ into  from the Stokes region ${\cal X}$ into the Stokes region ${\cal Y}$ (left), a singularity $-s_b(\xu)$ appears on the first sheet from under the cut starting at $-s_a(\xu)$. The cut $[-s_b(\xu),\infty)$ may have two flaps (center) or just one (right).} \label{Paper3p8} \end{figure} 

In the notation of the previous paragraph, most frequently it happens that when $\xu$ is in the region ${\cal Y}$, the cut $[-s_b(\xu),+\infty)$  on the first sheet of ${\mathcal S}'_{\xu}$ has two flaps, see fig.\ref{Paper3p8}. Then, for $\xu\in {\cal U}_{L,\delta/2}$, we will attach to ${\mathcal S}'_{\xu}$ two subsets according to the procedure which we are going to describe on the example of the blue singularity $-s_1(x)$ appearing from under the cut starting at the red singularity $S(\xu)$ when $\xu$ crosses $L_1$ from ${\cal A}$ to ${\cal B}$. 
For other such pairs of singularities, one appearing from the cut starting at the other when $\xu$ crosses a Stokes curve, similar subsets must be attached, up to maybe reversing the roles of blue and red singularities and maybe reflecting all pictures with respect to a horizontal line.

 \begin{figure} \includegraphics{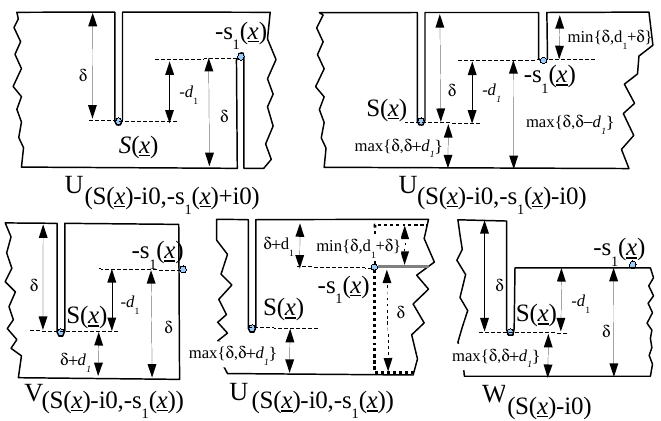} \caption{Attaching additional strips to ${\mathcal S}'$ -- a ``generic'' situation.} \label{Paper3p6} \end{figure}

Consider the subsets of the complex plane $U_{(S(\xu)-i0,-s_{1}(\xu)+i0)}$ and $U_{(S(\xu)-i0,-s_{1}(\xu)-i0)}$, defined on page \pageref{PMiCharts} and shown on figure \ref{Paper3p6}, top.  Glue these two subsets together along a subset $V_{(S(\xu)-i0,-s_{1}(\xu))}$ in their intersection, obtain a set $U_{(S(\xu)-i0,-s_{1}(\xu))}$ whose natural projection to the complex plane will no longer be one-to-one, figure \ref{Paper3p6}, center bottom. Now let consider a subset $W_{(S(\xu)-i0)}\subset \C$ which naturally identifies with subsets of both ${\mathcal S}'$ and $U_{(S(\xu)-i0,-s_{1}(\xu))}$. Attach $U_{(S(\xu)-i0,-s_{1}(\xu))}$ to ${\mathcal S}'$ along $W_{(S(\xu)-i0)}\subset \C$.  

\begin{figure} \includegraphics{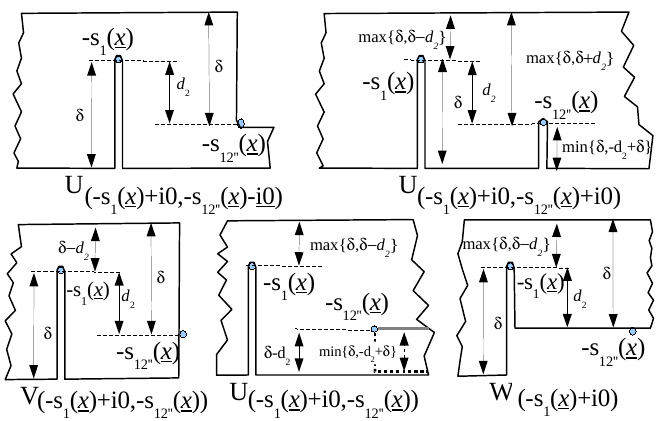} \caption{Attaching additional strips to ${\mathcal S}'$ -- a case when the singularity appearing from under the cut has only one flap.} \label{Paper3p7} \end{figure}

However, there are also cases when the singularity that appears from under the cut has only one flap, e.g.,  the singularity $-s_{12''}(\xu)$ appearing from under the cut starting at $-s_{1}(\xu)$ when  $\xu$ crosses $L''_2$ from ${\cal B}$ to ${\cal E}$. Let $\xu\in {\cal U}_{L''/2,\delta/2}$. Consider the subsets of the complex plane $U_{(-s_1(\xu)+i0,-s_{12''}(\xu)-\underline{i0})}$ and $U_{(-s_1(\xu)+i0,-s_{12''}(\xu)+i0)}$, figure \ref{Paper3p7}. Glue these two subsets along a subset $V_{(-s_1(\xu)+i0,-s_{12''}(\xu))}$ in their intersection, obtain a set $U_{(-s_1(\xu)+i0,-s_{12''}(\xu))}$. A subset $W_{(-s_1(\xu)+i0)}\subset \C$  naturally identifies with subsets of both ${\mathcal S}'$ and $U_{(-s_1(\xu)+i0,-s_{12''}(\xu))}$. Attach $U_{(-s_1(\xu)+i0,-s_{12''}(\xu))}$ to ${\mathcal S}'$ along $W_{(-s_1(\xu)+i0)}\subset \C$. 

In particular, for $\xu\in {\cal U}_{L_0,\delta/2}$ satisfying $\Re S(\xu)\ge 0$ we obtain ${\cal S}_\xu$ from ${\cal S}'_{\xu}$ by attaching to it four sets $U_{(-S(\xu)\pm i0,S(\xu) \pm i0)}$; for $\xu\in {\cal U}_{L_0,\delta/2}$ satisfying $\Re S(\xu)\le 0$ we obtain ${\cal S}_\xu$ from ${\cal S}'_{\xu}$ by attaching to it four sets $U_{(S(\xu)\pm i0,-S(\xu) \pm i0)}$.

This finishes the description of ${\mathcal S}$.  

{\bf Remarks.} 
\begin{enumerate} 
\item As we see, the size of the flaps is controlled by the parameter $\delta$. For  $\delta=0$,   the fibers ${\mathcal S}_{\xu}$ are subsets of the complex plane and constructions of this article simplify; it may be helpful for a reader to first understand the rest of the article in this case. 
\item There is a tempting idea to prove the existence of analytic continuations of terms of \eqref{vNsimpl} to the first sheet of ${\mathcal S}$ only (which is equivalent to setting $\delta=0$ and therefore simpler), and then repeat the same construction for what ~\cite{DP99} call {\it a (re)summation direction $\alpha$} for a small enough angle $\alpha$. That would involve drawing cuts on the $s$-plane in the direction $e^{i\alpha}$ (rather than in the positive real direction $e^{i0}$), defining Stokes curves by conditions $S(x)-S(x_j)\in e^{i\alpha}\R$ and correspond to the asymptotics of type \eqref{hyperasexpn} but for $|h|\to 0+$, $\arg h=\alpha$. This idea, as far as we can make it work, indeed provides analytic continuations beyond the first sheet of ${\mathcal S}_{\xu}$, but only for those $\xu$ that are far enough from the Stokes curves.
\item The ``Riemann surface" ${\mathcal S}$ constructed above has the following property: If the canonical distance from  $\xu$  to the points $x_1,x_2$ and to all curves among $\{ L_0\} \cap {\cal L}_{eS}$ is greater than $\delta/2$, every singularity present on the first sheet of ${\mathcal S}_\xu$ has a sector around it of radius $\delta$ and aperture $\ge \frac{5\pi}{2}$ inside ${\mathcal S}_{\xu}$. This uniform estimate of the ``size'' of ${\mathcal S}$ can potentially be useful for the applications of the Watson's lemma to the calculation of hyperasymptotic expansions (and estimating their error terms) as described in the Introduction and on the fig.\ref{Paper3p3}. In particular, although some singularities on ${\cal S}$ have a flap attached only on one side, this is still sufficient for the purposes of deforming the contour of the Laplace integral (\ref{LaplaceTransf}).
\end{enumerate}

\section{Construction of analytic continuations.} \label{PfMresSec}


\subsection{Strategy of the proof. } \label{StrategySec}

The definitions  in the sections \ref{Notation}, \ref{StructureRS} have given a precise sense to content of the section \ref{RedToPaths}; we are continuing now where we stopped at the end of the section \ref{RedToPaths}. 

In order to prove theorem  \ref{RjGmainTh}, we will construct the integration paths for $R_j$ and $(s,\xu)\in {\cal S}$ from $x_0$ to $\xu\in A$, $A\in {\mathbb S}$, cf. \eqref{tiling}, step by step. First we will consider the case $A= {\cal A}_{int}$. After that, if $x\in A\in {\mathbb S}$, $A\ne {\cal A}_{int}$, we will construct a piece of integration path leading from an $\xu'\in B$ to $\xu$, where $B\in {\mathbb S}$ is closer to $x_0$ than $A$. This will provide an argument similar to inductive. \label{IndArgProposed}

In the case when $A={\cal X}_{int}$ for some Stokes region ${\cal X}$, $\xu \in A$, the set $U\subset {\cal S}_\xu$ will be chosen as one of the horizontal strips in ${\cal S}_\xu$ in the sense of section \ref{FSS}. Let us explain now how we will choose the sets $U$ in the other case, i.e. if $\xu\in A = {\cal U}_{L,\delta/2}$ for an internal Stokes curve $L\in {\cal L}_{iS}$. 

Fix one of the choices $j=1$ or $2$ and suppose the Stokes curve $L$ starts at the turning point $x_t$ and separates two Stokes regions ${\cal X}$ and ${\cal Y}$, and ${\cal X}$ is closer to $x_0$ than ${\cal Y}$. Consider those singularities from the list \eqref{ListOfSing} which appear on the first sheet of ${\cal S}_\xu$ when $\xu\in {\cal Y}$, and group them as follows, using definition \ref{StatMovingDef}:\\
a) pairs $\sigma_\nu(x), \sigma'_\nu(x)$ ($\nu=1,2,...,N_0$), where $\sigma_\nu$ is stationary and $\sigma'_\nu$ is moving and $\sigma_\nu(x_t)=\sigma'_\nu(x_t)$; \\
b) moving singularities $\tau_\nu(x)$ ($\nu =N_0+1,N_0+2,..,N_1$) and stationary singularities $\upsilon_\nu(x)$ ($\nu=N_1+1,N_1+2,...,N_2$) which cannot be included into pairs as in item a).

For definiteness suppose that $\Re S_j$ is growing away from $x_t$ along $L$ and that one has to go clockwise around $x_t$ to cross $L$ from ${\cal X}$ to ${\cal Y}$; other cases can be considered analogously. Then $\Im \sigma'_\nu(x)>\Im \sigma_\nu(x)$ for $x\in {\cal U}_{L,\delta/2}\cap {\cal X}$. Let for $x\in {\cal U}_{L,\delta/2}$
$$ E_{\nu,x} = \left\{ s\in \C \ : \ \begin{array}{l} \Re s\ge \min \{ \Re \sigma_\nu(x), \Re \sigma'_\nu (x) \}; \\
 \ \Im s\le \max \{ \Im \sigma_\nu(x), \Im \sigma'_\nu(x)\} \\ \Im s \ge \Im (\sigma_\nu(x))-\delta \end{array} \right\}, \ \text{for} \ \nu=1,...,N_0;  $$
$$ E_{\nu,x} = \left\{ s\in \C \ : \ \begin{array}{l} \Re s\ge \min \{ \Re 2\tau(x_t)-\tau_\nu(x), \Re \tau_\nu (x) \}; \\
 \ \Im s\le \max \{ \Im 2\tau(x_t)-\tau_\nu(x), \Im \tau_\nu(x)\} \\ \Im s \ge \Im( 2\tau_\nu(x_t) - \tau_\nu(x))-\delta \end{array} \right\}, \ \text{for} \ \nu=N_0+1,...,N_1;  $$
$$ E_{\nu,x} = \upsilon_\nu(x) + \R_{\ge 0}, \ \text{for} \ \nu=N_1+1,...,N_2.$$
Our conditions on $\delta$ in \eqref{DeltaIsSmall} imply that $E_{\nu,x}\cap E_{\mu,x}=\emptyset$ if $\mu\ne \nu$. Let $E_{0,x}=\C\backslash \bigcup_{\nu=1}^{N_2} E_{\nu,x}$ and identify $E_{0,x}$ with a subset of ${\cal S}_x$. 

\private{ Let $E_{0,x}^{\delta'}$ be defined analogously to $E_{0,x}$ with $\delta$ replaced everywhere by $0<\delta'<\delta$. Then $E_{0,x} = \bigcup_{\delta'<\delta} E^{\delta'}_{0,x}$. \textcolor{blue}{Insert this comment where it belongs} }

If $\xu\in {\cal U}_{L,\delta/2}$, consider $U=E_{0,\xu}\subset {\cal S}_\xu$ and take as a piece of an integration path a path $y(t)$ with $\Im S(y(t))=\Im S(\xu)$ ending at $\xu$ and starting at some $\xu'\in {\cal X}\cap {\cal U}_{L,\delta} \backslash {\cal U}_{L,\delta/2}$. Using conditions on $\delta$ in \eqref{DeltaIsSmall}, one finds $U$ can be transported within ${\cal S}$ along $y(t)$ parallel to $-S_j$. 

\private{\textcolor{blue}{Saying something about thickening of the strips $E_\nu$ twice might be helpful}}

\private{\textcolor{blue}{In the following paragraph, one should say that shrinking $\delta$ no more than twice will still allow us to draw an integration path to ${\cal X}_\int$.}}

After replacing $\delta$ in the definition of $E_{0,\xu}$ with a $\delta'<\delta$ such that ${\cal U}_{L,\delta'/2}$ still contains $\xu$, one can cover ${\cal S}_\xu\backslash E_{0,\xu}$ by charts of the type $U_{(\sigma_\nu(\xu)\pm i0, \sigma'_\nu(x))}$, etc. introduced on page \pageref{PMiCharts}. For $\xu\in {\cal U}_{L,\delta/2}$, we will take $U_{(\sigma_\nu(\xu)\pm i0, \sigma'_\nu(x))}$, etc., as $U$ when presenting a construction of the integration paths. 

\label{Combina} \paragraph{Combinatorial complexity of the problem.} Let us make a rough guess how much work is needed to prove the theorem \ref{RjGmainTh} by this method for a general potential $V(x)$ (not necessarily the same as we are studying in this article) and a corresponding piece of the ``Riemann surface" ${\cal S}$. Let ${\cal Y}$ be a Stokes region and $\ell$ be a Stokes curve separating it from an earlier Stokes region.  For $\xu\in {\cal Y}_{int}$, the sets $U$ are mostly strips between pairs of singularities: the upper boundary of a strip is defined  either by a moving or by a stationary singularity, the lower boundary of the strip is defined either by a moving or by a stationary singularity, and $\Re S_j(x)$ may be either increasing or decreasing along $\ell$, this makes $2\times 2\times 2=8$ cases. For $x\in {\cal U}_{\ell,\delta}$, the charts $U_{(\sigma_1(x)\pm i0, \sigma_2(x)\pm i0)}$ have two independently chosen $\pm$ signs, and either $\sigma_1$ is a moving and $\sigma_2$ is a stationary singularity, or $\sigma_1$ is stationary and $\sigma_2$ is moving, and again $\Re S_j$ is either increasing or decreasing along $\ell$, this makes $2^2\times 2\times 2=16$ more cases. On the one hand, not all of these possibilities are actually realized, but on the other hand, there are also some degenerate cases such as semi-infinite strips, so in total we expect some twenty or thirty separate cases to consider. Some of these cases can be treated together by the same argument, some subdivide a little further, but qualitatively this estimate justifies the number of lemmas in the section \ref{modelcases}. 






\subsection{Proof for the region ${\cal A}_{int}$.} \label{ProofAint}

We will work with $R_1G$, the argument for $R_2G$ being completely analogous.

\begin{figure} \includegraphics{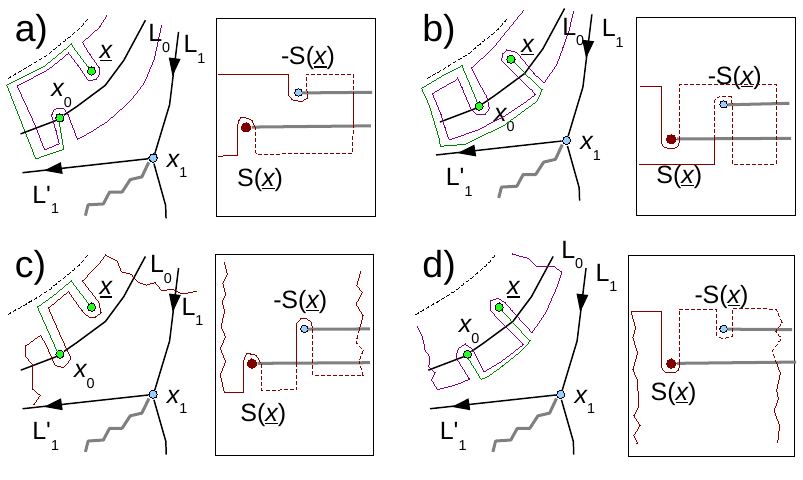} \caption{The sets  a) ${\tilde U}^{\varepsilon,N}_{(S(\xu)+i0,-S(\xu)-i0)}$, b) $ {\tilde U}^{\varepsilon,N}_{(S(\xu)-i0,-S(\xu)+i0)}$, c) $T^{\varepsilon}_{(S(\xu)+i0,-S(\xu)+i0)}$, d) $T^{\varepsilon}_{(S(\xu)-i0,-S(\xu)-i0)}$ in ${\cal S}_\xu$, the corresponding sets $V$ in ${\cal O}$, and the integration paths from $x_0$ to $\xu$. The curve $\Im S(x)=-\delta/2$ is drawn with a dashed line.}  \label{Paper3p40} \end{figure}

Let $\xu\in {\cal U}_{L_0,\delta/2}$; we will construct an integration path $y(t)$ from $x_0$ to $\xu$ for every $(s,\xu)\in {\cal S}_\xu$. Suppose for definiteness that $\Re S(\xu)\le 0$. Cover ${\cal S}_\xu$ as follows:
$$ {\cal S}_\xu = \bigcup_{N,\varepsilon>0} \left( {\tilde U}^{\varepsilon,N}_{(S(\xu)+i0,-S(\xu)-i0)} \cup {\tilde U}^{\varepsilon,N}_{(S(\xu)-i0,-S(\xu)+i0)} \cup  T^{\varepsilon}_{(S(\xu)+i0,-S(\xu)+i0)} \cup T^{\varepsilon}_{(S(\xu)-i0,-S(\xu)-i0)} \right), $$
where 
$$ {\tilde U}^{\varepsilon,N}_{(S(\xu)+i0,-S(\xu)-i0)} =\left\{ s\in\C \ : \begin{array}{l} \Im(s) < \min \{\Im (S(\xu)), \Im (-S(\xu))\} + \delta-\varepsilon \\
                                                          \Im(s) > \max \{\Im (S(\xu)), \Im (-S(\xu))\} - \delta +\varepsilon\\
                                                           (S(\xu)\in i\R_{\ge 0}) \Rightarrow (\Re s< \Re S(\xu)) \\ 
\Re s < N  \end{array}  \right\}\backslash $$ $$ \ \ \ \ \ \ \backslash (Sl^\cap_\varepsilon(S(\xu))\cup Sl^\cup_\varepsilon(-S(\xu))), $$
$$ {\tilde U}^{\varepsilon,N}_{(S(\xu)-i0,-S(\xu)+i0)} =\left\{ s\in\C \ : \begin{array}{l} \Im(s) < \min \{\Im (S(\xu)), \Im (-S(\xu))\} + \delta-\varepsilon \\
                                                          \Im(s) > \max \{\Im (S(\xu)), \Im (-S(\xu))\} - \delta +\varepsilon\\
                                                           (S(\xu)\in -i\R_{\ge 0}) \Rightarrow (\Re s< \Re S(\xu)) \\ 
\Re s < N  \end{array}  \right\}\backslash $$ $$ \ \ \ \ \ \ \backslash(Sl^\cup_\varepsilon(S(\xu))\cup Sl^\cap_\varepsilon(-S(\xu))), $$
$$ T^{\varepsilon}_{(S(\xu)-i0, -S(\xu)-i0)} = \left\{ s\in\C \ : \begin{array}{l} \Im(s) < \Im (S(\xu)) + \delta \\
                                                          (\Re s > \Re (-S(\xu)) \Rightarrow (\Im (s)< \Im(-S(\xu)) + \delta)   \end{array}  \right\}\backslash $$ 
$$ \ \ \ \ \ \backslash (Sl^\cup(S(\xu)) \cup Sl^\cup(-S(\xu)) ); $$
$$ T^{\varepsilon}_{(S(\xu)+i0, -S(\xu)+i0)} = \left\{ s\in\C \ : \begin{array}{l} \Im(s) > \Im (S(\xu)) - \delta \\
                                                          (\Re s > \Re (-S(\xu)) \Rightarrow (\Im (s)> \Im(-S(\xu)) - \delta)  \end{array}  \right\}\backslash $$ 
$$ \ \ \ \ \ \backslash (Sl^\cap(S(\xu)) \cup Sl^\cap(-S(\xu)) ). $$
We will construct $y(t)$ in the four cases of $s$ contained in the sets  ${\tilde U}^{\varepsilon, N}_{(S(\xu)+i0,-S(\xu)-i0)}$, $ {\tilde U}^{\varepsilon,N}_{(S(\xu)-i0,-S(\xu)+i0)}$, $T^{\varepsilon}_{(S(\xu)+i0,-S(\xu)+i0)}$, or $T^{\varepsilon}_{(S(\xu)-i0,-S(\xu)-i0)}$. If $U$ is one of these four sets, consider $V=S^{-1}\left( \frac{U+S(\xu)}{2} \right)\subset \tilde{\cal O}$, where the branch of $S^{-1}$ is chosen so that $S^{-1}(0)=x_0$ (cf. \eqref{SjInverse}). We have made our definitions in such a way that in all four cases the set ${\cal A}\backslash V$ is connected, and there is a path $y(t)$ from $x_0$ to $\xu$ in ${\cal A}\backslash V$ which can be taken as an integration path, figure  \ref{Paper3p40}.

\label{WhyAnal}
It is easy to see that this defines an {\it analytic} function in a neighborhood of a point $(s,\xu)$, where $s$ belongs, say, to $U=\tilde U^{\varepsilon,N}_{(S(\xu)+i0,-S(\xu)-i0)}$. If $\xu'$ is another point in ${\cal U}_{L_0,\delta/2}$, then the set $U$ can also be transported parallel to $-S_1$ along any path $z(t)$ from $\xu$ to $\xu'$ if $z(t)$ is contained in ${\cal U}_{L_0,\delta/2}\backslash V$. There is an open contractible set ${\cal N}_\xu\subset \tilde {\cal O}$ such that $\xu\in {\cal N}_\xu\subset {\cal U}_{L_0,\delta/2}\backslash V$.  For any $s'\in U$ and $x\in {\cal N}_\xu$, the function  $(R_1G)(s'+S_1(\xu)-S_1(x),x)$ is holomorphic with respect to $x$ because it is an integral of a holomorphic function and with respect to $s'$ because $s'$ is a holomorphic parameter of the integrand, hence, by Osgood theorem, this function is holomorphic in both $s'$ and $x$, and so, after a change of variables, is $(R_1G)(s',x)$ is holomorphic in $(s',x)$ in a neighborhood of $(s,\xu)$.

\begin{figure} \includegraphics{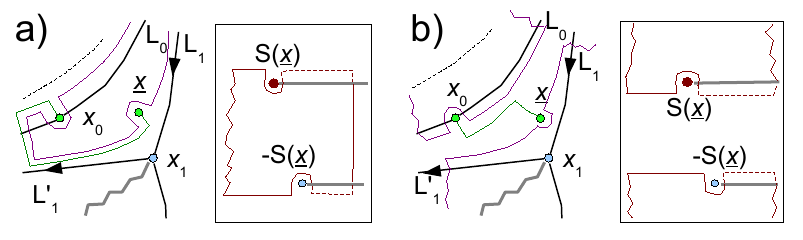} \caption{The sets  a) ${\tilde U}^{\varepsilon,N}_{(S(\xu)+i0,-S(\xu)-i0)}$, b)$W^{+,\varepsilon}_\xu$ and $ W^{-,\varepsilon}_\xu$, the corresponding sets $V$ in $\tilde{\cal O}$, and the integration paths from $x_0$ to $\xu$. The curve $\Im S(x)=-\delta/2$ is drawn with a dashed line.} \label{Paper3p41} \end{figure}

If $\xu\in {\cal A}_{int}\backslash {\cal U}_{L_0,\delta/2}$, 
then 
$$ {\cal S}_x = \bigcup_{N,\varepsilon>0} \left( \tilde U^{\varepsilon,N}_{(S(\xu)+i0,-S(\xu)-i0)} \cup W^{+,\varepsilon}_\xu \cup W^{-,\varepsilon}_\xu \right), $$
where $U^{\varepsilon,N}_{(S(\xu)+i0,-S(\xu)-i0)}$ is the same as above and 
$$ W^{+,\varepsilon}_\xu \ = \ \left\{ s\in \C : \Im s > \Im (S(\xu))+\delta-\varepsilon \right\} \backslash Sl^\cap_\varepsilon(S(\xu)),  $$
$$ W^{-,\varepsilon}_\xu \ = \ \left\{ s\in \C : \Im s < \Im (-S(\xu))-\delta+\varepsilon \right\} \backslash Sl^\cup_\varepsilon(-S(\xu)).  $$
Analogously to the case of $\xu\in{\cal U}_{L_0,\delta/2}$, we will construct $y(t)$ in the three cases of $s$ contained in the sets  ${\tilde U}^{\varepsilon, N}_{(S(\xu)+i0,-S(\xu)-i0)}$, $W^{+,\varepsilon}_\xu$, or $ W^{-,\varepsilon}_\xu$. If $U$ is one of these three sets, consider $V=S^{-1}\left( \frac{U+S(\xu)}{2} \right)\subset \tilde{\cal O}$, where the branch of $S^{-1}$ is chosen so that $S^{-1}(0)=x_0$.  We have made our definitions in such a way that in all four cases the set ${\cal A}\backslash V$ is connected, and there is a path $y(t)$ from $x_0$ to $\xu$ in ${\cal A}\backslash V$ which can be taken as an integration path, figure  \ref{Paper3p41}. Analyticity of $R_1G$ can be shown the same way as above.

\subsection{Continuation to the further Stokes regions -- examples.} \label{examples}

The basic method in the rest of the paper is the same as in section \ref{ProofAint}.


By way of an example, let us take $\xu\in {\cal D}_{int}\backslash {\cal U}_{\tilde L''_2,\delta/2}$ such that $\Im[S(\xu)-S(x_1)]>\delta/2$ and let $U\subset {\cal S}_\xu$ be the strip (see section \ref{FSS}), slightly narrowed and truncated from the right, between the singularities $-s_2(\xu)$ and $S(\xu)$:
$$ U = \left\{ s\in \C \ : \ \begin{array}{l} \Re s < N \\ \Im s < \Im (-s_2(\xu))+\delta-\varepsilon \\ \Im s>\Im (S(\xu)) +\varepsilon \end{array} \right\} \backslash $$
$$ \ \ \ \ \backslash (Sl^\cup_\varepsilon(-s_2(\xu))\cup Sl^\cap_\varepsilon(S(\xu)) ) $$
for a small $\varepsilon>0$ and a large $N>0$. 
For $R_2$  and $(s,\xu)$, $s\in U$, let us construct a piece of the integration path $y(t)$ that starts in some region on the list ${\mathbb S}$ which is closer to $x_0$ than ${\cal C}_{int}$ and ends at $\xu$.  

\begin{figure} \includegraphics{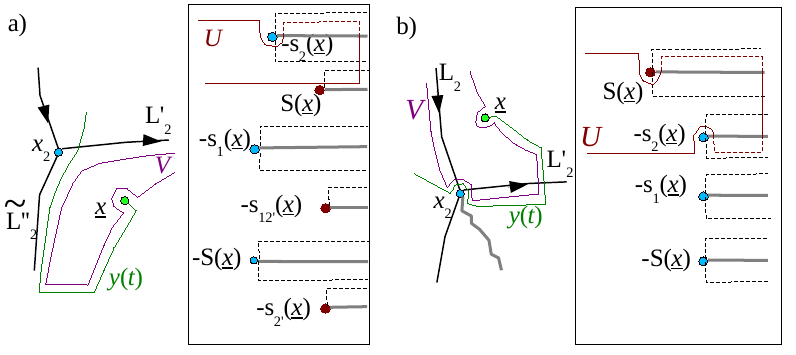} \caption{Integration paths a) for $R_2$ and for $(s,\xu)$ with  $\xu\in {\cal D}_{int}$ and $s$ in the strip between $-s_2$ and $S$. b) for $R_1$ and for $\xu\in {\cal C}_{int}$ and $s$ in the strip between $S$ and $-s_1$.} \label{Paper3p43} \end{figure}

Draw on ${\cal O}$ a set $V=S_2^{-1}\left( \frac{U+S_2(\xu)}{2}\right)$, where the branch of $S_2^{-1}$ is chosen in such a way that $S_2^{-1}(S_2(\xu))=\xu$. A possible choice of a path $y(t)$ starting in ${\cal C}_{int}$ is shown on the figure \ref{Paper3p43},a. Carefully looking at the definition of the ``Riemann surface" ${\cal S}$, one easily checks that $U$ can be transported within ${\cal S}$ parallel to $-S_2$ along this path $y(t)$.

A similar argument can be used for $R_2$, for  $\xu \in {\cal D}_{int}$ and $s$ in the strips between $-s_1$ and $-s_{12'}$, between $-S$ and $-s_{2'}$, for $\xu\in {\cal E}_{int}$ and $s$ in the strips between $-s_{2''}$ and $-S$, for $\xu$ in ${\cal F}$ and $s$ in the strip between $-s_1$ and $-S$, etc. This argument is therefore phrased as lemma \ref{L5bis}. A slightly modified argument, lemma \ref{L5bis}.A and fig.\ref{Paper3p43},b, is necessary, e.g. for $R_1$, $x\in {\cal C}_{int}$ and $s$ in the strip $S(\xu)$ and $-s_2(\xu)$. Lemma \ref{L5bis} is suited for a Stokes region that does border an external Stokes curve, lemma \ref{L5bis}.A -- for a Stokes region that does not. 

\private{In the cases where we apply lemma \ref{L5bis}, we draw the integration path near an external Stokes curve which we may not cross and therefore have a strip with only one flap in ${\cal S}_x$; we apply lemma \ref{L5bis}.A when we draw in integration path near an internal Stokes curve, may cross it and can therefore prove something for a strip in ${\cal S}$ with two flaps. 
}

\private{\textcolor{blue}{\paragraph{Example 2.} \textcolor{red}{Lemma 5.2}}}


\subsection{Continuation to the further Stokes regions -- model cases.} \label{modelcases}


We will now start carrying out the idea described in the section \ref{IndArgProposed}, namely, we will construct integration paths from one element of ${\mathbb S}$ to the next. We will formulate here the basic building blocks of this construction.

Recall that when constructing an integration path for $R_j$, we call singularities from \eqref{ListOfSing} of the form $-S_j(x)+const$ {\it stationary} and those of the form $S_j(x)+const$ {\it moving}.

\begin{figure} \includegraphics{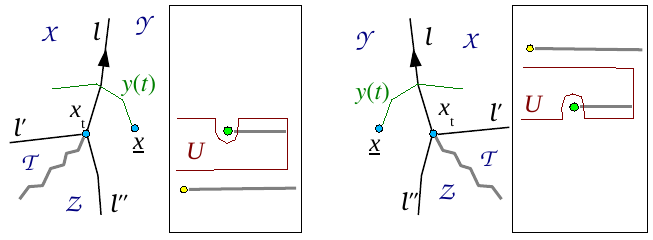} \captionempty 
\label{ShStFixp202a} \end{figure}

In the lemmas below, ${\cal Y}$ denotes one of the Stokes regions ${\cal B}$,..., ${\cal G}$ and $\ell$ is a Stokes curve starting at a turning point $x_t$, $\ell\subset \overline{\cal Y}$ and such that ${\cal U}_{\ell,\delta/2}$ is closer to $x_0$ than ${\cal Y}$; $\ell''$ denotes the other (internal or external) Stokes curve starting from $x_t$ and bordering ${\cal Y}$.  If ${\cal X}$ is mentioned in a lemma, we take it to be the Stokes region on the other side of $\ell$; it is necessarily closer to $x_0$ than ${\cal U}_{\ell,\delta/2}$. If $\ell'$ or ${\cal T}$ are mentioned in a lemma, then $\ell'\subset \overline{\cal X}$, $\ell'$ starts at $x_t$, $\ell'\ne \ell$, and ${\cal T}$ is the Stokes region on the other side of $\ell'$ from ${\cal X}$; we make no assumption whether ${\cal T}$ is closer to $x_0$ than ${\cal X}$. Finally, if ${\cal Z}$ is mentioned, it is the Stokes region on the other side of $\ell''$ from ${\cal Y}$. The lemmas are formulated for the case when the order of ${\cal T},{\cal X},{\cal Y},{\cal Z}$ around $x_t$ is clockwise, fig. \ref{ShStFixp202a}, left. The similar statements for the counterclockwise order are obtained by reversing the direction of the imaginary axis in the $s$-plane, fig. \ref{ShStFixp202a}, right.

Suppose $\Re S(x)$ increases along $\ell$ in the direction away from $x_t$. We will say that ${\cal Y}\subset \tilde{\cal O}$ has {\it canonical width} $\ge C$, where $C\in \R_{\ge 0}$, and write $c.w.({\cal Y})\ge C$ if the function $x\mapsto S(x)-S(x_t)$ defines a bijective map between a subset of ${\cal Y}$ and a set $\{ z\in \C : 0<\Im z < C\}$. Similarly define $c.w.$ for the case when $\Re S(x)$ decreases along $\ell$. As an exception, put $c.w.({\cal A})=\Im S(x_1)$.  It follows from \eqref{DeltaIsSmall} that the canonical width of all Stokes regions is $>3\delta$, so this assumption will not be explicitly added to the lemmas.

In line what we said about the combinatorics of the problem on page \pageref{Combina}, the lemmas we are going to present now can be systematized in the following manner. Let ${\cal Y}$ be a Stokes region and $\ell\subset \overline{\cal Y}$ be a Stokes curve starting at $x_t$ and  such that ${\cal U}_{\ell,\delta/2}$ is closer to $x_0$ than ${\cal Y}$. Firstly, every lemma may pertain to an analytic continuation of $R_jG$ to a set of type ${\cal Y}_{int}$, or to a region of type ${\cal U}_{\ell, \delta/2}$. Secondly, in the situation of every particular lemma, the function $\Re S_j$ can either be increasing or decreasing along a Stokes curve $\ell$ away from the turning point. Some lemmas are useful in several such situations. We have the following table :

\begin{tabular}{p{3.5cm}||p{4cm}|p{4cm}}
& region of type ${\cal Y}_{int}$ & region of type ${\cal U}_{\ell,\delta/2}$  \\ \hline \hline
$\Re S_j$ increasing away from $x_t$ 
& Lemmas \ref{TrivialReason}, \ref{Lemma18}, \ref{Le23}, \ref{Le23}.A, \ref{Lm1bis}, \ref{Lm1bis}.A, \ref{Lm1bis}.B, \ref{C3bo}, \ref{C3T} 
&  
Lemmas \ref{TrivialReason}.C, \ref{Lm1bis}, \ref{Lm1bis}.A, \ref{Lm1bis}.B,   \ref{LN5thin}, \ref{LC10},  \ref{LN23thin}, \ref{LN23thin}.A, \ref{LN23thin}.B 
\\ \hline
$\Re S_j$ decreasing away from $x_t$ & 
Lemmas \ref{TrivialReason}, \ref{TrivialReason}.A, \ref{TrivialReason}.B, ,  \ref{Lemma104}, \ref{LC632}, \ref{LC632}.A, \ref{L5bis}, \ref{L5bis}.A   
& 
Lemmas \ref{LN5thin}, \ref{LN5thin}.A,\ref{LC632}, \ref{LC632}.A,  \ref{L63LB},  \ref{Lemma19}
 \end{tabular}
Organizing the proof by splitting it into sections similar to section \ref{examples} is also possible but it would produce a more verbose document.

As before, $\pi:{\cal S}\to \tilde{\cal O}$ denotes the obvious projection. We will often and freely identify a subset $D\subset {\cal S}$ with a subset of $\C\times \tilde {\cal O}$ or with a subset of $\C\times \C$; we hope that this will not cause any confusion. 

We will use an abbreviation ``a function $G$ is C.A.I. in $D$" to mean that $G$ is continuous on $D$ and analytic in its interior. 

Also recall that $C_s$ and $\C_x$ denote the complex planes of the variables $s,x$, respectively.

With this, let us start presenting the lemmas. 

Construction of the  integration path is especially easy in the following case:

\begin{Lemma} \label{TrivialReason}  
 Suppose $\sigma_0(x)$ is a stationary singularity.
We assume that the function $G$ is constructed in a connected set $D\subset \C_s\times\C_x$ with $\pi(D)\subset {\cal X}\cup \ell \cup {\cal Y}$, and such that for $x\in \pi(D)$, $\pi^{-1}(x)\cap D$ is $\sigma_0(x)+E$ for some fixed subset $E\subset \C$. \\
If $R_jG$ is defined and  C.A.I. on $\pi^{-1}({\cal X})\cap D\ne \emptyset$, 
then $R_j G$ is defined and  C.A.I. on the whole $D$.  \end{Lemma}

\textsc{Proof.} Any path within $D$ starting at a point in ${\cal X}$ can be chosen as the integration path. $\Box$

There are three slight modifications of this lemma:

\paragraph*{Lemma \ref{TrivialReason}.A.} {\it Suppose $\sigma_0(x)$ is a stationary singularity and $\Re \sigma_0$ decreases along $\ell$ in the direction away from $x_t$. Let $\{E_t\}_{t\in [0,\delta]}$ be a family of subsets in $\C$ such that $E_{t_1}\subset E_{t_2}$ if $t_1<t_2$. Let $A>\delta$ be a number, $c.w.{\cal Y}\ge{\cal A}$. 
We assume that the function $G$ is C.A.I. in a connected subset $D\subset \C_s\times\C_x$ with 
$$\pi(D)\subset {\cal X}\cup \ell \cup {\cal Y} ; \ \ \ \pi(D)\cup {\cal Y} = \{ x\in {\cal Y} : \Im \sigma_0(x)-\sigma_0(x_t) <A \} $$  and such that for $x\in \pi(D)$, $\pi^{-1}(x)\cap D$ is described as follows: 
\begin{itemize} 
\item For $x\not \in {\cal U}_{\ell'',\delta/2}$, let $\pi^{-1}(x)\cap D = \sigma_0(x)+E_\delta$;  
\item For $x\in {\cal U}_{\ell'',\delta/2}$, let $\pi^{-1}(x) \cap D = \sigma_0(x) + E_{\frac{1}{2}\Im [\sigma_0(x)-\sigma_0(x_t)]}$.
\end{itemize} 
If $R_jG$ is defined and C.A.I. on $\pi^{-1}({\cal X})\cap D\ne\emptyset$, 
then $R_j G$ is defined and C.A.I. on the whole $D$. }

\paragraph*{Lemma \ref{TrivialReason}.B.} {\it  Suppose $\sigma_0(x)$ is a stationary singularity and $\Re \sigma_0$ can either increase or decrease along $\ell$ in the direction away from $x_t$. Let $\{E_t\}_{t\in [0,\delta]}$ be a family of subsets in $\C$ such that $E_{t_1}\subset E_{t_2}$ if $t_1<t_2$. Let $A>\delta$ be a number, $c.w.({\cal Y})\ge A$. 
We assume that the function $G$ is C.A.I. in a connected subset $D\subset \C_s\times\C_x$ with 
$$\pi(D)\subset {\cal X}\cup \ell \cup {\cal Y} ; \ \ \ \pi(D)\cup {\cal Y} = \{ x\in {\cal Y} : |\Im \sigma_0(x)-\sigma_0(x_t)| <A \} $$  and such that for $x\in \pi(D)$, $\pi^{-1}(x)\cap D$ is described as follows: 
\begin{itemize} 
\item For  $|\Im \sigma_0(x)-\sigma_0(x_t)|\le A-\delta$, let $\pi^{-1}(x)\cap D = \sigma_0(x)+E_\delta$;  
\item For $|\Im \sigma_0(x)-\sigma_0(x_t)|> A-\delta$, let $\pi^{-1}(x) \cap D = \sigma_0(x) + E_{\frac{1}{2}(A-\Im [\sigma_0(x)-\sigma_0(x_t)])}$.
\end{itemize} 
If $R_jG$ is defined and C.A.I. on $\pi^{-1}({\cal X})\cap D\ne \emptyset$, 
then $R_j G$ is defined and C.A.I. on the whole $D$. }

\paragraph*{Lemma \ref{TrivialReason}.C.} {\it  Under the same assumption as in the lemma \ref{TrivialReason}, 
If $R_jG$ is defined and C.A.I. on $\pi^{-1}({\cal X}\backslash {\cal U}_{\ell,\delta/2})\cap D\ne \emptyset$, 
then $R_j G$ is defined and C.A.I. on the whole $D$.  }



\begin{figure} \includegraphics{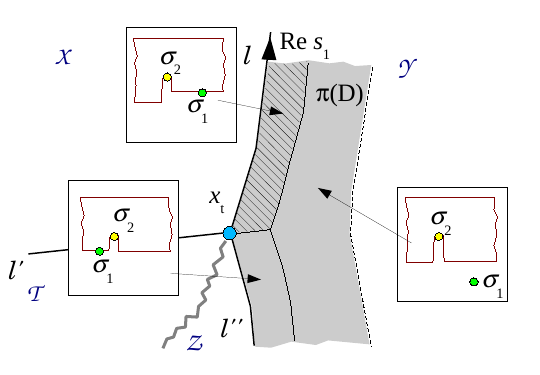} \caption{ Lemma \ref{Lemma18}. The set $\pi(D)$ is shown, as well as $D_x$ for $x\in {\cal Y}\cap{\cal U}_{\ell,\delta/2}$, $x\in {\cal Y}\cap{\cal U}_{\ell'',\delta/2}$, and for $x\in \pi(D)\backslash ({\cal U}_{\ell,\delta/2}\cap {\cal U}_{\ell'',\delta/2})$. The set ${\cal Y}\cap {\cal U}^\circ_{\ell,\delta/2}$ is hatched.
} \label{Paper3p21} \end{figure}

\begin{Lemma} \label{Lemma18} 
Suppose $\sigma_2$ is a stationary and $\sigma_1$ is a moving singularity and $\sigma_2(x_t)=\sigma_1(x_t)$, and $\Re \sigma_1$ grows along $\ell$ in the direction away from $x_t$.
 Let $2c.w.({\cal Y})\ge A>\delta$, $B>0$. We assume that the function $G$ is C.A.I. in the set $D\subset \C_s\times\C_x$ so that
$$ \pi(D) \ = \  \{ x\in {\cal Y}\ : \ \Im \sigma_2(x)-\sigma_1(x)<A \} ,$$ 
and so that for $x\in \pi(D)$ the fiber $D_x=\pi^{-1}(x)\cap D$ is described as follows: (see fig.\ref{Paper3p21})
   \begin{itemize}
   \item  
for $x \in {\cal Y} \backslash ( {\cal U}_{\ell,\delta/2} \cup {\cal U}_{\ell'',\delta/2} ) $,
$$ D_x = \left\{ s\in \C \ :  \begin{array}{c}  \Im (\sigma_2(x))-\delta<\Im s <\Im (\sigma_2(x))+B \end{array} \right\} \backslash (\sigma_2(x)-i\R_{\ge 0}); $$ 
   \item for $x\in {\cal Y}\cap {\cal U}_{\ell,\delta/2}$,
$$ D_x = \left\{ s\in \C \ :  \begin{array}{c}  \Im (\sigma_2(x))-\delta<\Im s <\Im (\sigma_2(x))+B \\ 
(\Re s > \Re \sigma_2(x)) \Rightarrow (\Im s > \Im \sigma_1(x)) \end{array} \right\} \backslash (\sigma_2(x)-i\R_{\ge 0}); $$  
%
   \item for $x\in {\cal Y}\cap {\cal U}_{\ell'',\delta/2}$, 
$$ D_x = \left\{ s\in \C \ :  \begin{array}{c}  \Im (\sigma_1(x))<\Im s <\Im (\sigma_2(x))+B \end{array} \right\} \backslash (\sigma_2(x)-i\R_{\ge 0}) .$$ 
       \end{itemize}
If  $R_jG$ is defined and C.A.I. on $D\cap \pi^{-1}({\cal U}^\circ_{\ell,\delta/2})$, then $R_j G$ is also defined and C.A.I. on the whole $D$. 
\end{Lemma}

\private{\textcolor{blue}{Is $D$ open? Can I achieve the same in other cases?}}

\textsc{Proof.} Let us construct the function $R_jG(s,\xu)$ where $(s,\xu)\in D$ and $\xu\in \pi(D)\backslash {\cal U}_{\ell,\delta}$. For $\varepsilon \in \R$, let $D^{\varepsilon}_\xu = D_\xu \cap (D_\xu + i\varepsilon)$,  then $D_\xu = \bigcup_{\varepsilon>0} D^{\varepsilon}_\xu$.

If $\xu\in {\cal Y} \backslash ( {\cal U}_{\ell,\delta/2} \cup {\cal U}_{\ell'',\delta/2} )$, then 
 it is enough to construct, for each small enough  $\varepsilon>0$,  an integration path $y(t)$ starting at a point in ${\cal Y}\cap {\cal U}^\circ_{\ell,\delta/2}$ and ending at $\xu$ such that  $D^{\varepsilon}_\xu$ can be transported along $y(t)$ parallel to $-S_j$ in such a way that the set $D^{\varepsilon}_\xu+S_j(\xu)-S_j(y(t))$ will remain inside $D_{y(t)}$ for all $t$. Any path $y(t)$ within $\pi(D) \backslash ( {\cal U}_{\ell,(\delta-\varepsilon)/2} \cup {\cal U}_{\ell'',\delta/2} )$ from any point in ${\cal U}^\circ_{\ell,\delta/2}\backslash{\cal U}_{\ell,(\delta-\varepsilon)/2}$ to $\xu$ can be chosen as an integration path.

If $\xu\in {\cal Y} \cap {\cal U}_{\ell'',\delta/2}$, let $D^\varepsilon_\xu = D^{\varepsilon}_\xu \cap (D_\xu+i\varepsilon) \backslash Sl^\cap_{\varepsilon}(\sigma_1(\xu))$. 
Consider the set $P=\pi(D) \backslash \xi(D^\varepsilon_\xu)$  
where $\xi(s)$ is the branch of the  function $\sigma_1^{-1}\left( \frac{s+\sigma_1(\xu)}{2} \right)$ such that $\xi(\sigma_1(\xu))=\xu$ (cf. \eqref{SjInverse}). Then any set within $P$  from any point in ${\cal U}^\circ_{\ell,\delta'/2}\backslash{\cal U}_{\ell,\delta''/2}$ to $\xu$ can be chosen as an integration path. 


Continuity and analyticity of the function $R_jG$ can be checked as on page \pageref{WhyAnal}.
$\Box$

\begin{figure} \includegraphics{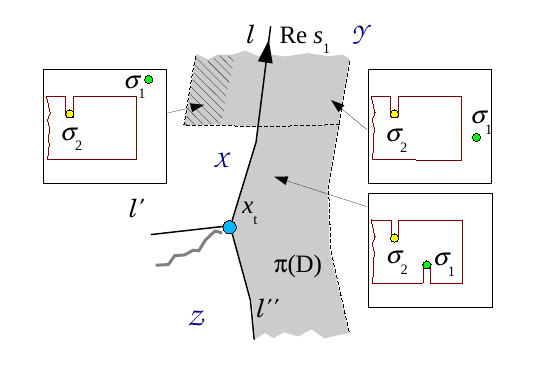} \caption{Lemma \ref{Le23}. The set $\pi(D)$ is shown, as well as $D_x$ for different $x$. The set $\pi(D)\cap ({\cal X}\backslash {\cal U}_{\ell,\delta/2})$ is hatched.} \label{Paper3p35} \end{figure}


\begin{Lemma} \label{Le23} 
Suppose $\sigma_2$ is a stationary and $\sigma_1$ is a moving singularity and $\sigma_2(x_t)=\sigma_1(x_t)$, and $\Re \sigma_1$ grows along $\ell$ in the direction away from $x_t$. Assume that there are also moving singularities $\sigma'_1,...,\sigma'_k$ ($k\ge 0$)  such that $\Im \sigma_1-\sigma'_j>\delta$ for all $j$.
Let $2c.w.({\cal Y})\ge A>\delta$ and $N>0$. 
We assume that the function $G$ is C.A.I. in a set $D\subset \C_s\times\C_x$ such that: 
$$\pi(D) = \{ x\in {\cal Y} \ : \  \Im \sigma_2(x)-\sigma_1(x) < A \} \ \cup \ \ \ \ \ \ \ $$ $$ \ \cup \ \{ x\in {\cal U}_{\ell,\delta} \ : \  \Re \sigma_1(x)-\sigma_2(x)>N  \ \text{and for $\forall j$} \ \Re \sigma'_j(x)-\sigma_2(x)>N \} $$
and for $x\in \pi(D)$, the fiber $D_x=\pi^{-1}(x)\cap D$ is described as follows:  
   \begin{enumerate}
   \item \label{Abr} if $x\in {\cal Y}$ satisfies  $\Re \sigma_1-\sigma_2\le N$ or, for some $j$, $\Re \sigma'_j-\sigma_2\le N$, 
$$ D_x = \left\{ s\in \C\ : \ \begin{array}{l} \Re s < \Re \sigma_2(x) + N \\ 
                                                \max\{\Im (\sigma_1(x))-\delta,\Im (\sigma_2(x))-A\} < \Im s < \Im (\sigma_2(x))+\delta \end{array} \right\} \backslash $$ 
$$ \ \ \ \ \backslash ((\sigma_2(x)+i\R_{\ge 0})\cup (\sigma_1(x)-i\R_{\ge 0}) ); $$
   \item if $x$ is such that  $\Re \sigma_1-\sigma_2> N$ and  $\Re \sigma'_j-\sigma_2>N$ for all $j$ 
$$ D_x = \left\{ s\in \C\ : \ \begin{array}{l} \Re s < \Re \sigma_2(x) + N \\ 
                                                \Im (\sigma_2(x))-A < \Im s < \Im (\sigma_2(x))+\delta \end{array} \right\} \backslash  (\sigma_2(x)+i\R_{\ge 0}) ; $$
     \end{enumerate}
If $R_jG$ is defined and C.A.I. on $D\cap \pi^{-1}({\cal X}\backslash {\cal U}_{\ell,\delta/2})$, then $R_j G$ is defined and C.A.I. on the whole $D$. \end{Lemma}

\textsc{Proof} is done by the same method as for other lemmas. $\Box$

We will also need the following modification of the lemma:

{\bf Lemma \ref{Le23}.A.} {\it  The same statement as in Lemma \ref{Le23}, with the item \ref{Abr} replaced by the following:
\begin{enumerate}
\item[(\ref{Abr}')] if $x\in {\cal Y}$ satisfies  $\Re \sigma_1-\sigma_2\le N$ or, for some $j$, $\Re \sigma'_j-\sigma_2\le N$,
$$ D_x^{\text{Lemma \ref{Le23}.A}} = D_x^{\text{Lemma \ref{Le23}.A}} \backslash  (\sigma_1(x)+(\R_{\ge 0}\times (-i)\R_{\ge 0}) ). $$
\end{enumerate}
}

\begin{figure} \includegraphics{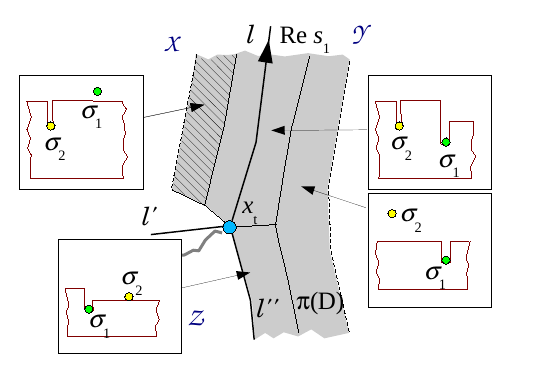} \caption{Lemma \ref{Lm1bis}. The set $\pi(D)$ is shown, as well as $D_x$ for different $x$. The set $\pi(D)\cap ({\cal X}\backslash {\cal U}_{\ell,\delta/2})$ is hatched.}  \label{Paper3p36} \end{figure}


\begin{Lemma} \label{Lm1bis} 
Suppose $\sigma_2$ is a stationary and $\sigma_1$ is a moving singularity and $\sigma_2(x_t)=\sigma_1(x_t)$, and $\Re \sigma_1$ grows along $\ell$ in the direction away from $x_t$.
 Assume $2c.w.({\cal Y})\ge A>\delta$, $+\infty\ge B>0$.  We assume that the function $G$ is C.A.I. in a set $D\subset \C_s\times\C_x$ such that 
$$ \pi(D) = \{ x\in {\cal Y} \ : \ \Im \sigma_2(x)-\sigma_1(x)<A \} \ \cup \ ( {\cal X} \cap {\cal U}_{\ell,\delta/2}) , $$
and for $x\in \pi(D)$, the fiber $D_x=\pi^{-1}(x)\cap D$ is described as follows: 
  \begin{enumerate} 
  \item If $x\in {\cal Y}\backslash ({\cal U}_{\ell,\delta/2}\cup {\cal U}_{\ell'',\delta/2})$, 
$$ D_x = \left\{ s\in \C : \begin{array}{l} \Re s < N + \Re \sigma_2(x) \\
                                            \Im (\sigma_2(x))- B < \Im s < \Im (\sigma_1(x))+\delta \end{array} \right\} \backslash (\sigma_1(x)+i\R_{\ge 0}); $$
  \item \label{Lm1bisN2} If $x\in ({\cal X}\cap {\cal U}_{\ell,\delta}) \cup  {\cal U}_{\ell,\delta/2}$, 
$$ D_x = \left\{ s\in \C : \begin{array}{l} \Re s < N + \Re \sigma_2(x) \\
                                            \Im (\sigma_2(x))- B < \Im s < \Im (\sigma_2(x))+\delta \\
                                            (\Re s > \Re \sigma_1(x))\Rightarrow (\Im s < \Im (\sigma_1(x))+\delta) \end{array} \right\} \backslash $$ $$ \ \ \  \ \backslash ((\sigma_1(x)+i\R_{\ge 0})\cup(\sigma_2(x)+i\R_{\ge 0})); $$
  \item If $x \in {\cal Y}\cap {\cal U}_{\ell'',\delta/2}$,
$$ D_x = \left\{ s\in \C : \begin{array}{l} \Re s < N + \Re \sigma_2(x) \\
                                            \Im (\sigma_2(x))- B < \Im s < \Im (\sigma_1(x))+\delta \\
                                            (\Re s> \Re \sigma_1(x))\Rightarrow (\Im s<\Im \sigma_2(x)) \end{array} \right\} \backslash (\sigma_1(x)+i\R_{\ge 0}). $$
  \end{enumerate} 
If  $R_jG$ is defined and C.A.I. on $D\cap \pi^{-1}( {\cal X}\backslash {\cal U}_{\ell,\delta/2})$, then $R_j G$ is defined and C.A.I. on all of $D$.
 \end{Lemma}

\textsc{Proof} is done by the same method as in other lemmas. \private{Care needs to be taken not to draw the integrations path for $\xu\not \in  {\cal U}_{\ell'',\delta/2}$  through ${\cal U}_{\ell,\delta/2}\cap {\cal U}_{\ell'',\delta/2}$ because of a small flap size there.}  $\Box$

\paragraph*{} There are a couple of easy modification of this lemma:

\paragraph*{} {\bf Lemma \ref{Lm1bis}.A.} {\it  Same statement as in lemma \ref{Lm1bis}, but with $D_x$ modifies as follows: $D_x^{\text{Lemma \ref{Lm1bis}.A} } = D_x^{\text{Lemma \ref{Lm1bis}} } \backslash (\sigma_3(x)-i\R_{\ge 0}) $ for some stationary singularity $\sigma_3(x)$. }

\paragraph*{} {\bf Lemma \ref{Lm1bis}.B.}  {\it 
Suppose $\sigma_1$ is a moving singularity,  $\Re \sigma_1$ grows along $\ell$ in the direction away from $x_t$, and let $\sigma_2(x)=2\sigma_1(x_t)-\sigma_1(x)$. \footnote{$\sigma_2(x)$ is not thought of as a first sheet singularity here}. The rest of the statement is the same as in Lemma \ref{Lm1bis}, except the item (\ref{Lm1bisN2}) is replaced by:
\begin{enumerate}
\item[(\ref{Lm1bisN2}')] 
   If $x\in ({\cal X}\cap {\cal U}_{\ell,\delta}) \cup  {\cal U}_{\ell,\delta/2}$,
$$ D_x = \left\{ s\in \C : \begin{array}{l} \Re s < N + \Re \sigma_2(x) \\
                                            \Im (\sigma_2(x))- B < \Im s < \Im (\sigma_2(x))+\delta \\
                                            (\Re s > \Re \sigma_1(x))\Rightarrow (\Im s < \Im (\sigma_1(x))+\delta) \end{array} \right\} \backslash (\sigma_1(x)+i\R_{\ge 0}) $$
(i.e. unlike in lemma \ref{Lm1bis}, $(\sigma_2(x)+i\R_{\ge 0})$ is not subtracted). 
\end{enumerate}
  }



\begin{figure} \includegraphics{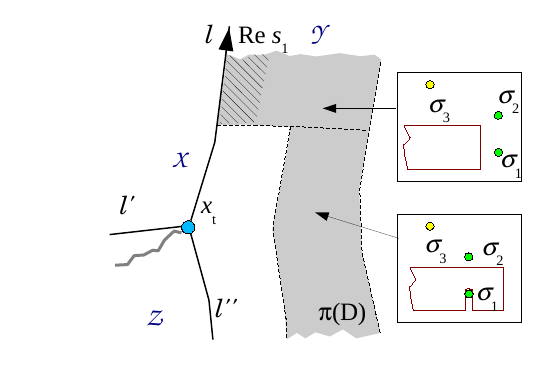} \caption{Lemma \ref{C3bo}. The set $\pi(D)$ is shown, as well as $D_x$ for different $x$. The set $\pi(D)\cap {\cal U}_{\ell,\delta/2}$ is hatched.} \label{Paper3p39} \end{figure}


\begin{Lemma} \label{C3bo}  
Suppose $\sigma_1$, $\sigma_2$ are moving and $\sigma_3$ is a stationary singularity, $\sigma_2(x_t)=\sigma_3(x_t)$, $\Im \sigma_2(x)>\Im \sigma_1(x)$, and $\Re \sigma_1$ grows along $\ell$ in the direction away from $x_t$. 
Assume $2c.w.({\cal Y})\ge B>A>0$, $N>0$.
We assume that the function $G$ is C.A.I. in a set $D\subset \C_s\times\C_x$ with the projection to $\C_x$: 
$$ \pi(D) = \{ x\in{\cal Y}  \ : \  A<\Im \sigma_3-\sigma_1<B \ \text{and} \ \Im \sigma_3-\sigma_2<A \} \cup $$
$$ \ \ \ \ \cup \{ x\in {\cal Y} \ : \ \Im \sigma_3-\sigma_1<B \ \text{and} \ \Re \sigma_1-\sigma_3>N \}, $$
and for $x\in \pi(D)$, the fiber $D_x=\pi^{-1}(x)\cap D$ is described as follows:
$$ D_x = \left\{ s\in \C :  \begin{array}{l} \Re s < \Re \sigma_3(x)+N \\
                            \Im(\sigma_3(x))-B <\Im s < \Im (\sigma_3(x))-A \\
                            (\Re \sigma_1(x)<\Re \sigma_3(x)+N) \Rightarrow \\ \ \ \Rightarrow ( \Im s > \max\{\Im (\sigma_1(x))-\delta, \Im (\sigma_1(x)-\sigma_2(x)+\sigma_3(x))- A\}) \end{array} \right\}\backslash $$ $$ \ \ \ \ \ \backslash (\sigma_1(x)-i\R_{\ge 0}) $$
If  $R_jG$ is defined and C.A.I. on ${\cal D}\cap\pi^{-1}({\cal U}_{\ell,\delta/2})$, then it is also defined and C.A.I. on the whole $D$.  \end{Lemma}

\textsc{Proof } by the same method as for other lemmas. $\Box$ 

\begin{figure} \includegraphics{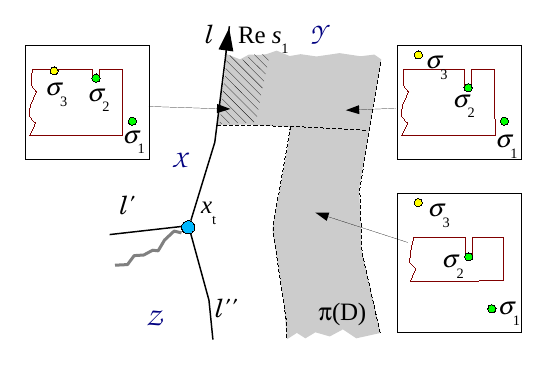} \caption{
Lemma \ref{C3T}. The set $\pi(D)$ is shown, as well as $D_x$ for different $x$. The set $\pi(D)\cap {\cal U}_{\ell,\delta/2}$ is hatched.} \label{Paper3p37} \end{figure}


\begin{Lemma} \label{C3T} 
Suppose $\sigma_1$, $\sigma_2$ are moving and $\sigma_3$ is a stationary singularity, $\sigma_2(x_t)=\sigma_3(x_t)$, $\Im \sigma_1(x)<\Im \sigma_2(x)$, and $\Re \sigma_1$ grows along $\ell$ in the direction away from $x_t$.
Assume $2c.w.({\cal Y})\ge B>A\ge 0$, $N>0$.  
We assume that the function $G$ is C.A.I. in a set $D\subset \C_s\times\C_x$ with the projection to $\C_x$: 
$$ \pi(D) = \{ x\in {\cal Y} \ : \ \Im(\sigma_3(x)-\sigma_1(x))>B \ \text{and} \  A<\Im(\sigma_3(x)-\sigma_2(x))<B \}  \ \cup \ $$
$$ \ \ \cup \  \{ x\in {\cal Y} \ : \ \Im(\sigma_3(x)-\sigma_2(x))<B \ \text{and} \ \Re(\sigma_3(x)-\sigma_1(x))>N  \} $$
and for $x\in \pi(D)$, the fiber $D_x=\pi^{-1}(x)\cap D$ is described as follows:
$$ D_\xu = \left\{ s\in \C \ : \ \begin{array}{l} \Re s < \Re \sigma_3(x)+ N; \\ 
                                                  \Im s < u(x) ; \\
                                                  \Im s > \Im (\sigma_3(x))-B   \end{array} \right\} \backslash (\sigma_2(x)+i\R), $$
where $u(x)=\min\{ \Im (\sigma_3(x))-A, \Im(\sigma_2(x)+\sigma_3(x)-\sigma_1(x))-B, \Im (\sigma_2(x))+\delta\}$.
If $R_jG$ is defined and C.A.I. on $D\cap \pi^{-1}({\cal U}_{\ell,\delta/2})$, then $R_jG$ is defined and C.A.I. on the whole $D$.  \end{Lemma} 

\begin{figure}[h] \includegraphics{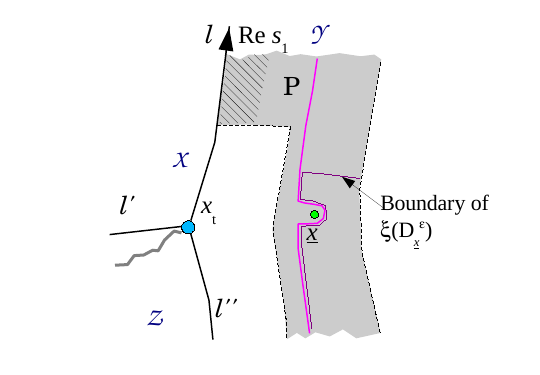}  \caption{Proof of lemma \ref{C3T} } \label{Paper3p38} \end{figure}

\textsc{Proof.} Let $(s,\xu)\in D$; we will limit ourselves to considering case of $\Re(\sigma_3(x)-\sigma_1(x))\le N$. 
Consider 
$$ D^\varepsilon_\xu = D_\xu \cap (D_\xu - i\varepsilon) \backslash (Sl^\cup_{\varepsilon}(\sigma_2(\xu)) ; $$
clearly $D_\xu=\bigcup_{\varepsilon>0} D_\xu^\varepsilon$. Thus it is enough to construct, for each sufficiently small $\varepsilon>0$, an integration path $y(t)$ starting at a point in $\pi(D)\cap {\cal U}_{\ell,\delta/2}$ and ending at $\xu$ such that the $D^{\varepsilon}_\xu$ can be transported along $y(t)$ parallel to $-S_j$ in such a way that set $D^{\varepsilon}_\xu+S_j(\xu)-S_j(y(t))$ will remain inside $D_{y(t)}$ for all $t$.

Consider the set $$P=\pi(D) \backslash \xi\left( \left\{ s\in \C \ : \ \begin{array}{l}  
                                                  \Im s < u(\xu) -\varepsilon  \end{array} \right\} \backslash Sl^\cup_{\varepsilon}(\sigma_2(\xu)) \right),$$ where $\xi(s)$ is the branch of the  function  $\sigma_2^{-1}\left( \frac{s+\sigma_2(\xu)}{2} \right)$ such that $\xi(\sigma_2(\xu))=\xu$ (cf. \eqref{SjInverse}).
Any path in $P$ from $\pi(D)\cap {\cal U}_{\ell,\delta/2}$ to $\xu$ can be taken as $y(t)$. Continuity and analyticity of the function $R_jG$ can be checked as on page \pageref{WhyAnal}. $\Box$


\begin{figure} \includegraphics{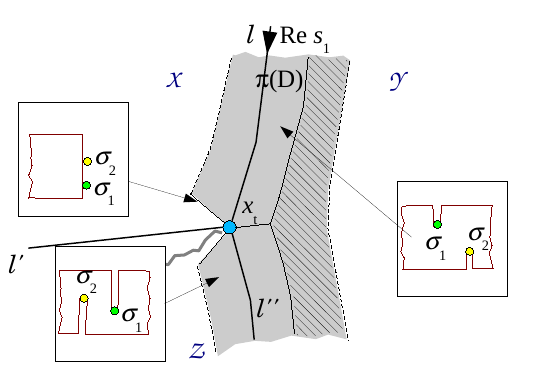} \caption{Lemma \ref{LN5thin}. The set $\pi(D)$ is shown, as well as $D_x$ for different $x$. The set $\pi(D) \cap (\{x\in{\cal Y} \ : \ \delta<\Im (\sigma_1(x)-\sigma_2(x))\})$ is indicated by hatching.} \label{Paper3p30} \end{figure}


\begin{Lemma} \label{LN5thin} 
Suppose $\sigma_1$ is a moving and $\sigma_2$ is a stationary singularity and $\sigma_1(x_t)=\sigma_2(x_t)$, and $\Re \sigma_1$ decreases along $\ell$ in the direction away from $x_t$. 
 We assume that the function $G$ is C.A.I. in a set $D\subset \C_s\times\C_x$ with 
$$ \pi(D) = {\cal U}_{\ell,\delta/2}  \ \cup \ {\cal U}_{\ell'',\delta/2} \ \cup \ \{x\in{\cal Y} \ : \ \Im (\sigma_1(x)-\sigma_2(x))<2\delta \}, $$
 so that for $x\in \pi(D)$ the fiber $D_x = D\cap \pi^{-1}(x)$ is described as follows:
\begin{itemize}
\item if $x\in {\cal X}\cup {\cal Z}$ satisfies $\Re \sigma_1(x) = \Re \sigma_2(x)$, then 
$$ D_x = \{ s\in\C \ : \ \Re s < \Re \sigma_1(x), \ \Im (\sigma_2(x))-\delta < \Im s < \Im (\sigma_1(x))+\delta\}; $$
\item otherwise
$$ D_x = \{ s\in\C \ : \  \Im (\sigma_2(x))-\delta < \Im s < \Im (\sigma_1(x))+\delta\} \backslash [(\sigma_2(x)-i\R_{\ge 0})\cup (\sigma_1(x)+i\R_{\ge 0})]. $$
\end{itemize}
If $R_jG$ is defined and C.A.I. on  $D\cap \pi^{-1}(\{x\in{\cal Y} \ : \ \delta<\Im (\sigma_1(x)-\sigma_2(x))\})$, 
then $R_jG$ is defined and C.A.I. on all of $D$.  \end{Lemma}

\textsc{Proof} proceeds by the same method as in the other lemmas. $\Box$

Also the following variant of this lemma will be used:

\paragraph*{} {\bf Lemma \ref{LN5thin}.A } {\it
Suppose $\sigma_1$ is a moving and $\sigma_2$ is a stationary singularity and $\sigma_1(x_t)=\sigma_2(x_t)$, and $\Re \sigma_1$ decreases along $\ell$ in the direction away from $x_t$. 
 We assume that the function $G$ is C.A.I. in a set $D\subset \C_s\times\C_x$ with 
$$ \pi(D) = {\cal U}_{\ell,\delta/2}  \ \cup \ \{x\in{\cal Y} \ : \ \Im (\sigma_1(x)-\sigma_2(x))<2\delta \}, $$
 so that for $x\in \pi(D)$ the fiber $D_x = D\cap \pi^{-1}(x)$ is described as follows:
\begin{itemize}
\item if $x\in {\cal X}$ satisfies $\Re \sigma_1(x) = \Re \sigma_2(x)$, then 
$$ D_x = \{ s\in\C \ : \ \Re s < \Re \sigma_1(x), \ \Im (\sigma_2(x))-\delta < \Im s < \Im (\sigma_1(x))+\delta\}; $$
\item otherwise
$$ D_x = \{ s\in\C \ : \  \Im (\sigma_2(x))-\delta < \Im s < \Im (\sigma_1(x))+\delta\} \backslash \ \ \ \ $$
$$\ \ \ \ \backslash [(\sigma_2(x)+[\R_{\ge 0}\times(-i\R_{\ge 0})])\cup (\sigma_1(x)+i\R_{\ge 0})]. $$
\end{itemize}
If $R_jG$ is defined and C.A.I. on $D\cap \pi^{-1}(\{x\in{\cal Y} \ : \ \delta<\Im (\sigma_1(x)-\sigma_2(x))\})$, 
then $R_jG$ is defined and C.A.I. on all of $D$.  }



\begin{figure} 
\includegraphics{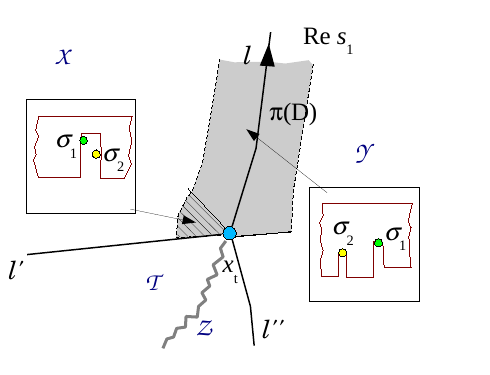}
\caption{Lemma \ref{LC10}: the set $\pi(D)$, as well as $\pi^{-1}(\xu)\cap D$ for $x\in {\cal U}_{\ell,\delta/2}$ and for $x\in \pi(D)\cap {\cal U}^\circ_{\ell',\delta/2}$. The set $\pi(D)\cap {\cal U}^\circ_{\ell',\delta/2}$  is hatched. } \label{Paper3p24} \end{figure}


\begin{Lemma} \label{LC10} 
Suppose $\sigma_1$ is a moving and $\sigma_2$ is a stationary singularity, $\sigma_1(x_t)=\sigma_2(x_t)$, and $\Re \sigma_1$ increases along $\ell$ in the direction away from $x_t$. 
Let $B>\delta$. We assume that the function $G$ is C.A.I. in a set $D\subset \C_s\times\C_x$ so that:
$$ \pi(D) \ = \  {\cal U}_{\ell,\delta/2} \  \cup \ \{ x\in {\cal X} : |\sigma_1(x)-\sigma_1(x_t)|<\delta/2 \} , $$
and for $x\in \pi(D)$, the fiber $D_x=\pi^{-1}(x)\cap D$ is described as follows: 
  \begin{itemize} 
    \item if $x\in \pi(D)\backslash {\cal U}^\circ_{\ell',\delta/2}$,
$$ D_x = \left\{ s\in \C : \begin{array}{l} \Im (\sigma_2(x)) - \delta< \Im s< \Im (\sigma_2(x))+B \\
                                            (\Re s> \sigma_1(x))\Rightarrow (\Im s >\Im (\sigma_1(x))-\delta ) \end{array} \right\} \backslash $$ $$ \ \ \ \ \ \backslash ((\sigma_1(x)-i\R_{\ge 0})\cup (\sigma_2(x)-i\R_{\ge 0})); $$
    \item if $x\in \{ x\in {\cal X}: |\sigma_1(x)-\sigma_1(x_t)|<\delta/2\}\backslash {\cal U}^\circ_{\ell,\delta/2}$,
if $x\in \pi(D)\backslash {\cal U}^\circ_{\ell',\delta/2}$,
$$ D_x = \left\{ s\in \C : \begin{array}{l} \Im (\sigma_1(x)) - \delta< \Im s< \Im (\sigma_2(x))+B \\
                                            (\Re s \in [\Re \sigma_1(x), \Re \sigma_2(x)])\Rightarrow (\Im s >\Im (\sigma_1(x)) ) \end{array} \right\}. $$ 
     \end{itemize} 
If  $R_jG$ is defined and C.A.I. on $D\cap \pi^{-1} ( {\cal U}^\circ_{\ell',\delta/2})$, then  $R_j G$ is defined and C.A.I. on the whole $D$. \end{Lemma}

\textsc{Proof} is very similar to the proof of Lemma \ref{Lemma19}. $\Box$



\begin{figure} \includegraphics{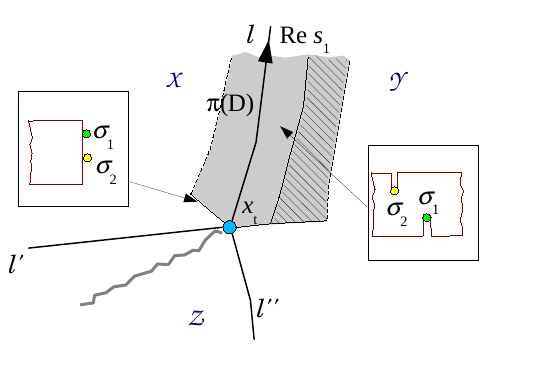} \caption{Lemma \ref{LN23thin}. The set $\pi(D)$ is shown, as well as $D_x$ for different $x$. The set $\pi(D) \cap ({\cal Y}\backslash {\cal U}_{\ell,\delta/2})$ is indicated by hatching.} \label{Paper3p32} \end{figure}


\begin{Lemma} \label{LN23thin} 
Suppose $\sigma_2$ is a stationary and $\sigma_1$ is a moving singularity and $\sigma_2(x_t)=\sigma_1(x_t)$, and $\Re \sigma_1$ grows along $\ell$ in the direction away from $x_t$.
We assume that the function $G$ is C.A.I. in a set $D\subset \C_s\times\C_x$ with the projection to $\C_x$: 
$$ \pi(D) = {\cal U}_{\ell,\delta/2} \cup ({\cal Y}\cap {\cal U}_{\ell,\delta} )$$
\begin{itemize}
\item if $x\in {\cal X}$ satisfies $\Re \sigma_1(x) = \Re \sigma_2(x)$, then 
$$ D_x = \{ s\in\C \ : \ \Re s < \Re \sigma_1(x), \ \Im (\sigma_1(x))-\delta < \Im s < \Im (\sigma_2(x))+\delta\}; $$
\item otherwise
\begin{equation} D_x = \{ s\in\C \ : \  \Im (\sigma_1(x))-\delta < \Im s < \Im (\sigma_2(x))+\delta\} \backslash [(\sigma_1(x)-i\R_{\ge 0})\cup (\sigma_2(x)+i\R_{\ge 0})]. \label{eq23thin06} \end{equation}
\end{itemize}
If $R_jG$ is defined and C.A.I. on  $D\cap \pi^{-1}({\cal Y}\backslash {\cal U}_{\ell,\delta/2})$, 
then $R_jG$ is defined and C.A.I. on all of $D$.  \end{Lemma}

\textsc{Proof} is done similarly to other lemmas in this section. $\Box$ \\

There are variants of this lemma:

\paragraph*{} {\bf Lemma \ref{LN23thin}.A.} {\it  The same statement as Lemma \ref{LN23thin}, but with \eqref{eq23thin06} replaced with
$$ D_x = \{ s\in\C \ : \  \Im (\sigma_1(x))-\delta < \Im s < \Im (\sigma_2(x))+\delta\} \backslash \ \ \ \ $$ $$ \ \ \ \backslash [(\sigma_1(x)+ [\R_{\ge 0}\times(-i\R_{\ge 0})])\cup (\sigma_2(x)+i\R_{\ge 0})] .$$ }

\begin{figure} \includegraphics{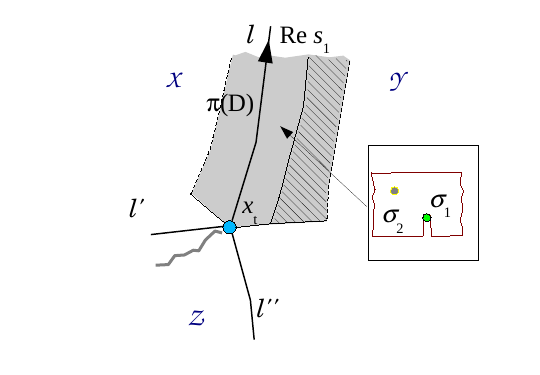} \caption{Lemma \ref{LN23thin}.B. The set $\pi(D)$ is shown, as well as $D_x$ for different $x$. The set $\pi(D) \cap ({\cal Y}\backslash {\cal U}_{\ell,\delta/2})$ is indicated by hatching.} \label{Paper3p33} \end{figure}

\paragraph*{} {\bf Lemma \ref{LN23thin}.B} {\it
Suppose $\sigma_1$ is a moving singularity; let $\sigma_2(x)=2\sigma_1(x_t)-\sigma_1(x)$, \footnote{here $\sigma_2(x)$ is not thought of as a first-sheet singularity} and $\Re \sigma_1$ grows along $\ell$ in the direction away from $x_t$.
We assume that the function $G$ is C.A.I. in a set $D\subset \C_s\times\C_x$ with the projection to $\C_x$: 
$$ \pi(D) = {\cal U}_{\ell,\delta/2} \cup ({\cal Y}\cap {\cal U}_{\ell,\delta}, $$
and 
$$ D_x = \{ s\in\C \ : \  \Im (\sigma_1(x))-\delta < \Im s < \Im (\sigma_2(x))+\delta\} \backslash (\sigma_1(x)-i\R_{\ge 0}). $$
If $R_jG$ is defined and C.A.I. on  $D\cap \pi^{-1}({\cal Y}\backslash {\cal U}_{\ell,\delta/2})$, 
then $R_jG$ is defined and C.A.I. on all of $D$. }


\begin{figure} \includegraphics{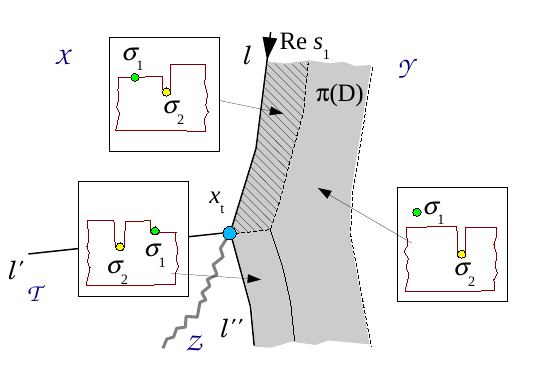} \caption{Lemma \ref{Lemma104}. The set $\pi(D)$ is shown, as well as $D_x$ for different $x$. The set $\pi(D)\cap {\cal U}^\circ_{\ell,\delta/2}$ is hatched.} \label{Paper3p22} \end{figure}


\begin{Lemma} \label{Lemma104} 
Suppose $\sigma_1$ is a moving and  $\sigma_2$ is a  stationary singularity and $\sigma_1(x_t)=\sigma_2(x_t)$, and $\Re \sigma_1$ decreases along $\ell$ in the direction away from $x_t$.  
Let $A,B>\delta$,  $2c.w.({\cal Y})\ge B$. We assume that the function $G$ is C.A.I. in a set $D\subset \C_s\times\C_x$ such that:
$$ \pi(D) =  \{ x\in {\cal Y} \ : \ \Im \sigma_1(x)-\sigma_2(x)<B \} , $$
and for $x\in \pi(D)$, the fiber $\pi^{-1}(x)\cap D$ is described as follows: 
  \begin{itemize} 
  \item if $x\in {\cal Y}\backslash {\cal U}^\circ_{\ell'',\delta/2}$,
$$ D_x = \left\{ s\in \C \ :  \begin{array}{c} \Im (\sigma_2(x))-A<\Im s <\Im (\sigma_2(x))+\delta \\
                                               (\Re s<\Re \sigma_2(x)) \Rightarrow (\Im s < \Im \sigma_1(x)) \end{array} \right\} \backslash (\sigma_2+i\R_{\ge 0}). $$
  \item for $x\in {\cal Y}\cap {\cal U}^\circ_{\ell'',\delta/2} $,
$$ D_x = \left\{ s\in \C \ :  \begin{array}{c} \Im (\sigma_2(x))-A<\Im s <\Im (\sigma_2(x))+\delta \\ 
(\Re s \ge \Re \sigma_1(x) ) \Rightarrow (\Im s < \Im \sigma_1(x))  \end{array} \right\} \backslash (\sigma_2+i\R_{\ge 0}).  $$
  \end{itemize} 
If $R_jG$ is defined and C.A.I. on $D\cap \pi^{-1}({\cal U}^\circ_{\ell,\delta/2})$ , then  $R_j G$ is defined and C.A.I. on the whole $D$. \end{Lemma}

\textsc{Proof.}  Let us construct the function $R_jG(s,\xu)$ where $(s,\xu)\in D$ and $\xu\in \pi(D)\backslash {\cal U}_{\ell,\delta/2}$. Denote $D_{x}=D\cap \pi^{-1}(x)$ and identify it with a subset of $\C$. For sufficiently small $\varepsilon>0$, let $D^{\varepsilon}_\xu = D_\xu \cap (D_\xu-i\varepsilon)\backslash Sl^\cup_\varepsilon(\sigma_1(\xu))$; then $D_\xu = \bigcup_{\varepsilon>0 } D^{\varepsilon}_\xu$.

If $\xu\in {\cal Y} \backslash ( {\cal U}_{\ell,\delta/2} \cup {\cal U}_{\ell'',\delta/2} )$, then, rephrasing the assumptions of the lemma,
$$ D_\xu = \left\{ s\in \C \ :  \begin{array}{c} \Im (\sigma_2(\xu))-A<\Im s <\Im (\sigma_2(\xu))+\delta \end{array} \right\} \backslash (\sigma_2+i\R_{\ge 0}). $$
 It is enough to construct, for each fixed $\varepsilon\in (0,\delta)$, an integration path $y(t)$ starting at a point in ${\cal Y}\cap {\cal U}^\circ_{\ell,\delta/2}$ and ending at $\xu$ such that the $D^{\varepsilon}_\xu$ can be transported along $y(t)$ parallel to $-S_j$ in such a way that set $D^{\varepsilon}_\xu+S_j(\xu)-S_j(y(t))$ will remain inside $D_{y(t)}$ for all $t$.

\begin{figure} \includegraphics{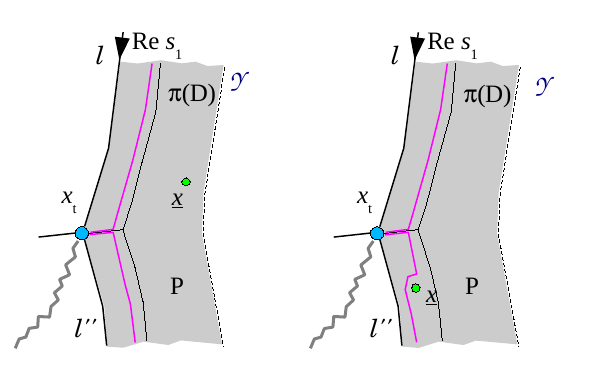}  \caption{Proof of Lemma \ref{Lemma104}. Left: for $ \xu\in {\cal Y} \backslash ( {\cal U}_{\ell,\delta/2} \cup {\cal U}_{\ell'',\delta/2} )$; right: for $ \xu\in {\cal Y} \cap  {\cal U}_{\ell'',\delta/2} $. The set $P$ is bounded by the purple curve and contains $\xu$. } \label{Paper3p23} \end{figure}

Consider the set $P=\pi(D) \backslash \xi(D^\varepsilon_\xu)$  
where $\xi(s)$ is the branch of the  function $\sigma_1^{-1}\left( \frac{s+\sigma_1(\xu)}{2} \right)$ such that $\xi(\sigma_1(\xu))=\xu$ (cf. \eqref{SjInverse}).

 Any path in $P$ from $\pi(D)\cap {\cal U}^\circ_{\ell,\delta/2}$ to $\xu$ can be taken as $y(t)$, fig.\ref{Paper3p23}.Continuity and analyticity of the function $R_jG$ can be checked as on page \pageref{WhyAnal}.

If $\xu\in {\cal Y} \cap {\cal U}_{\ell'',\delta/2}$, 
proceed analogously. 
 $\Box$





\begin{figure} 
\includegraphics{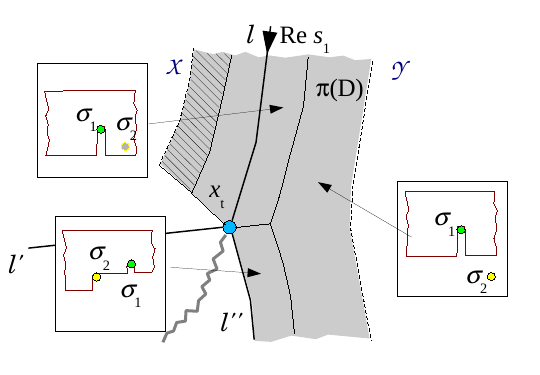} \caption{Lemma \ref{LC632}. The set $\pi(D)$ is shown, as well as $D_x$ for different $x$. The set $\pi(D)\cap ({\cal X}\backslash{\cal U}_{\ell,\delta/2})$ is hatched.} \label{Paper3p25} 
\end{figure}


\begin{Lemma} \label{LC632} 
Suppose $\sigma_1$ is a moving and $\sigma_2$ is a stationary singularities, $\sigma_1(x_t)=\sigma_2(x_t)$, and $\Re \sigma_1$ decreases along $\ell$ in the direction away from $x_t$. Let $+\infty\ge C\ge A>\delta$, $2c.w.({\cal Y})\ge A$.  
We assume that the function $G$ is C.A.I. in a set $D\subset \C_s\times\C_x$ with the projection to $\C_x$: 
$$ \pi(D) = {\cal U}_{\ell,\delta} \ \cup \ \{ x\in {\cal Y} \ : \ \Im \sigma_1(x)-\sigma_2(x) <A\},  $$ 
and for $x\in \pi(D)$, the fiber $D_x=\pi^{-1}(x)\cap D$ is described as follows:: 
   \begin{itemize}
   \item for  $x\in ({\cal Y}\backslash {\cal U}^\circ_{\ell',\delta/2})\cup {\cal U}_{\ell,\delta}$,
$$  D_x = \left\{ s\in \C \ : \ \begin{array}{l} \Im (\sigma_1(x))+\delta <\Im s < \Im (\sigma_2(x))+C    \end{array} \right\} \backslash (\sigma_1(x)+i\R). $$
   \item for $x\in {\cal U}^\circ_{\ell',\delta/2}$, 
$$ D_x = \left\{ s\in \C \ : \ \begin{array}{l} \Im (\sigma_1(x))+\delta <\Im s < \Im (\sigma_2(x))+C  \\
(\Re s\ge\Re \sigma_2(x)) \Rightarrow (\Im (s-\sigma_2(x))>0)   \end{array} \right\} \backslash (\sigma_1(x)+i\R).
 $$
     \end{itemize}
If $R_jG$ is defined and C.A.I. on $D\cap \pi^{-1}({\cal X}\backslash {\cal U}_{\ell,\delta/2})$, then  $R_j G$ is also defined and C.A.I. on the whole $D$. \end{Lemma}

\textsc{Proof}
Let $(s,\xu)\in D$; we will limit ourselves to considering the least trivial case of $\xu\in {\cal U}^\circ_{\ell',\delta/2}$. Consider  
$$ D^\varepsilon_\xu = D_\xu \cup (D_\xu + i\varepsilon) \backslash (Sl^\cup_{\varepsilon}(\sigma_1(\xu)) \cup Sl^\cup_{\varepsilon}(\sigma_2(\xu))); $$
clearly $D_\xu=\bigcup_{\varepsilon>0} D_\xu^\varepsilon$. Thus it is enough to construct, for each sufficiently small $\varepsilon>0$, an integration path $y(t)$ starting at a point in $\pi(D)\cap ({\cal X}\backslash {\cal U}^\circ_{\ell,\delta/2})$ and ending at $\xu$ such that the $D^{\varepsilon}_\xu$ can be transported along $y(t)$ parallel to $-S_j$ in such a way that set $D^{\varepsilon}_\xu+S_j(\xu)-S_j(y(t))$ will remain inside $D_{y(t)}$ for all $t$. 

\begin{figure} \includegraphics{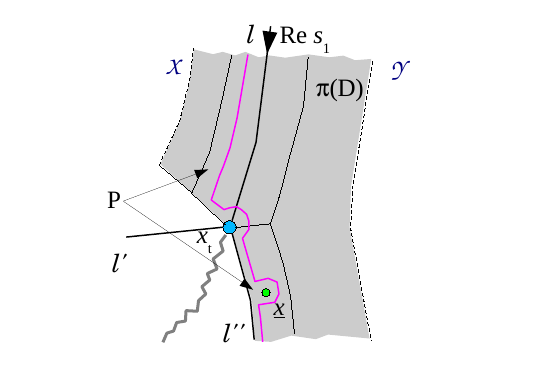}  \caption{Proof of lemma \ref{LC632} } \label{Paper3p26} \end{figure}

Consider the set $P=\pi(D) \backslash \xi(D^\varepsilon_\xu)$ where $\xi(s)$ is the branch of the  function  $\sigma_1^{-1}\left( \frac{s+\sigma_1(\xu)}{2} \right)$ such that $\xi(\sigma_1(\xu))=\xu$ (cf. \eqref{SjInverse}). Any path in $P$ from $\pi(D)\cap ({\cal X}\backslash{\cal U}^\circ_{\ell,\delta/2})$ to $\xu$ can be taken as $y(t)$. $\Box$

\begin{figure} 
\includegraphics{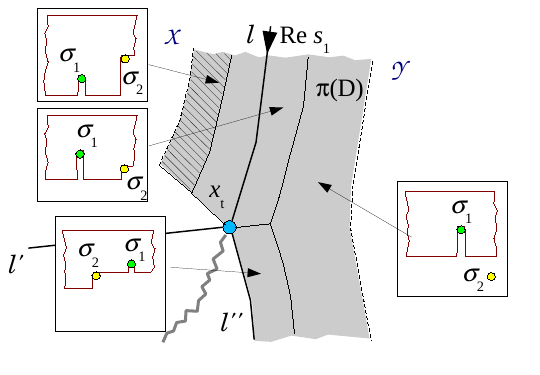} \caption{Lemma \ref{LC632}.A. The set $\pi(D)$ is shown, as well as $D_x$ for different $x$. The set $\pi(D)\cap ({\cal X}\backslash{\cal U}_{\ell,\delta/2})$ is hatched.} \label{Paper3p27} 
\end{figure}

\paragraph{} Here is a variant of this lemma:


\paragraph*{} {\bf Lemma \ref{LC632}.A.} 
{\it Suppose $\sigma_1$ is a moving and $\sigma_2$ is a stationary singularities, $\sigma_1(x_t)=\sigma_2(x_t)$, and $\Re \sigma_1$ decreases along $\ell$ in the direction away from $x_t$. The rest of the statement is the same as in the lemma \ref{LC632} with the only difference that:
$$ D_x^{\text{Lemma \ref{LC632}.A}} = D_x^{\text{Lemma \ref{LC632}}}\backslash (s_2+(\R_{\ge 0}\times (-i)\R_{\ge 0})). $$
}

\begin{figure} \includegraphics{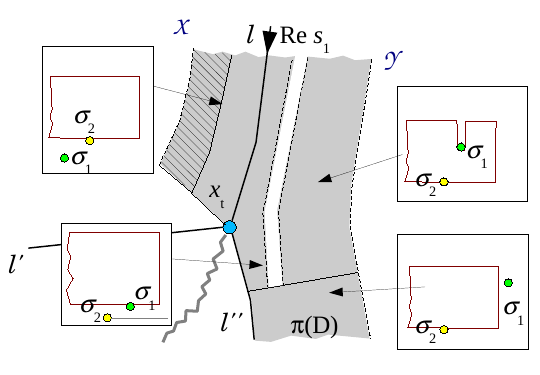} \caption{Lemma \ref{L5bis}. The set $\pi(D)$ is shown, as well as $D_x$ for different $x$. The set $\pi(D)\cap ({\cal X}\backslash {\cal U}_{\ell,\delta/2})$ is hatched.} \label{Paper3p28} \end{figure}





\begin{Lemma} \label{L5bis} 
Suppose $\sigma_1$ is a moving and $\sigma_2$ is a stationary singularity and $\sigma_1(x_t)=\sigma_2(x_t)$, and $\Re \sigma_1$ decreases along $\ell$ in the direction away from $x_t$. 
Let  $c.w.({\cal Y})\ge \frac{A}{2}>\delta>2\eta>0$, $N\in \R$. 
We assume that the function $G$ is C.A.I. in a set $D\subset \C_s\times\C_x$ with the projection to $\C_x$: 
$$ \pi(D) = (\overline{\cal X}\cap {\cal U}_{\ell,\delta} ) \ \cup \   \{ x \in {\cal Y}: \Im \sigma_1(x)-\sigma_2(x) < A \}  \ \backslash $$
$$ \ \ \ \ \backslash \{ x\in {\cal Y}: \ \eta<\Im (\sigma_1(x)-\sigma_2(x))<\delta-\eta ; \ \Re (\sigma_1(x)-\sigma_2(x))<N \} $$ 
and for $x\in \pi(D)$, the fiber $D_x=\pi^{-1}(x)\cap D$ is described as follows: 
\begin{itemize}
\item if $x\in {\cal Y}$ and $\delta-\eta \le \Im (\sigma_1(x)-\sigma_2(x)) < A$ and $\Re (\sigma_1(x)-\sigma_2(x)) < N$, 
$$ D_x = \left\{ s\in \C \ : \ \begin{array}{l} \Re s < \Re (\sigma_2(x)) + N; 
\\ \Im \sigma_2(x)< \Im s < \min\{ \Im(\sigma_1(x))+\delta, \Im (\sigma_2(x))+A\}  \end{array} \right\} \backslash (\sigma_1(\xu)+i\R_{\ge 0});
 $$
\item if $x\in {\cal Y}$ and $\Re (\sigma_1(x)-\sigma_2(x)) \ge  N$, 
$$ D_x = \left\{ s\in \C \ : \ \begin{array}{l} \Re s < \Re (\sigma_2(x)) + N; 
\\ \Im \sigma_2(x)< \Im s <  \Im (\sigma_2(x))+A   \end{array} \right\} ;
 $$
\item otherwise,
$$ D_x = \left\{ s\in \C \ : \ \begin{array}{l} \Re s < \Re (\sigma_2(x)) + N; 
\\ \max\{ \Im \sigma_1(x),\Im \sigma_2(x)\}< \Im s <  \Im (\sigma_2(x))+A   \end{array} \right\} .
 $$
\end{itemize}
If  $R_jG$ is defined and C.A.I. on ${\cal D}\cap \pi^{-1}({\cal X}\backslash {\cal U}_{\ell,\delta/2})$, then $R_jG$ is defined and C.A.I. on all of $D$.  \end{Lemma}

\begin{figure}[h] \includegraphics{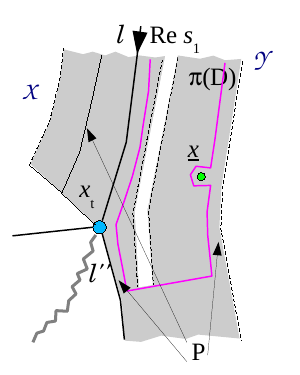}  \caption{Proof of lemma \ref{L5bis} } \label{Paper3p29} 
\end{figure}

\textsc{Proof} 
Let $(s,\xu)\in D$; we will limit ourselves to considering the least trivial case of $\xu\in {\cal Y}$ and $\delta -\eta < \Im (\sigma_1(\xu)-\sigma_2(\xu)) < A$ and $\Re (\sigma_1(\xu)-\sigma_2(\xu)) \le N$. 
Consider, for sufficiently small $\varepsilon>0$, 
$$ D^\varepsilon_\xu =  \left\{ s\in \C \ : \ \begin{array}{l} \Re s < \Re (\sigma_2(\xu)) + N; 
\\ \Im (\sigma_2(\xu))+\varepsilon< \Im s < \min\{ \Im(\sigma_1(\xu))+\delta, \Im (\sigma_2(\xu))+A\} -\varepsilon  \end{array} \right\} \backslash Sl^\cup_\varepsilon(\sigma_1(\xu)); $$ 
clearly $D_\xu=\bigcup_{\varepsilon>0} D_\xu^\varepsilon$. Thus it is enough to construct, for each sufficiently small $\varepsilon>0$, an integration path $y(t)$ starting at a point in $\pi(D)\cap ({\cal X}\backslash {\cal U}_{\ell,\delta/2})$ and ending at $\xu$ such that the $D^{\varepsilon}_\xu$ can be transported along $y(t)$ parallel to $-S_j$ in such a way that set $D^{\varepsilon}_\xu+S_j(\xu)-S_j(y(t))$ will remain inside $D_{y(t)}$ for all $t$.

Consider the set $P=\pi(D) \backslash \xi(D^\varepsilon_\xu)$ where $\xi(s)$ is the branch of the  function  $\sigma_1^{-1}\left( \frac{s+\sigma_1(\xu)}{2} \right)$ such that $\xi(\sigma_1(\xu))=\xu$ (cf. \eqref{SjInverse}). Any path in $P$ from $\pi(D)\cap ({\cal X}\backslash {\cal U}_{\ell,\delta/2})$ to $\xu$ can be taken as $y(t)$, fig.\ref{Paper3p29}. Analyticity of the resulting function can be checked as on page \ref{WhyAnal}. $\Box$

There is also the following variant of this lemma:



\begin{figure} \includegraphics{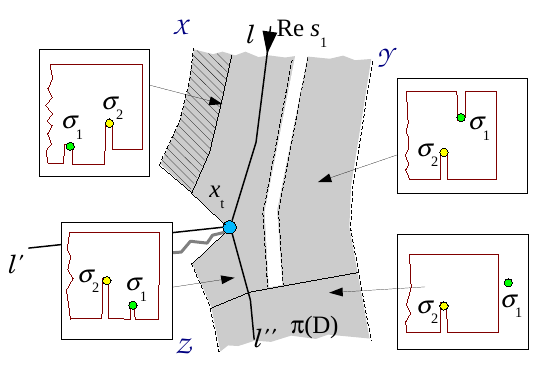} \caption{Lemma \ref{L5bis}.A. The set $\pi(D)$ is shown, as well as $D_x$ for different $x$. The set $\pi(D)\cap ({\cal X}\backslash {\cal U}_{\ell,\delta/2})$ is hatched. } \label{Paper3p31} \end{figure}




\paragraph*{} {\bf Lemma \ref{L5bis}.A.} {\it 
Suppose $\sigma_1$ is a moving and $\sigma_2$ is a stationary singularity and $\sigma_1(x_t)=\sigma_2(x_t)$, and $\Re \sigma_1$ decreases along $\ell$ in the direction away from $x_t$. 
Let  $c.w.({\cal Y})\ge \frac{A}{2}>\delta>2\eta>0$, $N\in \R$. 
We assume that the function $G$ is C.A.I. in a set $D\subset \C_s\times\C_x$ with the projection to $\C_x$: 
$$ \pi(D) = (\overline{\cal X}\cap {\cal U}_{\ell,\delta} ) \ \cup \ (\overline{\cal Z}\cap {\cal U}_{\ell'',\delta/2} )  \ \cup \  \{ x \in {\cal Y}: \Im \sigma_1(x)-\sigma_2(x) < A \}  \ \backslash $$
$$ \ \ \ \ \backslash \{ x\in {\cal Y}: \ \eta<\Im (\sigma_1(x)-\sigma_2(x))<\delta-\eta ; \ \Re (\sigma_1(x)-\sigma_2(x))<N \} $$
and for $x\in \pi(D)$, the fiber $D_x=\pi^{-1}(x)\cap D$ is described as follows: 
\begin{itemize}
\item if $x\in {\cal Y}$ and $\delta \le \Im (\sigma_1(x)-\sigma_2(x)) < A$ and $\Re (\sigma_1(x)-\sigma_2(x)) < N$, 
$$ D_x = \left\{ s\in \C \ : \ \begin{array}{l} \Re s < \Re (\sigma_2(x)) + N; 
\\ \Im (\sigma_2(x))-\delta< \Im s < \min\{ \Im(\sigma_1(x))+\delta, \Im (\sigma_2(x))+A\}  \end{array} \right\} \backslash $$ $$ \ \ \ \ \backslash ( (\sigma_1(\xu)+i\R_{\ge 0}) \cup (\sigma_2(\xu)-i\R_{\ge 0}));
 $$
\item if $x\in {\cal Y}$ and $\Re (\sigma_1(x)-\sigma_2(x)) \ge  N$, 
$$ D_x = \left\{ s\in \C \ : \ \begin{array}{l} \Re s < \Re (\sigma_2(x)) + N; 
\\ \Im (\sigma_2(x))-\delta< \Im s <  \Im (\sigma_2(x))+A   \end{array} \right\} \backslash  (\sigma_2(\xu)-i\R_{\ge 0}) ;
 $$
\item if $x\in {\cal U}_{\ell'',\delta/2}$,
$$ D_x = \left\{ s\in \C \ : \ \begin{array}{l} \Re s < \Re (\sigma_2(x)) + N; 
\\ \Im (\sigma_2(x))-\delta< \Im s <  \Im (\sigma_2(x))+A   \end{array} \right\} \backslash $$
$$ \ \ \ \ \backslash ((\sigma_1(\xu)-i\R_{\ge 0})\cup (\sigma_2(\xu)-i\R_{\ge 0})) ; $$
\item otherwise,
$$ D_x = \left\{ s\in \C \ : \ \begin{array}{l} \Re s < \Re (\sigma_2(x)) + N; 
\\\Im \sigma_1(x)-\delta< \Im s <  \Im (\sigma_2(x))+A   \\
(\Re s > \Re \sigma_2(x)) \Rightarrow (\Im s > \Im (\sigma_2(x)) - \delta)\end{array} \right\}\backslash $$ $$ \ \ \ \ \ \backslash ((\sigma_1(\xu)-i\R_{\ge 0})\cup (\sigma_2(\xu)-i\R_{\ge 0})) .
 $$
\end{itemize}
If  $R_jG$ is defined and C.A.I. on ${\cal D}\cap \pi^{-1}({\cal X}\backslash {\cal U}_{\ell,\delta/2})$, then $R_jG$ is defined and C.A.I. on all of $D$.  }


\textsc{Proof} analogous to the proof of lemma \ref{L5bis}. $\Box$

\begin{figure} \includegraphics{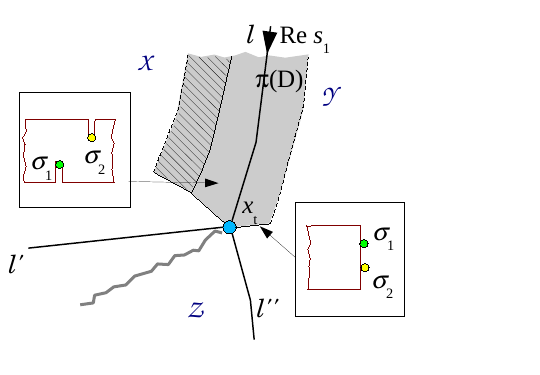} \caption{Lemma \ref{L63LB}. The set $\pi(D)$ is shown, as well as $D_x$ for different $x$. The set $\pi(D)\cap({\cal X}\backslash{\cal U}_{\ell,\delta/2})$ is hatched.} \label{Paper3p34} 
\end{figure}


\begin{Lemma} \label{L63LB}  
Suppose $\sigma_1$ is a moving and $\sigma_2$ is a stationary singularity and $\sigma_1(x_t)=\sigma_2(x_t)$, and $\Re \sigma_1$ decreases along $\ell$ in the direction away from $x_t$. 
We assume that the function $G$ is C.A.I. in a set $D\subset \C_s\times\C_x$ with the projection to $\C_x$: 
$$ \pi(D) = {\cal U}_{\ell,\delta/2} \cup ({\cal X}\cap {\cal U}_{\ell,\delta}), $$
and for $x\in \pi(D)$, the fiber $D_x=\pi^{-1}(x)\cap D$ is described as follows: 
\begin{itemize}
\item if $x\in {\cal Y}$ satisfies $\Re \sigma_1(x) = \Re \sigma_2(x)$, then 
$$ D_x = \{ s\in\C \ : \ \Re s < \Re \sigma_1(x), \ \Im (\sigma_1(x))-\delta < \Im s < \Im (\sigma_2(x))+\delta\}; $$
\item otherwise
$$ D_x = \{ s\in\C \ : \  \Im (\sigma_1(x))-\delta < \Im s < \Im (\sigma_2(x))+\delta\} \backslash [(\sigma_1(x)-i\R_{\ge 0})\cup (\sigma_2(x)+i\R_{\ge 0})]. $$
\end{itemize}
If $R_jG$ is defined and C.A.I. on  $D\cap \pi^{-1}({\cal X}\backslash{\cal U}_{\ell,\delta/2})$, 
then $R_jG$ is defined and C.A.I. on all of $D$.  \end{Lemma}

\textsc{Proof} by the same method as in the other lemmas. $\Box$

\begin{figure} 
\includegraphics{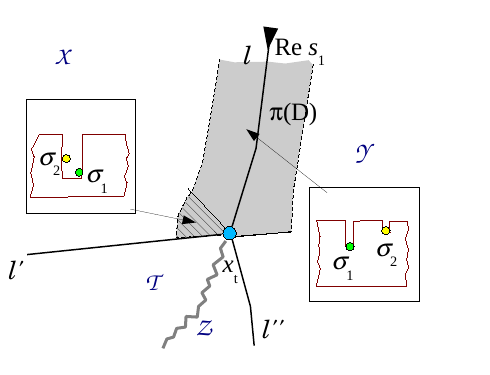}
\caption{Lemma \ref{Lemma19}. The set $\pi(D)$ is shown, as well as $\pi^{-1}(\xu)\cap D$ for $x\in {\cal U}_{\ell,\delta/2}$ and for $x\in \pi(D)\cap {\cal U}^\circ_{\ell',\delta/2}$. The set $\pi(D)\cap {\cal U}^\circ_{\ell',\delta/2}$ is hatched.} \label{Paper3p18} \end{figure}

\begin{Lemma} \label{Lemma19} 
Suppose $\sigma_1$ is a moving and $\sigma_2$ is a stationary singularity, $\sigma_1(x_t)=\sigma_2(x_t)$, and $\Re \sigma_1$ decreases along $\ell$ in the direction away from $x_t$. 
Let $B>\delta$. We assume that the function $G$ is C.A.I. in a set $D\subset \C_s\times\C_x$ so that:
$$ \pi(D) \ = \  {\cal U}_{\ell,\delta/2} \  \cup \ \{ x\in {\cal X} : |\sigma_1(x)-\sigma_1(x_t)|<\delta/2 \} , $$
and for $x\in \pi(D)$, the fiber $D_x=\pi^{-1}(x)\cap D$ is described as follows: 
  \begin{itemize} 
    \item if $x\in \pi(D)\backslash {\cal U}^\circ_{\ell',\delta/2}$, 
$$ D_x = \left\{ s\in \C \ : \ \begin{array}{l} \Im s > \Im \sigma_2(x)-B ; \\
\Im (s-\sigma_1(x))<\delta  \\
(\Re s>\Re \sigma_2(x)) \Rightarrow (\Im (s-\sigma_2(x))<\delta)   \end{array} \right\} \backslash ((\sigma_1(x)+i\R)\cup (\sigma_2(x)+i\R)).
 $$
    \item if $x\in \{ x\in {\cal X}: |\sigma_1(x)-\sigma_1(x_t)|<\delta/2\}\backslash {\cal U}^\circ_{\ell,\delta/2}$,
$$ D_x = \left\{ s\in \C \ : \ \begin{array}{l} \Im s > \Im \sigma_2(x)-B ; \\
\Im (s-\sigma_1(x))<\delta  \\
(\Re s \in [\Re \sigma_2(x), \Re \sigma_1(x)]) \Rightarrow (\Im s < \Im \sigma_1(x))   \end{array} \right\}.   $$
     \end{itemize} 
If  $R_jG$ is defined and C.A.I. on $D\cap \pi^{-1} ( {\cal U}^\circ_{\ell',\delta/2})$, then  $R_j G$ is defined and C.A.I. on the whole $D$. \end{Lemma}

\textsc{Proof.}  Let $(s,\xu)\in D$ and $\xu\in \pi(D)\backslash {\cal U}^\circ_{\ell',\delta/2}$. 
Consider 
$$ D^\varepsilon_\xu = D_\xu \cap (D_\xu - i\varepsilon) \backslash (Sl^\cup_{\varepsilon}(\sigma_1(\xu)) \cup Sl^\cup_{\varepsilon}(\sigma_2(\xu))); $$
clearly $D_\xu=\bigcup_{\varepsilon>0} D_\xu^\varepsilon$. Thus it is enough to construct, for each sufficiently small $\varepsilon>0$, an integration path $y(t)$ starting at a point in $\pi(D)\cap {\cal U}^\circ_{\ell',\delta/2}$ and ending at $\xu$ such that the $D^{\varepsilon}_\xu$ can be transported along $y(t)$ parallel to $-S_j$ in such a way that set $D^{\varepsilon}_\xu+S_j(\xu)-S_j(y(t))$ will remain inside $D_{y(t)}$ for all $t$. 

\begin{figure} 
 \includegraphics{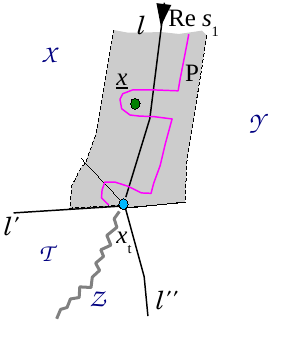} \caption{Proof of Lemma \ref{Lemma19}.} \label{Paper3p19} \end{figure}

Consider the set $P=\pi(D) \backslash \xi(D^\varepsilon_\xu)$ where $\xi(s)$ is the branch of the  function  $\sigma_1^{-1}\left( \frac{s+\sigma_1(\xu)}{2} \right)$ such that $\xi(\sigma_1(\xu))=\xu$ (cf. \eqref{SjInverse}). Any path in $P$ from $\pi(D)\cap {\cal U}^\circ_{\ell',\delta/2}$ to $\xu$ can be taken as $y(t)$. Continuity and analyticity of the function $R_jG$ can be checked as on page \pageref{WhyAnal}. $\Box$


\subsection{Applications of lemmas from section \ref{modelcases} ...} \label{ApplyManyLemmas}

We will now consider one after the other the two operators $R_j$, $j=1,2$, and for each of them -- sets $X\in {\mathbb S}\backslash \{ {\cal A}_{int}\}$. 

If $X={\cal Y}_{int}$ for some Stokes region ${\cal Y}$, then for each fixed $\xu\in{\cal Y}_{int}$ we will cover ${\cal S}_\xu$ by horizontal strips as in section \ref{FSS}. E.g., for $\xu\in {\cal C}_{int}$ these strips are: the strip above $S(\xu)$, the strips between $S(\xu)$ and $-s_2(\xu)$, between $-s_2(\xu)$ and $-s_1(\xu)$, between $-s_1(\xu)$ and $-S(\xu)$, and the strip below $-S(\xu)$. For each such strip $U\subset {\cal S}_\xu$ we will indicate a method of constructing the integration path from an earlier element of ${\mathbb S}$ to $\xu$ for $R_j$ and all $(s,\xu)\in U$. 

If $X={\cal U}_{L,\delta/2}$ for a Stokes curve $L$, we argued in the section \ref{StrategySec} that it remains to construct the integration path from an earlier element of ${\mathbb S}$ to $\xu\in {\cal U}_{L,\delta/2}$ for $(s,\xu)\in {\cal S}_{\xu}$ where $s$ belongs to charts $U_{\sigma_1\pm i0, \sigma_2\pm i0}\subset {\cal S}_{\xu}$ etc, where sigmas denote singularities of ${\cal S}_\xu$. We will indicate how this construction is done for each of these charts. 

\subsubsection{... towards constructing $R_1G$.} 

In the lemmas we are going to use for studying $R_1$, one should interpret blue singularities as stationary and red singularities as moving.

For $\xu$ along the curve $L_1$, i.e. for $\xu\in {\cal U}_{L_1,\delta/2}$:

\begin{itemize}

\item chart $(S+i0)$ -- Lemma \ref{LC632}.

\item chart $(S-i0,-s_1+i0)$ -- Lemma \ref{LN5thin} reduces the question to the existence of $R_1G(s,\xu)$ for $\xu\in {\cal B}_{int}$, $s$ in strip between $S(\xu)$ and $-s_1(\xu)$; below we will see how that in turn reduces to the existence of $R_1G(s,x)$ for $x\in {\cal A}_{int}$. 

\item chart $(S-i0,-s_1-i0)$ -- Lemma \ref{Lemma19} reduces the question to existence of $R_1G$ along the curve $L'_1$; considerations along the curve $L'_1$ reduce the question to existence on $R_1G$ in ${\cal A}_{int}$.

\end{itemize}

For $\xu$ in the region ${\cal B}_{int}$:

\begin{itemize}

\item above $S(x)$ -- Lemma \ref{LC632}.

\item between $S(\xu)$ and $-s_1(\xu)$ -- Lemma \ref{L5bis}.A. 
\private{ \textcolor{blue}{with $A=2\Im[S(x_2)-S(x_1)]$.}}



\item between $-s_1$ and $-S$, below $-S$ -- Lemma \ref{TrivialReason}.

\end{itemize}


For $\xu$ around the curve $L_2$, i.e. $\xu\in {\cal U}_{L_2,\delta/2}$:

\begin{itemize}

\item chart $(S+i0)$ -- Lemma \ref{LC632}.

\item chart $(S-i0,-s_2+i0)$ -- Lemma \ref{LN5thin} reduces the question to the existence of $R_1G(s,\xu)$ for $\xu\in {\cal C}_{int}$, $s$ in strip between $S(\xu)$ and $-s_1(\xu)$; below we will see how that in turn reduces to the existence of $R_1G(s,x)$ for $x\in {\cal B}_{int}$. 

 application of this lemma uses the existence of $R_1G$  in region ${\cal B}$, in the strip between $S$ and $-s_2$, which will be shown below.

\item chart $(S-i0,-s_2-i0)$ -- Lemma \ref{Lemma19} reduces the question to existence of $R_1G$ along the curve $L'_2$; considerations along the curve $L'_2$ reduce the question to existence on $R_1G$ in ${\cal B}_{int}$.
. 

\end{itemize}

For $\xu$ in the region ${\cal C}_{int}$:

\begin{itemize}

\item above $S(x)$ -- Lemma \ref{LC632}.

\item between $S$ and $-s_2$ -- Lemma \ref{L5bis}.A.


\item between $-s_2$ and $-s_1$, between $-s_1$ and $-S$ -- Lemma \ref{TrivialReason}.B.

\item below $-S$ -- Lemma \ref{TrivialReason}. 

\end{itemize}

For $\xu$ along the curve $L'_2$, i.e. $\xu\in {\cal U}_{L'_2,\delta/2}$:

\begin{itemize}

\item chart $(-s_2+i0,S-i0)$ -- Lemma \ref{LN5thin}.

\item chart $(-s_2+i0,S+i0)$ -- Lemma \ref{LC10}.

\item charts $(-s_1+i0)$, $(-S+i0)$ -- Lemma \ref{TrivialReason}.C.

\item charts $(-s_2-i0,S+\underline{i0})$, $(-s_1-i0, -s_{12'}+\underline{i0})$, $(-S-i0,-s_{2'}+\underline{i0})$ -- Lemma \ref{LN23thin}.A reduces the question to constructing the analytic continuation of $R_1G$ in the region ${\cal D}_{int}$; consideration in the region ${\cal D}_{int}$ described below will further reduce it to the situation in the region ${\cal C}_{int}$.

\item charts $(-s_2-i0,S-i0)$, $(-s_1-i0, -s_{12'}-i0)$, $(-S-i0,-s_{2'}-i0)$ -- The proof follows from Lemmas \ref{Lm1bis}.A and \ref{Lm1bis} and can be done together with the case of region ${\cal D}$, strips between $S$ and $-s_1$, between $-s_{12'}$ and $-s_{2'}$, and under $-s_2$, respectively.

\end{itemize}


For $\xu$ in the region ${\cal D}_{int}$: 

\begin{itemize}

\item above $-s_2$ -- by Lemma \ref{Lemma18}. 

\item between $-s_2$ and $S$, between $-s_1$ and $-s_{12'}$, between $-S$ and $-s_{2'}$ --  Lemma \ref{Le23}.A.

\item between $S$ and $-s_1$, between $-s_{12'}$ and $-S$, under $-s_{2'}$ -- Lemmas \ref{Lm1bis}.A and \ref{Lm1bis}.

\end{itemize}


For $\xu$ along the curve $L''_2$, i.e. $\xu\in {\cal U}_{L''_2,\delta/2}$:

\begin{itemize}

\item charts  $(-s_1+i0,-s_{12''}-\underline{i0})$, $(-S+i0,-s_{2''}-\underline{i0})$ --  Lemma \ref{LN23thin}.A reduces the question to constructing the analytic continuation of $R_1G$ in the region ${\cal E}_{int}$; consideration in the region ${\cal E}_{int}$ described below will further reduce it to the situation in the region ${\cal B}_{int}$. Lemma \ref{LN23thin}.B does the same for the chart $(S-i0)$. 

\item charts $(S+i0)$,  -- Lemma \ref{Lm1bis}.B.

\item chars $(-s_1+i0,-s_{12''}+i0)$, $(-S+i0,-s_{2''}+i0)$ -- Lemma \ref{Lm1bis} or \ref{Lm1bis}.A.

\item charts $(-s_{12'}-i0)$, $(-s_{2'}-i0)$ -- Lemma \ref{TrivialReason}.C. 

\end{itemize}

For $\xu$ in the region ${\cal E}_{int}$:

\begin{itemize}

\item above $S$ -- Lemma \ref{Lm1bis}.B.

\item between $S$ and $-s_{12''}$ -- cover this strip by two thinner horizontal strips: one from $\Im [S(\xu)]+\delta$ to $\Im [-s_{12''}(\xu)]+\delta$ is dealt with by Lemma \ref{C3bo}; the strip from $S(\xu)$ to the flap below the cut $[-s_{12''}(x),+\infty)$ is treated with Lemma \ref{C3T} if $\xu\not\in {\cal U}_{\tilde L'_2,\delta/2}$ and by Lemma \ref{Lm1bis} if $\xu\in {\cal U}_{\tilde L'_2,\delta/2}$. 

\item between $-s_{12''}$ and $-s_1$, between $-s_{2''}$ and $-S$ -- Lemma \ref{Le23}.A. 
Note that in order to apply the lemma \ref{Le23}.A, we need to know existence of $R_1G$ on a certain subset of $\pi^{-1}({\cal B}\cup {\cal U}_{L''_2,\delta/2})$;  existence of $R_1G$ there follows from the considerations in charts $ (-s_1+i0,-s_{12''}+i0)$, $(-S+i0, -s_2+i0)$ along $L''_2$.

\item between $-s_1$ and $-s_{2''}$ -- lemma \ref{Lm1bis}.A.

\item below $-S$ -- Lemma \ref{TrivialReason}. 
\private{\textcolor{red}{what is $\pi(D)$?}.}

\end{itemize}

For $\xu$ the curve $L''_1$, i.e. for $\xu\in {\cal U}_{L''_1,\delta/2}$:

\begin{itemize}

\item chart $(-s_1+i0,S-i0)$ -- Lemma \ref{LN5thin}.

\item chart $(-S+i0)$ -- Lemma \ref{TrivialReason}.C. 

\item chart $(-s_1+i0,S+i0)$ -- Lemma \ref{LC10}.

\item charts $(-s_1-i0,S+\underline{i0})$, $(-S-i0,-s_{1''}+\underline{i0})$ -- Lemma \ref{LN23thin}.A reduces the question to constructing the analytic continuation of $R_1G$ in the region ${\cal F}_{int}$; consideration in the region ${\cal F}_{int}$ described below will further reduce it to the situation in the region ${\cal B}_{int}$.

\item charts $(-s_1-i0,S-i0)$, $(-S-i0,-s_{1''}-i0)$ -- Lemma \ref{Lm1bis}.

\end{itemize}

For $\xu\in {\cal F}_{int}$:

\begin{itemize} 

\item above $-s_1$ -- Lemma \ref{Lemma18}.

\item between $-s_1$ and $S$, between $-S$ and $-s_{1''}$ -- Lemma \ref{Le23}.A.

\item between $S$ and $-S$: If $\xu\in {\cal F}$ is such that $\Im S(\xu)\le \delta$, draw an integration path along the curve $\{ t: \Re S(t)=\Re S(\xu)\}$ from $\xu$ to $\xu'\in {\cal F}$ where $\Re S(\xu')=\Re S(\xu)$ and $\Im S(\xu')=\frac{3}{2}\delta$.  

Thus, we can assume that $\Im S(\xu)>\delta$.  Now for the part of the strip from the flap along $(S,\infty)$ to $\Im (-S(\xu))+\delta$ use Lemma \ref{Lm1bis}, and for the part of the strip from $\Im S(\xu)-\eta$, for sufficiently small $\eta>0$, to $\Im (-S(\xu))-\delta$ use Lemma \ref{TrivialReason}. 
\private{ \\ \textcolor{blue}{Sort out these details}: with $\pi(D)=\{ x\in {\cal F} : \Im S(x)>\Im S(\xu)-\eta \} \ \cup {\cal U}_{L''_1,\delta, \textcolor{red}{\Box}}\backslash {\cal U}_{\tilde L'_1, \delta/2, \textcolor{red}{\Box}}$ \textcolor{blue}{\bf carefully deal with the slot size and with the closed part of $\partial U$.}}

\item below $-s_{1''}$ -- Lemma \ref{Lm1bis}.

\end{itemize}

For $\xu$ along the curve $L'_1$, i.e. for $\xu\in {\cal U}_{L'_1,\delta/2}$: 

\begin{itemize}

\item chart $(S+i0)$ -- Lemma \ref{TrivialReason}.C.

\item chart $(-S+i0,-s_{1'}-\underline{i0})$ -- Lemma \ref{LN23thin}.A reduces the question to constructing the analytic continuation of $R_1G$ in the region ${\cal G}_{int}$; consideration in the region ${\cal G}_{int}$ described below will further reduce it to the situation in the region ${\cal A}_{int}$. ; Lemma \ref{LN23thin}.B does the same for the chart $(S-i0)$.

\item chart $(-S+i0,-s_{1'}+i0)$ -- Lemma \ref{Lm1bis}.

\item chart $(S+i0)$ -- Lemma \ref{Lm1bis}.B.

\end{itemize}

For $\xu$ in the region ${\cal G}_{int}$:

\begin{itemize}

\item above $S$ -- Lemma \ref{Lm1bis}.B.

\item between $S$ and $-s_{1'}$ -- Cover this strip by two thinner horizontal strips: a strip from the flap above the cut $[S(\xu), +\infty)$ to $\Im [-s_{1'}(\xu)+\delta$ is dealt with by Lemma \ref{C3bo}; the strip from $S(\xu)$ to the flap below the cut $[-s_{1'}(x),+\infty)$ is treated with Lemma \ref{C3T} if $\xu\not\in {\cal U}_{\tilde L''_1,\delta/2}$ and by Lemma \ref{Lm1bis} if $\xu\in {\cal U}_{\tilde L''_1,\delta/2}$.

\item between $-s_{1'}$ and $-S$ -- Lemma \ref{Le23}.A can be used to obtain the result for $\xu$ satisfying $\Im [S(x_1)-S(\xu)]>\delta/2$. If one chooses $0<\delta'<\delta$ and applies Lemma \ref{Le23}.A with $\delta'$ instead of $\delta$, $R_1G$ can be constructed on this strip for $\xu$ satisfying $\Im [S(x_1)-S(\xu)]>\delta'/2$.

\item under $-S$ -- Lemma \ref{TrivialReason}.
\private{ \textcolor{blue}{ \\  More explanation is needed: Lemma \ref{TrivialReason} will guarantee the analytic continuation to the region which for $\Re s< \Re [-S(x)]$ goes up only until $\Im s=\Im [-S(x)]$; some words need to be said why this does not create a gap along the horizontal line to left of $-S(x)$.} } 

\end{itemize}


\subsubsection{... towards constructing $R_2G$.}

In the lemmas we are going to use for studying $R_2$, one should interpret red singularities as stationary and blue singularities as moving.

For $\xu$ along the curve $L_1$, i.e. $\xu\in {\cal U}_{L_1,\delta/2}$:

\begin{itemize}

\item chart $(S+i0)$ -- Lemma \ref{TrivialReason}.C.

\item chart $(S-i0,-s_1+i0)$ -- Lemma \ref{LN23thin} reduces the question to constructing the analytic continuation of $R_2G$ in the region ${\cal B}_{int}$; consideration in the region ${\cal B}_{int}$ described below will further reduce it to the situation in the region ${\cal A}_{int}$; Lemma \ref{LN23thin}.B does the same for the chart $(-S+i0)$.

\item chart $(S-i0,-s_1-i0)$ -- Lemma \ref{Lm1bis}.

\item chart  $(-S-i0)$ -- Lemma \ref{Lm1bis}.

\end{itemize}

For $\xu$ in the region ${\cal B}_{int}$:

\begin{itemize}

\item above $S$ -- Lemma \ref{TrivialReason}. \private{\textcolor{red}{what is $\pi(D)$?}}

\item between $S$ and $-s_1$ -- Lemma \ref{Le23}.

\item between $-s_1$ and $-S$ --  cover this strip by two thinner horizontal strips and apply Lemmas \ref{C3T} an \ref{C3bo}.

\item below $-S$ -- Lemma \ref{Lm1bis}.B.

\end{itemize}

For $\xu$ along the curve $L_2$, i.e. $\xu\in {\cal U}_{L_2,\delta/2}$: 

\begin{itemize}

\item chart $(S+i0)$ -- Lemma \ref{TrivialReason}.C.

\item chart $(S-i0,-s_2+i0)$ -- -- Lemma \ref{LN23thin} reduces the question to constructing the analytic continuation of $R_2G$ in the region ${\cal C}_{int}$; consideration in the region ${\cal C}_{int}$ described below will further reduce it to the situation in the region ${\cal B}_{int}$; Lemma \ref{LN23thin}.B does the same for the charts $(-s_1+i0)$ and $(-S+i0)$.

\item chart $(S-i0,-s_2-i0)$ -- Lemma \ref{Lm1bis}.

\item charts $(-s_1-i0)$ and $(-S-i0)$ -- Lemma \ref{Lm1bis}.B.

\end{itemize}

For $\xu$ in the region ${\cal C}_{int}$: 

\begin{itemize}

\item above $S$ -- Lemma \ref{TrivialReason}. 
\private{ \\ \textcolor{blue}{ with $\pi(D)={\cal C} \cup {\cal U}_{L_2,\delta, \textcolor{red}{\Box}}\backslash {\cal U}_{ L'_2, \delta/2, \textcolor{red}{\Box}}$ \textcolor{blue}{\bf carefully deal with the slot size and with the closed part of $\partial U$.}};}

\item between $S$ and $-s_2$ --  Lemma \ref{Le23}.A can be used to obtain the result for $\xu$ satisfying $ M-\Im[2S(\xu)-2S(x_2)]>\delta$. If one chooses $0<\delta'<\delta$ and applies Lemma \ref{Le23}.A with $\delta'$ instead of $\delta$, $R_2G$ can be constructed on this strip for $\xu$ satisfying $M-\Im [2S(\xu)-2S(x_2)]>\delta_1$.

\item between $-s_2$ and $-s_1$ , between $-s_1$ and $-S$ -- Cover each these strip by two thinner horizontal strips and apply Lemmas \ref{C3T} an \ref{C3bo}.

\item below $-S$ -- Lemma \ref{Lm1bis}.B.

\end{itemize}

For $\xu$ along the curve $L'_2$, i.e. $\xu\in {\cal U}_{L'_2,\delta/2}$:

\begin{itemize}

\item chart $(-s_2+i0,S-i0)$ -- Lemma \ref{L63LB}.

\item chart $(-s_2+i0,S+\underline{i0})$ -- Lemma \ref{LC632}.A.

\item charts $(-s_1+i0)$, $(-S+i0)$ -- Lemma \ref{LC632}.

\item charts $(-s_2-i0,S+\underline{i0})$, $(-s_1-i0,-s_{12'}+\underline{i0})$, $(-S-i0,-s_{2'}+\underline{i0})$ -- Lemma \ref{LN5thin}.A, using the construction in the region ${\cal D}$ to be performed below.

\item charts $(-s_2-i0,S-i0)$, $(-s_1-i0,-s_{12'}-i0)$, $(-S-i0,-s_{2'}-i0)$ -- Lemma \ref{Lemma19} reduces the question to existence of $R_2G$ along the curve $L_2$.

\end{itemize}

For $\xu$ in the region ${\cal D}_{int}$:

\begin{itemize}

\item above $-s_2$ -- Lemma \ref{LC632}.A.

\item between $-s_2$ and $S$, between $-s_1$ and $-s_{12}$, between $-S$ and $-s_{2'}$ -- Lemma \ref{L5bis}.

\item between $S(\xu)$ and $-s_1(\xu)$ -- The part of this strip from the flap above $S(\xu)$ to $\Im -s_1(\xu)$ is dealt with by Lemma \ref{Lemma104}, the part of the strip between from $\Im S(\xu)$ until the flap below $-s_1(\xu)$ is dealt with by Lemma \ref{LC632}. 

\item between $-s_{12'}(\xu)$ and $-S(\xu)$ -- The part of this strip from the flap above $-s_{12'}(\xu)$ to $\Im -S(\xu)$ is dealt with by Lemma \ref{Lemma104}, the part of the strip between from $\Im (-s_{12'}(\xu))$ until the flap below $-S(\xu)$ is dealt with by lemma \ref{LC632}. 

\item under $-s_{2'}$ -- Lemma \ref{Lemma104}.

\end{itemize}


For $\xu$ along the curve $L''_2$, i.e. $\xu\in {\cal U}_{L''_2,\delta/2}$:

\begin{itemize}

\item charts $(-s_1+i0,-s_{12''}-\underline{i0})$, $(-S+i0,-s_{2''}-\underline{i0})$  --Lemma \ref{LN5thin}.A reduces the question to construction of analytic continuation of $R_2G$ in the region ${\cal E}_{int}$; considerations for the region ${\cal E}$ reduce the question further to existence of $R_2G$ in ${\cal B}$.

\item charts $(-s_1+i0,-s_{12''}+i0)$, $(-S+i0,-s_{2''}+i0)$ -- Lemma \ref{Lemma19} reduces the question to existence of $R_2G$ along the curve $L_2$; considerations along the curve $L_2$ reduce the question to existence on $R_2G$ in ${\cal B}_{int}$. \footnote{Notice for comparison that when we are constructing $R_1G$, we reduce the question of existence of $R_1G$ along $L_2$ to existence of $R_1G$ along $L'_2$, and then -- to its existence in ${\cal B}_{int}$; here we have to proceed in the opposite order.}

\item charts $(-s_{1}-i0)$, $(-S-i0)$ -- Lemma \ref{LC632}.

\end{itemize}

For $\xu$ in the region ${\cal E}_{int}$:

\begin{itemize} 

\item above $S(x)$ -- Lemma \ref{TrivialReason} 
\private{\textcolor{red}{what is $\pi(D)$?};}

\item between $S$ and $-s_{12''}$ --  For the thinner strip from the flap along $(S,\infty)$ down to $\Im [-s_{12''}]+\delta$, use Lemma \ref{TrivialReason};
\private{ \textcolor{red}{what is $\pi(D)$?}} 
for the thinner strip from $\Im S-\delta$ down to flap along $(-s_{12''},\infty)$ -- Lemma \ref{Lemma104}.

\item between $-s_{12''}$ and $-s_1$, between $-s_{2''}$ and $-S$ -- Lemma \ref{L5bis}.

\item  between $-s_1$ and $-s_{2''}$ -- Without loss of generality, assume $\Im -s_1(\xu)-[-s_{2''}(\xu)]>2\delta$ (otherwise reduce the situation to this one by drawing a piece of the integration path such that $\Im S$ decreases and $\Re S$ stays constant along it). Then, for the part of the strip from the flap along $(-s_1,\infty)$ to $\Im (-s_{2''}(\xu))+\delta$ use Lemma \ref{LC632}, and for the part of the strip from $\Im -s_1(\xu)$ to the flap along $(-s_{2''}(\xu),\infty)$ use Lemma \ref{Lemma104}.

\item under $-S$ -- lemma \ref{LC632}. 

\end{itemize}

For $\xu$ along the curve $L''_1$, i.e. $\xu\in {\cal U}_{L''_1,\delta/2}$: 

\begin{itemize}

\item chart $(-s_1+i0,S-i0)$ -- Lemma \ref{L63LB}.

\item chart $(-s_1+i0,S+\underline{i0})$ -- Lemma \ref{LC632}.A.

\item charts $(-s_1-i0,S+\underline{i0})$, $(-S-i0,-s_{1'}+i0)$ --  Lemma \ref{LN5thin}.A reduces the question to construction of analytic continuation of $R_2G$ in the region ${\cal F}_{int}$; considerations for the region ${\cal F}$ reduce the question further to existence of $R_2G$ in ${\cal B}$.

\item chart $(-s_1-i0,S-i0)$, $(-S-i0,-s_{1'}-i0)$ -- Lemma \ref{Lemma19} reduces the question to existence of $R_2G$ along the curve $L_1$.

\item chart $(-S+i0)$ -- Lemma \ref{LC632}.

\end{itemize}

For $\xu$ in the region ${\cal F}_{int}$: 

\begin{itemize} 

\item above $-s_{1}$ -- Lemma \ref{LC632}.A.

\item  between $-s_{1}(x)$ and $-S(x)$, between $-S(\xu)$ and $-s_{1''}(\xu)$ -- Lemma \ref{L5bis}. 

\item between $S$ and $-S$ -- Without loss of generality, assume $\Im S(\xu)-[-S(\xu)]=2\Im S(\xu)>2\delta$ (otherwise reduce the situation to this one by drawing a piece of the integration path such that $\Im S$ increases and $\Re S$ stays constant along it). Now for the part of the strip from the flap along $(S,\infty)$ to $\Im (-S(\xu))+\delta$ use Lemma \ref{Lemma104}, and for the part of the strip from $\Im S(\xu)$ to $\Im -S(\xu)-\delta$ use Lemma \ref{LC632}.

\item below $-s_{1''}$ -- Lemma \ref{Lemma104}.

\end{itemize}

For $\xu$ along the curve $L'_1$, i.e. $\xu\in {\cal U}_{L'_1,\delta/2}$:

\begin{itemize}

\item chart $(-S-i0)$ -- Lemma \ref{LC632}.

\item chart $(-S+i0,-s_{1'}-\underline{i0})$ -- Lemma \ref{LN5thin}.A reduces the question to construction of analytic continuation of $R_2G$ in the region ${\cal G}_{int}$; considerations for the region ${\cal G}$ reduce the question further to existence of $R_2G$ in ${\cal A}$.

\item chart $(-S+i0,-s_{1'}+i0)$ -- Lemma \ref{Lemma19} reduces the question to existence of $R_2G$ along the curve $L_1$; considerations along the curve $L_1$ reduce the question to existence on $R_2G$ in ${\cal A}_{int}$.

\end{itemize}

For $\xu$ in the region ${\cal G}_{int}$:

\begin{itemize}

\item above $S$ -- Lemma \ref{TrivialReason}.A. 

\item between $S$ and $-s_{1'}$ -- For a thinner strip from the flap along $(S,\infty)$ to $\Im[-s_{1'}]+\delta$, use Lemma \ref{TrivialReason}.B, for a thinner strip from $\Im S -\delta$ to the flap along $(-s_{1'},\infty)$ -- Lemma \ref{Lemma104}.

\item between $-s_{1'}$ and $-S$ -- Lemma \ref{L5bis}.

\item below $-S$ -- Lemma \ref{LC632}.

\end{itemize}


{\it Remark.} In every case when $\xu \in {\cal U}_{\ell,\delta/2}\cap {\cal U}_{\ell'',\delta/2}$, where  $\ell$ and $\ell''$ are Stokes curves starting from the same turning point $x_t$, the analyticity of $R_jG$ at all points $(s,\xu)$, $s\in {\cal S}_\xu$, can be deduced from the combination of lemmas applied for $\xu$ in ${\cal U}_{\ell,\delta/2}$ and for $\xu$ in ${\cal U}_{\ell'',\delta/2}$. 

The above list provides a construction of $R_jG$ for every point of ${\cal S}$ and concludes the proof of theorem \ref{RjGmainTh}. $\Box$

\private{\textcolor{blue}{Need to write a good exposition on how to cover the situation in charts ${\cal U}_{\ell,\delta/2}\cap {\cal U}_{\ell'',\delta/2}$.}}


\section{Concluding remarks}

\paragraph*{Too many singularities for $\xu$ in ${\cal D}$ and ${\cal F}$. }

We see from the description of ${\mathcal S}$ that once $\xu$ goes one loop around the turning points $x_1$ or $x_2$, the locations of the singularities of the fiber ${\mathcal S}_\xu$ remains the same up to a permutation, except for one singularity in each case: namely, the singularity $S(\xu)$ is present in the regions ${\cal D}$ and ${\cal F}$, but the corresponding singularity is absent in the regions ${\cal E}$ and ${\cal G}$. 

Assuming that the series (\ref{vNseries}) converges and $\Phi(s,x)={\tilde Y}f(s)$ indeed gives a solution to the equation (\ref{MainEqu1}), we can hope to make a rigorous sense of the observation of ~\cite[page 243 and on]{V83} that the Laplace integral of $\Phi(s,x)$ gives a solution of (\ref{SchroeEq1}) which is unramified at $x_1$ and $x_2$, and hence show that $\Phi(s,x)$ has only a removable singularity at $S(\xu)$ for $\xu$ in  ${\cal D}$ and ${\cal E}$, and also show the relations between other singularities for $\xu$ in ${\cal D}$ and ${\cal E}$, ${\cal F}$  and ${\cal G}$ that would amount to Voros' connection formulas.

\paragraph*{Virtual turning points.}  [We admit that our treatment of virtual turning points in the earlier versions of this paper was incorrect.] If in section \ref{NotationSec} we decide to take a larger domain ${\cal O}\subset \C$, then we might have to discuss {\it virtual turning points}, e.g. a point $x_3$ on the boundary of Stokes region ${\cal D}$ satisfying $S(x_3)=S(x_1)$. At this point, the singularities $-s_1(x)$ and $S(x)$ coincide, but $V(x_3)$ does not have to vanish. According to ~\cite{V83}, at $x_3$ the two singularities are expected to pass through each other without creating any new singularities in ${\cal S}$. 

If one wishes to repeat the argument of this article in an example of domain ${\cal O}$  containing turning points $x_1,...,x_k$ (where $V$ has a zero)  and virtual turning points $x_{k+1},...,x_{k+\ell}$ (where different singularities of ${\cal S}$ coincide), we propose to:\\ 
a) take as $\tilde{\cal O}$ an appropriate subset the universal cover of ${\cal O}\backslash\{x_1,...x_{k+\ell} \}$; \\
b) define Stokes curves by a condition that $\Im \int_{x_t}^x p(y)dy=0$ where $x_t$ is either a turning point or a virtual turning point; these curves will split $\tilde{\cal O}$ into Stokes regions; \\ 
c) construct the fiber of ${\cal S}$ over each point of $\tilde{\cal O}$ bearing in mind that no new singularities should appear on the first sheet of ${\cal S}$ when we cross a Stokes curve starting at a virtual turning point. \\  Lemmas of section \ref{modelcases} work with minor modifications when $x_t$ is a virtual turning point.

{\bf \large Acknowledgments} 

This work was mainly carried out during the author's studies at the Department of Mathematics, Northwestern University, U.S.A., and during his stay at the Max Planck Institute for Mathematics in the Sciences, Leipzig, Germany. The author is profoundly grateful to Dmitry Tamarkin and Boris Tsygan for numerous discussions and for conscientious critique of the manuscript, and also to Ovidiu Costin, Stavros Garoufalidis, Rostislav Matveyev, Shinji Sasaki, Boris Shapiro, Yoshitsugu Takei, and Jared Wunsch for discussions, comments, and feedback. The author also appreciates valuable comments of the anonymous referee.



\vspace{3cm}


\end{document}